\documentclass[openany, amssymb, psamsfonts]{amsart}
\usepackage{amssymb,amsmath}
\usepackage{mathrsfs,comment}
\usepackage[usenames,dvipsnames]{color}
\usepackage[normalem]{ulem}
\usepackage{url}
\usepackage[all,arc,2cell]{xy}
\UseAllTwocells
\usepackage{enumerate}
\usepackage{tikz-cd}
\usepackage{hyperref}  
\hypersetup{%
	bookmarksnumbered=true,%
	bookmarks=true,%
	colorlinks=true,%
	linkcolor=blue,%
	citecolor=blue,%
	filecolor=blue,%
	menucolor=blue,%
	pagecolor=blue,%
	urlcolor=blue,%
	pdfnewwindow=true,%
	pdfstartview=FitBH}

\def\makeautorefname#1#2{\expandafter\def\csname#1autorefname\endcsname{#2}}

\def\equationautorefname~#1\null{(#1)\null}
\makeautorefname{footnote}{footnote}%
\makeautorefname{item}{item}%
\makeautorefname{figure}{Figure}%
\makeautorefname{table}{Table}%
\makeautorefname{part}{Part}%
\makeautorefname{appendix}{Appendix}%
\makeautorefname{chapter}{Chapter}%
\makeautorefname{section}{Section}%
\makeautorefname{subsection}{Section}%
\makeautorefname{subsubsection}{Section}%
\makeautorefname{theorem}{Theorem}%
\makeautorefname{thm}{Theorem}%
\makeautorefname{cor}{Corollary}%
\makeautorefname{lem}{Lemma}%
\makeautorefname{prop}{Proposition}%
\makeautorefname{pro}{Property}
\makeautorefname{conj}{Conjecture}%
\makeautorefname{defn}{Definition}%
\makeautorefname{notn}{Notation}
\makeautorefname{notns}{Notations}
\makeautorefname{rem}{Remark}%
\makeautorefname{quest}{Question}%
\makeautorefname{exmp}{Example}%
\makeautorefname{ax}{Axiom}%
\makeautorefname{claim}{Claim}%
\makeautorefname{ass}{Assumption}%
\makeautorefname{asss}{Assumptions}%
\makeautorefname{con}{Construction}%
\makeautorefname{prob}{Problem}%
\makeautorefname{warn}{Warning}%
\makeautorefname{obs}{Observation}%
\makeautorefname{conv}{Convention}%

\newtheorem{thm}{Theorem}[section]

\newtheorem{prop}{Proposition}[section]
\newtheorem{lem}{Lemma}[section]

\theoremstyle{definition}

\newtheorem{exmp}{Example}[section]

\newtheorem{rem}{Remark}[section]

\makeatletter
\let\c@obs=\c@thm
\let\c@cor=\c@thm
\let\c@prop=\c@thm
\let\c@lem=\c@thm
\let\c@prob=\c@thm
\let\c@con=\c@thm
\let\c@conj=\c@thm
\let\c@defn=\c@thm
\let\c@notn=\c@thm
\let\c@notns=\c@thm
\let\c@exmp=\c@thm
\let\c@ax=\c@thm
\let\c@pro=\c@thm
\let\c@ass=\c@thm
\let\c@warn=\c@thm
\let\c@rem=\c@thm
\let\c@sch=\c@thm
\let\c@equation\c@thm
\numberwithin{equation}{section}
\makeatother

\title{Generic doubling of rectangular pegs}
\author{Zhen Gao}
\date{\today}

\address{Institut für Mathematik, Universität Augsburg, Augsburg, Germany}
\email{zhen.gao@math.uni-augsburg.de}
\usepackage{lipsum}

\begin{document}	
	
\maketitle

\begin{center}
With an appendix joint with Urs Frauenfelder
\end{center}
	
\thispagestyle{empty}

\begin{abstract}
We prove a multiplicity result for rectangular pegs that there is a generic class of smooth Jordan curves in which every curve admits two geometrically distinct similar inscribed rectangles with aspect angle in $(0,\frac{\pi}{2})$, based on the existence of rectangular pegs in any smooth Jordan curve, which is first proved by Greene and Lobb \cite{Greene-Lobb2021} and we give an alternative Floer theoretical proof in this paper.  The key insight is that the rectangular peg problem is translated into finding intersection points of two Lagrangian tori. We present two distinct proofs for the multiplicity result: one involves Lagrangian Floer homology, and the other is differential topological in nature which employs a novel computation formula for the algebraic intersection number. Both rely crucially on certain generic geometric transversality of the two tori, and the correspondence between the intersection points and the inscribed rectangles. Moreover, such a generic doubling result can be further extended to cyclic quadrilateral pegs based on the existence result \cite{Greene-Lobb2023}. In Appendix \ref{Appendix_B}, written jointly with Urs Frauenfelder, we provide the detailed derivation of the computation formula for the algebraic intersection number.
\end{abstract}

\tableofcontents	

\section{Introduction}\label{Section 1}
The inscribed square problem a.k.a. square peg problem firstly posed by Otto Toeplitz in 1911, c.f. \cite{Toeplitz1911}, asks the following question:
\begin{center}
\textit{Does every continuous simple closed curve a.k.a. Jordan curve inscribe a square?}
\end{center}
More precisely, a square being inscribed on a Jordan curve means that there are four points in the Jordan curve forming four vertices of the square. The square peg problem has been affirmatively answered in various cases, yet remains open in its full generality to this day. We refer to \cite{Matschke2014} and \cite[\S 2]{Greene-Lobb2021} for historical survey and recent developments. There are many variants of the square peg problem, see e.g. \cite[p.349-351]{Matschke2014} and \cite{Matschke2022}. One of the closely related problems to the square peg problem is the so-called rectangular peg problem. We list some important breakthroughs: In 1977, Vaughan showed that every continuous
Jordan curve inscribes a rectangle (note that there is no knowlege about its aspect angle which is the angle between
the two diagonals); from 2018 to 2019, Hugelmeyer showed the existence of rectangular peg on any smooth Jordan curve with any aspect angle lying in a subset of $(0,\frac{\pi}{2}]$ whose Lebesgue measure is at least $\frac{\pi}{6}$, c.f. \cite{Hugelmeyer2021}; in 2020, the smooth rectangular peg problem in full generality, i.e. for any smooth Jordan curve and with any aspect angle in $(0,\frac{\pi}{2}]$, was solved by Greene and Lobb, c.f. \cite{Greene-Lobb2021}. In particular, the result of Greene and Lobb in \cite{Greene-Lobb2021} solves the square peg problem in smooth case. Further study by Greene and Lobb generalizes this result to cyclic quadrilateral pegs, c.f. \cite{Greene-Lobb2023}. The most important and interesting feature of their proofs is involving highly non-trivial symplectic topology. We review the geometric setup in their paper \cite{Greene-Lobb2021} below, and state the main result that we obtain based on their work.

Let $\gamma$ be a smooth simple closed curve in $\mathbb{C}$, i.e. a smooth Jordan curve. This is obviously a Lagrangian closed curve in $(\mathbb{C},\omega_\text{std})$ with $\omega_\text{std}:=\frac{i}{2}dz\wedge d\bar{z}$ in terms of coordinate $z\in\mathbb{C}$. The product $\gamma\times\gamma$ is then a smooth Lagrangian torus in $(\mathbb{C}^2,\omega_\text{std})$ where 
\[
\omega_\text{std}:=\frac{i}{2}(dz_1\wedge d\bar{z}_1 + dz_2 \wedge d\bar{z}_2)
\] 
in terms of coordinate $(z_1,z_2)\in\mathbb{C}^2$. Now we consider the map 
\begin{align}\label{l_map}
l:\mathbb{C}^2 &\to \mathbb{C}^2 \\
		(z_1,z_2) &   \mapsto \left(\frac{z_1+z_2}{2} , \frac{z_1-z_2}{2}  \right) \notag
\end{align}
which is a diffeomorphism and $l^*\omega_\text{std}=\frac{1}{2}\omega_\text{std}$. This yields that $L:=l(\gamma\times\gamma)$ is also a smooth Lagrangian torus in $(\mathbb{C}^2,\omega_\text{std})$. 

Applying the rotation on the second factor 
\begin{align}\label{rotation_phi}
R_\phi:\mathbb{C}^2&\to\mathbb{C}^2 \\
(z_1,z_2)&\mapsto(z_1,z_2 e^{i\phi})\notag
\end{align}
to the Lagrangian torus $L$ we obtain another Lagrangian torus $L_\phi:=R_\phi(L)$. They intersect cleanly along $\gamma\times\{0\}$, i.e.
\[
T_pL\cap T_p L_\phi= T_p(\gamma\times\{0\})\qquad \forall\,p\in \gamma\times\{0\}.
\]
There is a symplectic involution 
\begin{align}\label{symp_involution}
\varrho:\mathbb{C}^2&\to\mathbb{C}^2\\
(z_1,z_2)&\mapsto(z_1,-z_2)\notag
\end{align} 
preserving both $L$ and $L_\phi$ and fixing $\gamma\times\{0\}$.

The goal is to find an intersection point away from $\gamma\times\{0\}$, i.e. one needs to show that
\begin{equation}\label{a_priori_exis}
\# (L \cap L_\phi \setminus (\gamma\times\{0\}) ) \geq 1 .
\end{equation}
Suppose such an intersection point exists and reads in the coordinate as $(z_1,z_2 e^{i\phi})$, then it follows that four points $z_1\pm z_2$, $z_1\pm z_2 e^{i\phi}$ all lie on the curve $\gamma$ and are four vertices of a rectangle whose diagonals meet at the angle $\phi\in(0,\frac{\pi}{2}]$ in $\mathbb{C}$.

The union  $L\cup L_\phi$ has at least $\gamma\times\{0\}$ in the self-intersection set. After a Lagrangian smoothing surgery of $L\cup L_\phi$ along $\gamma\times\{0\}$ (c.f. \cite[Proposition 1.1]{Greene-Lobb2021}), one gets a smoothly immersed Lagrangian torus $L\#_{\gamma\times\{0\}} L_\phi$
in $\mathbb{C}^2$ that coincides with $L\cup L_\phi$ outside of a small open neighborhood $\mathcal{O}p(\gamma\times\{0\})$, and it is
disjoint from $\mathbb{C}\times\{0\}$. The involution $\varrho$ now acts as a fixed-point free involution on the immersed Lagrangian torus $L\#_{\gamma\times\{0\}}L_\phi$. 

In \cite{Greene-Lobb2021}, the following map is considered:
\begin{align*}
g:\mathbb{C}^2&\to\mathbb{C}^2\\
(z_1,z_2)&\mapsto\left(z_1,\frac{1}{\sqrt{2}}\frac{z_2^2}{|z_2|}\right)
\end{align*}
which is a symplectomorphism away from $\mathbb{C}\times\{0\}$. Applying to the Lagrangian tori $L$ and $L_\phi$, their image $g(L)$ and $g(L_\phi)$ are both homeomorphic to a Möbius strip $\operatorname{Sym}^2(\gamma)$, and are smooth Lagrangian away from $\mathbb{C}\times\{0\}$. The image $g(L\#_{\gamma\times\{0\}} L_\phi)$ is therefore a smoothly immersed Lagrangian Klein bottle $K\looparrowright\mathbb{C}^2$ such that 
\[
g(L\#_{\gamma\times\{0\}} L_\phi\setminus \mathcal{O}p(\gamma\times\{0\}))\cong g(L\cup L_\phi\setminus \mathcal{O}p(\gamma\times\{0\}))
\] 
where $\mathcal{O}p(\gamma\times\{0\})$ is a small open neighborhood of $\gamma\times\{0\}$ on which the immersion restricted is an smooth embedding in $\mathbb{C}^2$. Using the result of non-existence of embedded Lagrangian Klein bottle in $\mathbb{C}^2$, c.f. \cite{Shevchishin2009,Nemirovski2009}, they conclude that the set of self-intersection points of the immersed Lagrangian Klein bottle must be non-empty, and hence the existence of retangular peg is a corollary. 

\begin{thm}[\cite{Greene-Lobb2021}]\label{existence_result}
For every smooth Jordan curve and rectangle with any aspect angle $\phi\in(0,\frac{\pi}{2}]$ in the $\mathbb{R}^2$, there exists a similar rectangle whose vertices land on the curve.
\end{thm}

\noindent In \cite{Greene-Lobb2023}, the above existence result is generalized to arbitrary cyclic quadrilateral in a more general geometric setup. Instead of using the Lagrangian Klein bottle argument, they compute the minimal Maslov number of the surgered immersed Lagrangian torus, which is more straightforward, see \S\ref{Section 6}. We will give an alternative Floer theoretical proof of the existence theorem in \S\ref{Section 4} circumventing applying the Lagrangian surgery and using the result of the non-existence of embedded Lagrangian Klein bottle in $\mathbb{C}^2$ \cite{Shevchishin2009,Nemirovski2009} as well as the result of minimal Maslov number of embedded Lagrangian tori in $\mathbb{C}^2$ \cite{Viterbo1990,Polterovich1991}.

Beyond addressing the existence problem, the investigation into the multiplicity of inscribed rectangles remains largely unexplored in general. Leveraging the aforementioned existence result of the rectangular pegs by Greene and Lobb, our primary contribution lies in the following main theorem.

\begin{thm} \label{main}
For any fixed $\phi\in(0,\frac{\pi}{2})$, there is a generic class of smooth Jordan curves such that any rectangle with the fixed aspect angle $\phi$ in $\mathbb{R}^2$ admits two similar rectangles whose vertices land distinctly on any curve in this generic class. 
\end{thm}

\begin{rem}
A counter example of $\phi=\frac{\pi}{2}$ for the statement: consider an ellipse (noticing that ellipses are in the generic class of curves, see Example \ref{exa_ellipses}) and it is obvious that there is only one inscribed square. It is shown by Schnirelman \cite{Schnirelman1944} that any generic smooth Jordan curve inscribes an odd number of squares. See also \cite[p.348]{Matschke2014} for exposition. On the other hand, we do not know at the moment an example of smooth Jordan curve which carries only one inscribed rectangle with aspect angle $\phi\in(0,\frac{\pi}{2})$ due to the possible scenario discussed in Remark \ref{rem_trans_pt}.
\end{rem}

\begin{rem}
In fact, such generic doubling result can be further extended to inscribed cyclic quadrilaterals by adopting the geometric setup in \cite{Greene-Lobb2023}. See \S\ref{Section 6}, Theorem \ref{thm_doubling_cyc_quad}.
\end{rem}

We present two distinct proofs for Theorem \ref{main} in \S\ref{Section 5}. These proofs rely essentially on two propositions: the correspondence between intersection points in $L\cap L_\phi$ away from $\mathbb{C}\times\{0\}$, as described in \S\ref{Section 2}, and a generic geometric transversality property for $L$ and $L_\phi$, elaborated in \S\ref{Section 3}. We sketch the proofs as follows. The first proof involves computing the Euler characteristic of the Lagrangian Floer homology for the Lagrangian tori $L$ and $L_\phi$, generated from a generic curve $\gamma$, with the geometric transversality property. This approach applies the invariance of the Lagrangian Floer homology to obtain a contradiction when assuming the existence of only one intersection point away from $\mathbb{C}\times\{0\}$. The second proof follows a similar path of contradiction but employs a novel formula for the algebraic intersection number of cleanly intersecting Lagrangian submanifolds. In \S\ref{Section 6}, we further extend the main result to inscribed cyclic quadrilaterals. 

In Appendix \ref{App_Floer_homology}, we provide an overview of the Lagrangian Floer homology adapted to our geometric setup. Additionally, Appendix \ref{Appendix_B} furnishes a detailed derivation of the computation formula for the algebraic intersection number of cleanly intersecting submanifolds.

\subsubsection*{Acknowledgment}
I am deeply indebted to my Ph.D. supervisor, Urs Frauenfelder, for the enlightening conversations and insightful ideas provided throughout this research. I am grateful to Kai Cieliebak, Georgios Dimitroglou Rizell, Joshua Evan Greene, Felix Schlenk, and Felix Schmäschke for the valuable comments and fruitful discussions. I thank Daniel Cristofaro-Gardiner for raising interesting question on non-generic cases. I also thank Marián Poppr for pointing out the reference \cite{Karlsson2013} to me, and Martin Konrad for indicating typos and errors in the earlier draft. The authors acknowledge the partial support from DFG grant FR 2637/4-1.

\section{Intersection points and rectangular pegs}\label{Section 2}
Intersection points in $L\cap L_\phi$ away from $\mathbb{C}\times\{0\}$ do not one-to-one correspond to distinguished inscribed rectangles on the smooth Jordan curve $\gamma$.

The tori $L$ and $L_\phi$ in $\mathbb{C}^2$ are parameterized by $S^1\times S^1$:
\begin{align*}
S^1\times S^1 &\to L\subset\mathbb{C}^2\\ 
(t_1,t_2)&\mapsto \left(\frac{\gamma(t_1)+\gamma(t_2)}{2},\frac{\gamma(t_1)-\gamma(t_2)}{2} \right) , \\
S^1\times S^1 &\to L_\phi\subset\mathbb{C}^2\\
(t_3,t_4)&\mapsto \left(\frac{\gamma(t_3)+\gamma(t_4)}{2},\frac{\gamma(t_3)-\gamma(t_4)}{2} e^{i\phi}\right) .
\end{align*}
An intersection point away from $\gamma\times\{0\}$ is determined by the following condition for $t_1\neq t_2\neq t_3\neq t_4$ in $S^1$:
\begin{align}
\gamma(t_1)+\gamma(t_2)&=\gamma(t_3)+\gamma(t_4), \label{intsec_eq_1}\\
\gamma(t_1)-\gamma(t_2)    &=(\gamma(t_3)-\gamma(t_4)) e^{i\phi}. \label{intsec_eq_2}
\end{align}
The ordered quadruple $(t,t',t'',t''')$ solving above equations gives rise to the intersection point \[p=(z_1,z_2 e^{i\phi})=\left(\frac{\gamma(t'')+\gamma(t''')}{2},\frac{\gamma(t'')-\gamma(t''')}{2}e^{i\phi}\right)\in L\cap L_\phi\] and parameterizes the four points on the curve $\gamma$ which completely determine a geometrically distinct inscribed rectangle as their vertices $(A,B,C,D)$ in the counterclockwise order:
\[
(z_1+z_2,z_1+z_2 e^{i\phi},z_1-z_2,z_1-z_2e^{i\phi}).
\]
See Figure \ref{Fig_1}.

\begin{figure}
\[
\begin{tikzpicture}[x=0.75pt,y=0.75pt,yscale=-1,xscale=1]
	
	\draw   (248,27) -- (414,27) -- (414,103) -- (248,103) -- cycle ;
	\draw    [dash pattern={on 0.84pt off 2.51pt}] (248,27) -- (414,103) ;
	\draw    [dash pattern={on 0.84pt off 2.51pt}] (248,103) -- (414,27) ;
	
	\draw (391,106.4) node [anchor=north west][inner sep=0.75pt]    {$z_{1} +z_{2}=A$};
	\draw (387,3.4) node [anchor=north west][inner sep=0.75pt]    {$z_{1} +z_{2} \cdot e^{i\phi }=B$};
	\draw (215,4.4) node [anchor=north west][inner sep=0.75pt]    {$C=z_{1} -z_{2}$};
	\draw (203,105.4) node [anchor=north west][inner sep=0.75pt]    {$D=z_{1} -z_{2} \cdot e^{i\phi }$};
	\draw (350.01,55.72) node [anchor=north west][inner sep=0.75pt]    {$\phi $};

\end{tikzpicture}
\]
\caption{Inscribed rectangle in $\mathbb{C}$ with aspect angle $\phi$ and vertices $(A,B,C,D)$ ordered counterclockwise.}\label{Fig_1}
\end{figure}

\begin{figure}
	\[
	\begin{tikzpicture}[x=0.75pt,y=0.75pt,yscale=-1,xscale=1]
		
		\draw  [dash pattern={on 0.84pt off 2.51pt}]  (189.11,28.62) -- (279.11,82.62) ;
		\draw  [dash pattern={on 0.84pt off 2.51pt}]  (279.11,28.62) -- (189.11,82.62) ;
		\draw   (301.28,55.85) -- (364.94,55.85) ;
		\draw [shift={(366.94,55.85)}, rotate = 180] [color={rgb, 255:red, 0; green, 0; blue, 0 }  ][line width=0.75]    (10.93,-3.29) .. controls (6.95,-1.4) and (3.31,-0.3) .. (0,0) .. controls (3.31,0.3) and (6.95,1.4) .. (10.93,3.29)   ;
		\draw   (189.11,28.62) -- (279.11,28.62) -- (279.11,82.62) -- (189.11,82.62) -- cycle ;
		\draw   [dash pattern={on 0.84pt off 2.51pt}] (387.11,28.15) -- (477.11,82.15) ;
		\draw   [dash pattern={on 0.84pt off 2.51pt}] (477.11,28.15) -- (387.11,82.15) ;
		\draw   (387.11,28.15) -- (477.11,28.15) -- (477.11,82.15) -- (387.11,82.15) -- cycle ;
		
		\draw (257.71,87.87) node [anchor=north west][inner sep=0.75pt]    {$z_{1} +z_{2}=A$};
		\draw (252.94,3.38) node [anchor=north west][inner sep=0.75pt]    {$z_{1} +z_{2} e^{i\phi }=B$};
		\draw (160.06,4.42) node [anchor=north west][inner sep=0.75pt]    {$C=z_{1} -z_{2}$};
		\draw (147.52,86.42) node [anchor=north west][inner sep=0.75pt]    {$D=z_{1} -z_{2} e^{i\phi }$};
		\draw (464.71,89.4) node [anchor=north west][inner sep=0.75pt]    {$z_{1} +z_{2}=A'$};
		\draw (450.94,2.9) node [anchor=north west][inner sep=0.75pt]    {$z_{1} +z_{2} e^{i\phi }=B'$};
		\draw (358.06,3.95) node [anchor=north west][inner sep=0.75pt]    {$C'=z_{1} -z_{2}$};
		\draw (345.52,85.95) node [anchor=north west][inner sep=0.75pt]    {$D'=z_{1} -z_{2} e^{i\phi }$};

	\end{tikzpicture}
	\]
	\caption{Vertices relabeled: $(A',B',C',D')=(A,B,C,D)$}\label{Fig_2}
\end{figure}

\begin{figure}
	\[
	\begin{tikzpicture}[x=0.75pt,y=0.75pt,yscale=-1,xscale=1]
		
		\draw  [dash pattern={on 0.84pt off 2.51pt}]  (189.11,28.62) -- (279.11,82.62) ;
		\draw  [dash pattern={on 0.84pt off 2.51pt}]  (279.11,28.62) -- (189.11,82.62) ;
		\draw    (301.28,55.85) -- (364.94,55.85) ;
		\draw [shift={(366.94,55.85)}, rotate = 180] [color={rgb, 255:red, 0; green, 0; blue, 0 }  ][line width=0.75]    (10.93,-3.29) .. controls (6.95,-1.4) and (3.31,-0.3) .. (0,0) .. controls (3.31,0.3) and (6.95,1.4) .. (10.93,3.29)   ;
		\draw   (189.11,28.62) -- (279.11,28.62) -- (279.11,82.62) -- (189.11,82.62) -- cycle ;
		\draw  [dash pattern={on 0.84pt off 2.51pt}]  (387.11,28.15) -- (477.11,82.15) ;
		\draw  [dash pattern={on 0.84pt off 2.51pt}]  (477.11,28.15) -- (387.11,82.15) ;
		\draw   (387.11,28.15) -- (477.11,28.15) -- (477.11,82.15) -- (387.11,82.15) -- cycle ;
		
		\draw (257.71,87.87) node [anchor=north west][inner sep=0.75pt]    {$z_{1} +z_{2}=A$};
		\draw (252.94,3.38) node [anchor=north west][inner sep=0.75pt]    {$z_{1} +z_{2} e^{i\phi }=B$};
		\draw (160.06,4.42) node [anchor=north west][inner sep=0.75pt]    {$C=z_{1} -z_{2}$};
		\draw (147.52,86.42) node [anchor=north west][inner sep=0.75pt]    {$D=z_{1} -z_{2} e^{i\phi }$};
		\draw (362.71,5.4) node [anchor=north west][inner sep=0.75pt]    {$C'=z_{1} +z_{2}$};
		\draw (360.94,86.9) node [anchor=north west][inner sep=0.75pt]    {$D'=z_{1} +z_{2} e^{i\phi }$};
		\draw (460.06,85.95) node [anchor=north west][inner sep=0.75pt]    {$z_{1} -z_{2}=A'$};
		\draw (456.52,3.95) node [anchor=north west][inner sep=0.75pt]    {$z_{1} -z_{2} e^{i\phi }=B'$};

	\end{tikzpicture}
	\]
	\caption{Vertices relabeled: $(A',B',C',D')=(C,D,A,B)$}\label{Fig_3}
\end{figure}

\begin{figure}
	\[
	\begin{tikzpicture}[x=0.75pt,y=0.75pt,yscale=-1,xscale=1]
		
		\draw   (202.61,24.12) -- (265.61,24.12) -- (265.61,87.12) -- (202.61,87.12) -- cycle ;
		\draw  [dash pattern={on 0.84pt off 2.51pt}]  (202.61,24.12) -- (265.61,87.12) ;
		\draw  [dash pattern={on 0.84pt off 2.51pt}]  (265.61,24.12) -- (202.61,87.12) ;
		\draw    (301.28,55.85) -- (364.94,55.85) ;
		\draw [shift={(366.94,55.85)}, rotate = 180] [color={rgb, 255:red, 0; green, 0; blue, 0 }  ][line width=0.75]    (10.93,-3.29) .. controls (6.95,-1.4) and (3.31,-0.3) .. (0,0) .. controls (3.31,0.3) and (6.95,1.4) .. (10.93,3.29)   ;
		\draw   (394.61,25.1) -- (457.61,25.1) -- (457.61,88.1) -- (394.61,88.1) -- cycle ;
		\draw  [dash pattern={on 0.84pt off 2.51pt}]  (394.61,25.1) -- (457.61,88.1) ;
		\draw  [dash pattern={on 0.84pt off 2.51pt}]  (457.61,25.1) -- (394.61,88.1) ;
		
		\draw (247.71,88.87) node [anchor=north west][inner sep=0.75pt]    {$z_{1} +z_{2}=A$};
		\draw (248.94,2.38) node [anchor=north west][inner sep=0.75pt]    {$z_{1} +iz_{2}=B$};
		\draw (172.06,1.42) node [anchor=north west][inner sep=0.75pt]    {$C=z_{1} -z_{2}$};
		\draw (166.52,88.42) node [anchor=north west][inner sep=0.75pt]    {$D=z_{1} -iz_{2}$};
		\draw (444.94,88.4) node [anchor=north west][inner sep=0.75pt]    {$z_{1} +iz_{2}=A'$};
		\draw (441.06,4.4) node [anchor=north west][inner sep=0.75pt]    {$z_{1} -z_{2}=B'$};
		\draw (361.52,4.4) node [anchor=north west][inner sep=0.75pt]    {$C'=z_{1} -iz_{2}$};
		\draw (365.71,88.4) node [anchor=north west][inner sep=0.75pt]    {$D'=z_{1} +z_{2}$};

	\end{tikzpicture}
	\]
	\caption{Vertices relabeled: $(A',B',C',D')=(B,C,D,A)$}\label{Fig_4}
\end{figure}

\begin{figure}
	\[
	\begin{tikzpicture}[x=0.75pt,y=0.75pt,yscale=-1,xscale=1]
		
		\draw   (202.61,24.12) -- (265.61,24.12) -- (265.61,87.12) -- (202.61,87.12) -- cycle ;
		\draw  [dash pattern={on 0.84pt off 2.51pt}]  (202.61,24.12) -- (265.61,87.12) ;
		\draw  [dash pattern={on 0.84pt off 2.51pt}]  (265.61,24.12) -- (202.61,87.12) ;
		\draw    (301.28,55.85) -- (364.94,55.85) ;
		\draw [shift={(366.94,55.85)}, rotate = 180] [color={rgb, 255:red, 0; green, 0; blue, 0 }  ][line width=0.75]    (10.93,-3.29) .. controls (6.95,-1.4) and (3.31,-0.3) .. (0,0) .. controls (3.31,0.3) and (6.95,1.4) .. (10.93,3.29)   ;
		\draw   (394.61,25.1) -- (457.61,25.1) -- (457.61,88.1) -- (394.61,88.1) -- cycle ;
		\draw  [dash pattern={on 0.84pt off 2.51pt}]  (394.61,25.1) -- (457.61,88.1) ;
		\draw  [dash pattern={on 0.84pt off 2.51pt}]  (457.61,25.1) -- (394.61,88.1) ;
		
		\draw (247.71,88.87) node [anchor=north west][inner sep=0.75pt]    {$z_{1} +z_{2}=A$};
		\draw (248.94,2.38) node [anchor=north west][inner sep=0.75pt]    {$z_{1} +iz_{2}=B$};
		\draw (172.06,1.42) node [anchor=north west][inner sep=0.75pt]    {$C=z_{1} -z_{2}$};
		\draw (166.52,88.42) node [anchor=north west][inner sep=0.75pt]    {$D=z_{1} -iz_{2}$};
		\draw (442.71,3.4) node [anchor=north west][inner sep=0.75pt]    {$z_{1} +z_{2}=B'$};
		\draw (443.52,89.4) node [anchor=north west][inner sep=0.75pt]    {$z_{1} -iz_{2}=A'$};
		\draw (360.94,3.4) node [anchor=north west][inner sep=0.75pt]    {$C'=z_{1} +iz_{2}$};
		\draw (363.06,89.4) node [anchor=north west][inner sep=0.75pt]    {$D'=z_{1} -z_{2}$};

	\end{tikzpicture}
	\] 
	\caption{Vertices relabeled: $(A',B',C',D')=(D,A,B,C)$}\label{Fig_5}
\end{figure}

Suppose we have two intersection points $(z_1,z_2 e^{i\phi})$ and $(z'_1,z'_2 e^{i\phi})$ with respect to a fixed aspect angle $\phi$, and they constitute the same geometric distinct inscribed rectangle on $\gamma$. There are the following three cases:

\begin{itemize}
	\item[(1)] $z'_1+z'_2=z_1+z_2$: By the counterclockwise ordering of vertices, it yields $z'_1=z_1$, $z'_2=z_2$, i.e. these two intersection points are the same one. This means the intersection point $(z'_1,z'_2 e^{i\phi})$ is the image of identity map of $(z_1,z_2 e^{i\phi})$. See Figure \ref{Fig_2}.

	\item[(2)] $z'_1+z'_2=z_1-z_2$: By the counterclockwise ordering of vertices, it yields $z'_1=z_1$, $z'_2=-z_2$. This means the intersection point $(z'_1,z'_2 e^{i\phi})$ is the image under 
	\[
	\varrho:(z_1,z_2)\mapsto(z_1,-z_2)
	\] 
	of the intersection point $(z_1,z_2 e^{i\phi})\in L\cap L_\phi$ in the reflected Lagrangian tori $\varrho(L)\cap\varrho(L_\phi)$, and since $L$ and $L_\phi$ is invariant under $\varrho$, so it is another different intersection point in $L\cap L_\phi$, and it gives literally the same inscribed rectangle. See Figure \ref{Fig_3}. 

	\item[(3)] $z'_1+z'_2=z_1\pm z_2 e^{i\phi}$: By the counterclockwise ordering of vertices, this is only possible when $\phi=\frac{\pi}{2}$, otherwise 
	\[
	(A,B,C,D)=(z_1+z_2,z_1+z_2 e^{i\phi},z_1-z_2,z_1-z_2e^{i\phi})
	\]
	and
	\[
	(A',B',C',D')=(z'_1+z'_2,z'_1+z'_2 e^{i\phi},z'_1-z'_2,z'_1-z'_2e^{i\phi})
	\]
	would constitute two geometrically distinct rectangles in $\mathbb{C}$, in which $(A,B,C,D)$ inscribes on $\gamma$ whereas $(A',B',C',D')$ does not. In the case of $\phi=\frac{\pi}{2}$, there are $z'_1=z_1$, $z'_2=\pm z_2 e^{i\frac{\pi}{2}}$. 
	\begin{itemize}
		\item $(z'_1,z'_2)=(z_1,iz_2)$: It means the intersection point $(z'_1,z'_2 e^{i\phi})$ is the image of the map 
		\[
		\iota:(z_1,z_2)\mapsto(z_1,iz_2)
		\] 
		and $\iota^2=\varrho$. See Figure \ref{Fig_4}.

		\item $(z'_1,z'_2)=(z_1,-iz_2)$: It means the intersection point $(z'_1,z'_2 e^{i\phi})$ is the image of the map 
		\[
		\varrho\circ\iota:(z_1,z_2)\mapsto(z_1,-iz_2)
		\] 
		and $\iota\circ\varrho=\varrho\circ\iota$. See Figure \ref{Fig_5}.
	\end{itemize}
\end{itemize}

The above discussion in turns proves the following statement.

\begin{prop}\label{prop_aut_of_point}
	Let $\operatorname{Rect}(\gamma,\phi)$ be the set of inscribed rectangles with fixed aspect angle $\phi$ on the curve $\gamma$.
	\begin{itemize}
		\item[(i)] When $\phi\in(0,\frac{\pi}{2})$, there is a bijection between sets
		\[
		(L\cap L_\phi\setminus(\gamma\times\{0\}))/C_2\xrightarrow{\cong}\operatorname{Rect}(\gamma,\phi)
		\]
		where $C_2$ is presented by 
		\[
		\langle \operatorname{id},\varrho \mid \varrho^2=\operatorname{id}\rangle\cong\mathbb{Z}/2\mathbb{Z}.
		\]
		\item[(ii)] When $\phi=\frac{\pi}{2}$, there is a bijection between sets
		\[
		(L\cap L_\phi\setminus(\gamma\times\{0\}))/C_4\xrightarrow{\cong}\operatorname{Rect}(\gamma,\phi)
		\]
		where $C_4$ is presented by 
		\[
		\langle \operatorname{id},\varrho,\iota \mid \varrho^2=\operatorname{id},\iota^2=\varrho\rangle\cong\mathbb{Z}/4\mathbb{Z}.
		\]
	\end{itemize}
\end{prop}

\begin{rem}
Recall the ordered tuple $(t,t',t'',t''')$ determining a geometric distinct rectangle inscribed on the curve $\gamma$ with vertices ordered as 
\[
(A,B,C,D)=(z_1+z_2,z_1+z_2 e^{i\phi},z_1-z_2,z_1-z_2e^{i\phi})
\]
mentioned at the beginning of this section. The automorphism group of $\{t,t',t'',t'''\}$ as a set is the symmetric group $S_4$. Some of the elements in $S_4$ actually give rise to the reparameterizations of the vertices of the fixed geometrically distinct inscribed rectangle on $\gamma$, and further correspond to the symmetries of the intersection points in $L\cap L_\phi\setminus(\gamma\times\{0\})$, i.e. the elements in $C_2$ and $C_4$ above. We describe these symmetries in the following:
\begin{itemize}
	\item[(1)] $\begin{pmatrix}
		1 & 2 &  3 &	4\\
		1 & 2 &  3 &	4	
	\end{pmatrix}:(t,t',t'',t''')\mapsto(t,t',t'',t''')$ obviously corresponds to the identity map. 
	\item[(2)] $\begin{pmatrix}
		1 & 2 &  3 &	4\\
		2 & 1 &  4 &	3	
	\end{pmatrix}:(t,t',t'',t''')\mapsto(t',t,t''',t'')$. This corresponds to the symplectic involution $\varrho$ acting on the intersection point $p\in L\cap L_\phi$, and the image $q=\varrho(p)$ is another intersection point since $(t',t,t''',t'')$ also solves the equations \eqref{intsec_eq_1}, \eqref{intsec_eq_2} as $(t,t',t'',t''')$ does. The point $q\in L\cap L_\phi$ determines four point in $\mathbb{C}$ in counterclockwise order
	\[
	(A',B',C',D')=(z_1-z_2,z_1-z_2e^{i\phi},z_1+z_2,z_1+z_2e^{i\phi}).
	\]
	They form the exact same rectangle as the original one $(A,B,C,D)$ inscribed on the curve $\gamma$ with relabeled vertices. 
	\item[(3)] $\bullet$ $\begin{pmatrix}
		1 & 2 &  3 &	4\\
		3 & 4 &  1 &	2	
	\end{pmatrix}:(t,t',t'',t''')\mapsto(t'',t''',t,t')$ corresponds to $\iota$ acting on $p\in L\cap L_\phi$  when $\phi=\frac{\pi}{2}$. This can be seen as follows.  Since $t\neq t'\neq t''\neq t'''$, the pairs $(t'',t''')$ and $(t,t')$ give rise to two points $p'\in L$ and $p''\in L_\phi$ respectively. In particular, $p'=R_{-\phi}(p)$, $p''=R_{\phi}(p)$, and $p'\neq p''$. When $\phi\in(0,\frac{\pi}{2})$, this symmetric group element maps the original inscribed rectangle to a new rectangle in the plane which is no longer inscribed on the curve $\gamma$. When $\phi=\frac{\pi}{2}$, i.e. the rectangle is actually a square, we obtain two additional intersection points $p''=R_{\frac{\pi}{2}}(p)\in R_{\frac{\pi}{2}}( L\cap L_{\frac{\pi}{2}})=L_{\frac{\pi}{2}}\cap R_{\pi}(L)=L_{\frac{\pi}{2}}\cap L$ and $p'=R_{-\frac{\pi}{2}}(p)\in R_{-\frac{\pi}{2}}( L\cap L_{\frac{\pi}{2}})=R_{\pi}(L_{\frac{\pi}{2}})\cap L=L_{\frac{\pi}{2}}\cap L$ so that $p''=R_{\pi}(p')=\varrho(p')$, they commonly determine four points in $\mathbb{C}$ in counterclockwise order:
	\[
	(A',B',C',D')=(z_1+iz_2,z_1-z_2,z_1-iz_2,z_1+z_2).
	\] 
	This form the exact same square as the original one $(A,B,C,D)$ inscribed on the curve $\gamma$  with relabeled vertices.
		
\noindent $\bullet$ $\begin{pmatrix}
		1 & 2 &  3 &	4\\
		4 & 3 &  2 &	1	
	\end{pmatrix}:(t,t',t'',t''')\mapsto(t''',t'',t',t)$ corresponds to $\varrho\circ\iota$ acting on $p\in L\cap L_\phi$ when $\phi=\frac{\pi}{2}$, and note that  
	\begin{align*}
     \begin{pmatrix}
     	1 & 2 &  3 &	4\\
     	4 & 3 &  2 &	1	
     \end{pmatrix}=	\begin{pmatrix}
     	1 & 2 &  3 &	4\\
     	2 & 1 &  4 &	3	
     \end{pmatrix}\begin{pmatrix}
     	1 & 2 &  3 &	4\\
     	3 & 4 &  1 &	2	
     \end{pmatrix},\\
 	\begin{pmatrix}
 		1 & 2 &  3 &	4\\
 		2 & 1 &  4 &	3	
 	\end{pmatrix}\begin{pmatrix}
 		1 & 2 &  3 &	4\\
 		3 & 4 &  1 &	2	
 	\end{pmatrix}=\begin{pmatrix}
 		1 & 2 &  3 &	4\\
 		3 & 4 &  1 &	2	
 	\end{pmatrix}\begin{pmatrix}
 		1 & 2 &  3 &	4\\
 		2 & 1 &  4 &	3	
 	\end{pmatrix},
	\end{align*}
	such relation amounts to $\varrho\circ\iota=\iota\circ\varrho$. This can be seen as follows. Since $t\neq t'\neq t''\neq t'''$, the pairs $(t''',t'')$ and $(t',t)$ give rise to two points $q'\in L$ and $q''\in L_\phi$ respectively. In particular, $q'=R_{-\pi-\phi}(p)$, $q''=R_{\pi+\phi}(p)$, and $q'\neq q''$. When $\phi\in(0,\frac{\pi}{2})$, this symmetric group element maps the original inscribed rectangle to a new rectangle in the plane which is no longer inscribed on the curve $\gamma$. When $\phi=\frac{\pi}{2}$, i.e. the rectangle is actually a square, the point $q''=R_{\pi+\frac{\pi}{2}}(p)=R_\pi(p'')=p'\in L_{\frac{\pi}{2}}\cap L$ determines four points in $\mathbb{C}$ in counterclockwise order:
	\[
	(A',B',C',D')=(z_1-iz_2,z_1+z_2,z_1+iz_2,z_1-z_2).
	\] 
	They form the exact same square as the original one $(A,B,C,D)$ inscribed on the curve $\gamma$ with relabeled vertices. 
\end{itemize} 
A caveat is that the above set-theoretic bijection of the symmetric group elements in $S_4$, which form a Klein 4 group $K_4$, to the group elements in $C_4$ does not constitute a group homomorphism, since $K_4$ is not isomorphic to $C_4$.  
\end{rem}

It is worth noting that the symmetry of rectangular pegs is also addressed in \cite[\S 4.1]{Matschke2022}.

\section{Geometric transversality}\label{Section 3}

In general, we have no knowledge of the intersection set of $L$ and $L_\phi$ away from $\gamma\times\{0\}$. It can be badly degenerate, e.g. the example below.
\begin{exmp}
	In the most special case of $\gamma$ being the unit circle $S^1$ in $\mathbb{C}$, the tori $L$ and $L_\phi$ have badly degenerate (not clean) intersection in $\mathbb{C}^2$ away from $\gamma\times\{0\}$. This can be seen explicitly via the polar coordinate of the torus: 
	\begin{align*}
		\left(\frac{r e^{i\theta_1}+r e^{i\theta_2}}{2},\frac{r e^{i\theta_1}-r e^{i\theta_2}}{2}\right)&\in L, \\
		\left(\frac{r e^{i\theta_3}+r e^{i\theta_4}}{2} ,
		\frac{r e^{i\theta_3}-r e^{i\theta_4}}{2} e^{i\phi} \right)&\in L_\phi.
	\end{align*} 
where $\theta_1,\theta_2,\theta_3,\theta_4\in[0,2\pi]$, and $r$ is the radius of the standard circle $\gamma$. Despite the clean intersection $S^1\times\{0\}$ in $\mathbb{C}\times\{0\}$, there is another part of the intersection as a circle $\{0\}\times S^1$ in $\{0\}\times\mathbb{C}$ which is invariant under the symplectic involution \eqref{symp_involution}. Each point in the circle $\{0\}\times S^1$ gives rise to a geometrically distinct inscribed rectangle on $S^1\subset\mathbb{C}$ with aspect angle $\phi$. 
\end{exmp}

We show in the following that $L$ and $L_\phi$ intersect transversely away from $\gamma\times\{0\}$ under certain generic condition. The result plays a pivotal role in this paper.

\begin{prop}\label{prop_gen_imm}
For any fixed $\phi\in(0,\frac{\pi}{2}]$, there is a residual set of smooth Jordan curves such that for every $\gamma$ in the set, the corresponding smooth Lagrangian tori $L$ and $L_\phi$ constructed in above way have only finitely many transversal intersection points away from $\gamma\times\{0\}$. 
\end{prop}

\begin{proof} Let $\operatorname{Emb}^k(S^1,\mathbb{C})$ be the space of $C^k$-smooth embedded closed curves in $\mathbb{C}$ for $k\in\mathbb{N}$. There are open inclusions $\operatorname{Emb}^k(S^1,\mathbb{C})\subset\operatorname{Imm}^k(S^1,\mathbb{C})\subset C^k(S^1,\mathbb{C})$. We consider the $C^k$-smooth map between $C^\infty$-smooth Banach manifolds
\begin{align*}
\mathcal{F}: \operatorname{Emb}^k(S^1,\mathbb{C})\times ( (S^1)^4\setminus \Delta_{1234}^{(S^1)^4} ) &\to \mathbb{C}^2\times\mathbb{C}^2\\
(\gamma,t_1,t_2,t_3,t_4) & \mapsto \left( l(\gamma(t_1)\times\gamma(t_2)) ,  R_\phi (l(\gamma(t_3)\times\gamma(t_4))) \right)
\end{align*} 
so that
\begin{align*}
\operatorname{pr}_1( \mathcal{F}(\gamma,t_1,t_2,t_3,t_4))&= \frac{\gamma(t_1)+\gamma(t_2)}{2} ,\\ 
\operatorname{pr}_2( \mathcal{F}(\gamma,t_1,t_2,t_3,t_4))&=\frac{\gamma(t_1)-\gamma(t_2)}{2} , \\ 
\operatorname{pr}_3( \mathcal{F}(\gamma,t_1,t_2,t_3,t_4))&=\frac{\gamma(t_3)+\gamma(t_4)}{2} , \\ 
\operatorname{pr}_4( \mathcal{F}(\gamma,t_1,t_2,t_3,t_4))&=\frac{\gamma(t_3)-\gamma(t_4)}{2}e^{i\phi} ,
\end{align*}
and 
\[
\Delta_{1234}^{(S^1)^4}=\{ (t_1,t_2,t_3,t_4)\in (S^1)^4 : t_1=t_2=t_3=t_4 \}
\] 
is the small diagonal of the 4-torus $(S^1)^4$. Consider the diagonal in $\mathbb{C}^2\times\mathbb{C}^2$:
\[
\Delta_{\mathbb{C}^2}=\{ (z_1,z_2,z_3,z_4)\in\mathbb{C}^2\times\mathbb{C}^2: (z_1,z_2)= (z_3,z_4) \}
\] 
We want to show that $\mathcal{F}^{-1}(\Delta_{\mathbb{C}^2})$ is a manifold, and then further project to $C^{k}(S^1,\mathbb{C})$ in order to obtain the genericity result.

The linearization of map $\mathcal{F}$ at any point $(\gamma,\mathbf{t})=(\gamma,t_1,t_2,t_3,t_4)$ away from the small diagonal $\Delta_{1234}^{(S^1)^4}$ is given by
\begin{align*}
D_{(\gamma,\mathbf{t})}\mathcal{F}:T_{\gamma}C^k(S^1,\mathbb{C})\times T_{\mathbf{t}}(S^1)^4&\to T_{\mathcal{F}(\mathcal{\gamma,\mathbf{t}})}\mathbb{C}^4\\
(\xi,v_1,v_2,v_3,v_4)&\mapsto(V_1,V_2,V_3,V_4)
\end{align*}
where
\begin{align*}
V_1&=\frac{1}{2}(D_{t_1}\gamma(v_1)+D_{t_2}\gamma(v_2)+\xi_{\gamma(t_1)}+\xi_{\gamma(t_2)}),\\
V_2&=\frac{1}{2}(D_{t_1}\gamma(v_1)-D_{t_2}\gamma(v_2)+\xi_{\gamma(t_1)}-\xi_{\gamma(t_2)}),\\
V_3&=\frac{1}{2}(D_{t_3}\gamma(v_3)+D_{t_4}\gamma(v_4)+\xi_{\gamma(t_3)}+\xi_{\gamma(t_4)}),\\
V_4&=\frac{1}{2}(D_{t_3}\gamma(v_3)-D_{t_4}\gamma(v_4)+\xi_{\gamma(t_3)}-\xi_{\gamma(t_4)}),
\end{align*}
in which we regard $\xi\in T_\gamma C^k(S^1,\mathbb{C})$ as a section in $C^k(S^1,\gamma^*T\mathbb{C})$ via identification $C^k(S^1,\mathbb{C})\cong C^k(S^1,\gamma^*T\mathbb{C})$, thus $\xi_{\gamma(t)}\in T_{\gamma(t)}\mathbb{C}$ for all $t\in S^1$. We observe that at any point $(\gamma,\mathbf{t})\in\mathcal{F}^{-1}(\Delta_{\mathbb{C}^2})$, the composition of continuous linear maps 
\[
\begin{tikzcd}
	{T_{\gamma}C^k(S^1,\mathbb{C})\times T_{\mathbf{t}}(S^1)^4} \arrow[r, "{D_{(\gamma,\mathbf{t})}\mathcal{F}}"] & {T_{\mathcal{F}(\mathcal{\gamma,\mathbf{t}})}\mathbb{C}^4} \arrow[r, two heads] & {T_{\mathcal{F}(\mathcal{\gamma,\mathbf{t}})}\mathbb{C}^4/T_{\mathcal{F}(\mathcal{\gamma,\mathbf{t}})}\Delta_{\mathbb{C}^2}}
\end{tikzcd}
\]
between Banach spaces is clearly surjective, and the kernel is closed since it is the preimage of the closed subspace $T_{\mathcal{F}(\mathcal{\gamma,\mathbf{t}})}\Delta_{\mathbb{C}^2}\subset T_{\mathcal{F}(\mathcal{\gamma,\mathbf{t}})}\mathbb{C}^4$, thus it has complementary subspace in $T_{\gamma}C^k(S^1,\mathbb{C})\times T_{\mathbf{t}}(S^1)^4$. This implies that $\mathcal{F}\pitchfork\Delta_{\mathbb{C}^2}$ in $\mathbb{C}^4$. By implicit function theorem on Banach manifolds (c.f. \cite[Theorem A.3.3]{McDuff-Salamon2012}), we conclude that the real codimensional 4 submanifold $\Delta_{\mathbb{C}^2}\subset\mathbb{C}^2\times\mathbb{C}^2$ has preimage $\mathcal{F}^{-1}(\Delta_{\mathbb{C}^2})$ as a $C^k$-smooth Banach submanifold in $\operatorname{Emb}^k(S^1\times\mathbb{C})\times ((S^1)^4\setminus \Delta_{1234}^{(S^1)^4})$ of codimensional 4.

We next consider the projection map 
\[
\pi:\operatorname{Emb}^k(S^1,\mathbb{C})\times((S^1)^4\setminus \Delta_{1234}^{(S^1)^4})\to \operatorname{Emb}^k(S^1,\mathbb{C})
\] 
restricting to Banach submanifold $\mathcal{F}^{-1}(\Delta_{\mathbb{C}^2})$ which is $C^k$-smooth:
\begin{align*}
	\pi|_{\mathcal{F}^{-1}(\Delta_{\mathbb{C}^2})}:\mathcal{F}^{-1}(\Delta_{\mathbb{C}^2})&\to\operatorname{Emb}^k(S^1,\mathbb{C})\\
	(\gamma,\mathbf{t})&\mapsto\gamma
\end{align*}	
The linearized map reads
\begin{align*}
	D_{(\gamma,\mathbf{t})}\pi: T_{(\gamma,\mathbf{t})}\mathcal{F}^{-1}(\Delta_{\mathbb{C}^2})&\to T_{\gamma}C^k(S^1,\mathbb{C})\\
	(\xi,\mathbf{v})&\mapsto\xi
\end{align*}
where
\begin{align*}
T_{(\gamma,\mathbf{t})}\mathcal{F}^{-1}(\Delta_{\mathbb{C}^2})&=\{ (\xi,\mathbf{v})\in T_{\gamma}C^k(S^1,\mathbb{C})\times T_{\mathbf{t}}(S^1)^4 : \mathbf{V}=D_{(\gamma,\mathbf{t})}\mathcal{F}(\xi,\mathbf{v})\in T_{\mathcal{F}(\gamma,\mathbf{t})}\Delta_{\mathbb{C}^2} \} \\
&=  \{ (\xi,\mathbf{v})\in T_{\gamma}C^k(S^1,\mathbb{C})\times T_{\mathbf{t}}(S^1)^4 :  (V_1,V_2)=(V_3,V_4) , (\gamma,\mathbf{t})\in\mathcal{F}^{-1}(\Delta_{\mathbb{C}^2})  \}.
\end{align*}
A curve $\gamma$ being regular value of $\pi|_{\mathcal{F}^{-1}(\Delta_{\mathbb{C}^2})}$ is equivalent to the transversality condition $\mathcal{F}_\gamma\pitchfork\Delta_{\mathbb{C}^2}$ where $\mathcal{F}_\gamma:=\mathcal{F}(\gamma,\cdot):(S^1)^4\setminus \Delta_{1234}^{(S^1)^4}\to\mathbb{C}^4$, i.e.
\[
D_{(\gamma,\mathbf{t})}\mathcal{F}(\{0\}\times T_{\mathbf{t}}(S^1)^4)+ T_{\mathcal{F}(\gamma,\mathbf{t})}\Delta_{\mathbb{C}^2}=T_{\mathcal{F}(\gamma,\mathbf{t})}\mathbb{C}^4.
\]
This corresponds to the geometric condition in the statement:  $\operatorname{im}(\mathcal{F}_{\gamma})\pitchfork\Delta_{\mathbb{C}^2}$ means $L\times L_{\phi}\pitchfork\Delta_{\mathbb{C}^2}$ in $\mathbb{C}^2\times\mathbb{C}^2$ away from $(\gamma\times\{0\})\times(\gamma\times\{0\})$, which is equivalent to $L\pitchfork L_\phi$ in $\mathbb{C}^2$ away from $\gamma\times\{0\}$, and these transversal intersection points are isolated and hence finite due to compactness of the domain.

The linear map $D_{(\gamma,\mathbf{t})}\pi$ at every $(\gamma,\mathbf{t})\in\mathcal{F}^{-1}(\Delta_{\mathbb{C}^2})$ is a Fredholm operator: The kernel 
\begin{align*}
\ker(D_{(\gamma,\mathbf{t})}\pi)&= \left\{ (\xi,\mathbf{v})\in T_{\gamma}C^k(S^1,\mathbb{C})\times T_{\mathbf{t}}(S^1)^4 \left|  \substack{  \xi=0\\  D_{t_1}\gamma(v_1)=D_{t_3}\gamma(v_3)  \\ D_{t_2}\gamma(v_2)=D_{t_4}\gamma(v_4) }   \right\}\right.\\
&=T_{(\gamma,\mathbf{t})}\mathcal{F}^{-1}(\Delta_{\mathbb{C}^2})\cap (\{0\}\times T_{\mathbf{t}}(S^1)^4)\\
&\cong D_{(\gamma,\mathbf{t})}\mathcal{F}(\{0\}\times T_{\mathbf{t}}(S^1)^4)\cap T_{\mathcal{F}(\gamma,\mathbf{t})}\Delta_{\mathbb{C}^2}
\end{align*}
is finite dimensional. The cokernel is also finite dimensional by the following isomorphism
\[
\operatorname{coker}(D_{(\gamma,\mathbf{t})}\pi)=\frac{T_{\gamma}C^k(S^1,\mathbb{C})}{\operatorname{im}(D_{(\gamma,\mathbf{t})}\pi)}\cong\frac{T_{\mathcal{F}(\gamma,\mathbf{t})}\mathbb{C}^4}{D_{(\gamma,\mathbf{t})}\mathcal{F}(\{0\}\times T_{\mathbf{t}}(S^1)^4) +  T_{\mathcal{F}(\gamma,\mathbf{t})}\Delta_{\mathbb{C}^2} }
\]
due to the first isomorphism theorem:
\[
\begin{tikzcd}
	{T_{\gamma}C^k(S^1,\mathbb{C})\times T_{\mathbf{t}}(S^1)^4} \arrow[r, "{D_{(\gamma,\mathbf{t})}\mathcal{F}}", two heads] \arrow[d, "{D_{(\gamma,\mathbf{t})}\pi}"', two heads] \arrow[rrdd, two heads] & {T_{\mathcal{F}(\gamma,\mathbf{t})}\mathbb{C}^4} \arrow[r, "q_1", two heads] & {\frac{T_{\mathcal{F}(\gamma,\mathbf{t})}\mathbb{C}^4}{ D_{(\gamma,\mathbf{t})}\mathcal{F}(\{0\}\times T_{\mathbf{t}}(S^1)^4) + T_{\mathcal{F}(\gamma,\mathbf{t})}\Delta_{\mathbb{C}^2} }}                                  \\
	{T_{\gamma}C^k(S^1,\mathbb{C})} \arrow[d, "q_2"', two heads]                                                                                                                                       &                                                                              &                                                                                                                                                                                                                        \\
	{\frac{T_{\gamma}C^k(S^1,\mathbb{C})}{\operatorname{im}(D_{(\gamma,\mathbf{t})}\pi)}}                                                                                                                   &                                                                              & {\frac{T_{\gamma}C^k(S^1,\mathbb{C})\times T_{\mathbf{t}}(S^1)^4}{T_{(\gamma,\mathbf{t})}\mathcal{F}^{-1}(\Delta_{\mathbb{C}^2})+\{ 0\}\times T_{\mathbf{t}}(S^1)^4}} \arrow[ll, "\sim", dashed] \arrow[uu, "\sim"', dashed]
\end{tikzcd}
\]
since
\[
\ker(q_1\circ D_{(\gamma,\mathbf{t})}\mathcal{F})=T_{(\gamma,\mathbf{t})}\mathcal{F}^{-1}(\Delta_{\mathbb{C}^2})+\{0\}\times T_{\mathbf{t}}(S^1)^4=\ker(q_2\circ D_{(\gamma,\mathbf{t})}\pi ).
\]
Consequently, the set $\operatorname{im}(D_{(\gamma,\mathbf{t})}\pi)$ is closed in $T_{\gamma}C^k(S^1,\mathbb{C})$. Thus, the linear map $D_{(\gamma,\mathbf{t})}\pi$ is a Fredholm operator with index zero:
\begin{align*}
\operatorname{ind}(D_{(\gamma,\mathbf{t})}\pi)&=\dim\ker(D_{(\gamma,\mathbf{t})}\pi)-\dim\operatorname{coker}(D_{(\gamma,\mathbf{t})}\pi)\\
&=\dim(D_{(\gamma,\mathbf{t})}\mathcal{F}(\{0\}\times T_{\mathbf{t}}(S^1)^4)\cap T_{\mathcal{F}(\gamma,\mathbf{t})}\Delta_{\mathbb{C}^2})-\dim T_{\mathcal{F}(\gamma,\mathbf{t})}\mathbb{C}^4\\
&\qquad+\dim(D_{(\gamma,\mathbf{t})}\mathcal{F}(\{0\}\times T_{\mathbf{t}}(S^1)^4) +  T_{\mathcal{F}(\gamma,\mathbf{t})}\Delta_{\mathbb{C}^2} )\\
&=\dim(D_{(\gamma,\mathbf{t})}\mathcal{F}(\{0\}\times T_{\mathbf{t}}(S^1)^4)\cap T_{\mathcal{F}(\gamma,\mathbf{t})}\Delta_{\mathbb{C}^2})- \dim T_{\mathcal{F}(\gamma,\mathbf{t})}\mathbb{C}^4\\
&\qquad + \dim D_{(\gamma,\mathbf{t})}\mathcal{F}(\{0\}\times T_{\mathbf{t}}(S^1)^4)  + \dim T_{\mathcal{F}(\gamma,\mathbf{t})}\Delta_{\mathbb{C}^2}\\
&\qquad -  \dim(D_{(\gamma,\mathbf{t})}\mathcal{F}(\{0\}\times T_{\mathbf{t}}(S^1)^4)\cap T_{\mathcal{F}(\gamma,\mathbf{t})}\Delta_{\mathbb{C}^2})\\
&=\dim D_{(\gamma,\mathbf{t})}\mathcal{F}(\{0\}\times T_{\mathbf{t}}(S^1)^4)  + \dim T_{\mathcal{F}(\gamma,\mathbf{t})}\Delta_{\mathbb{C}^2}-\dim T_{\mathcal{F}(\gamma,\mathbf{t})}\mathbb{C}^4\\
&=4+4-8\\
&=0
\end{align*}
In conclusion, $\pi|_{\mathcal{F}^{-1}(\Delta_{\mathbb{C}^2})}$ is a Fredholm map of index $0$. By Sard-Smale theorem (c.f. \cite[Theorem A.5.1]{McDuff-Salamon2012}), the set of regular values of $\pi|_{\mathcal{F}^{-1}(\Delta_{\mathbb{C}^2})}$ is residual in $\operatorname{Emb}^k(S^1,\mathbb{C})$. Thus, we have proved that the set
\[
\operatorname{Emb}^k_\text{reg}(S^1,\mathbb{C}):=\left\{ \text{regular values of $\pi|_{\mathcal{F}^{-1}(\Delta_{\mathbb{C}^2})}$} \right\}
\]
is generic in $\operatorname{Emb}^k(S^1,\mathbb{C})$ with respect to $C^k$-topology.

Now for passing from $C^k$-genericity to $C^\infty$-genericity, we employ the argument in \cite[Lemma 6.2]{Staffa24} adapting to our setting here. Note that such type of argument is also known to Clifford Taubes earlier, see e.g. \cite[p.54-55]{McDuff-Salamon2012} and \cite[Appendix C.2, p.106]{Frauenfelder2003}. Since $\operatorname{Emb}^k_\text{reg}(S^1,\mathbb{C})\subseteq\operatorname{Emb}^k(S^1,\mathbb{C})$ is a residual subset with respect to $C^k$-topology for each fixed $k\geq 1$, and a residual set can be expressed as a countable intersection of open dense sets, c.f. \cite[Lemma 1.6.2 (vi)]{Buehler-Salamon2018}, we can therefore write
\[
\operatorname{Emb}^k_{\text{reg}}(S^1,\mathbb{C})=\bigcap_{n\in\mathbb{N}}\mathcal{V}^{k}_{n}
\]
where each $\mathcal{V}^{k}_{n}$ is open and dense in $\operatorname{Emb}^k(S^1,\mathbb{C})$ with the $C^k$-topology. For each integer $j\geq k$, we put
\[
\mathcal{V}^{j}_{n}:=\mathcal{V}^{k}_{n}\cap \operatorname{Emb}^{j}(S^1,\mathbb{C}),
\]
and we have
\begin{align*}
\bigcap_{n\in\mathbb{N}}\mathcal{V}^{j}_{n}&=\left( \bigcap_{n\in\mathbb{N}}\mathcal{V}^{k}_{n} \right)\cap\operatorname{Emb}^{j}(S^1,\mathbb{C}) \\
&=\operatorname{Emb}^k_{\text{reg}}(S^1,\mathbb{C}) \cap  \operatorname{Emb}^{j}(S^1,\mathbb{C}) \\
&=\operatorname{Emb}^{j}_{\text{reg}}(S^1,\mathbb{C}).
\end{align*}
Due to the fact that $\operatorname{Emb}^{j}_{\text{reg}}(S^1,\mathbb{C})$ is residual in $\operatorname{Emb}^{j}(S^1,\mathbb{C})$, by the Baire category theorem (c.f. \cite[Theorem 1.6.4]{Buehler-Salamon2018}), it is dense and thus each $\mathcal{V}^{j}_{n}$ is dense and open in $\operatorname{Emb}^{j}(S^1,\mathbb{C})$, since the $C^{j}$-topology is finer than the $C^{k}$-topology. We further put
\begin{align*}
\mathcal{V}^{\infty}_{n}&:=\mathcal{V}^{k}_{n}\cap \operatorname{Emb}^{\infty}(S^1,\mathbb{C})\\
&=\mathcal{V}^{k}_{n}\cap \left(\bigcap_{j\geq k}\operatorname{Emb}^{j}(S^1,\mathbb{C})  \right)\\
&=\bigcap_{j\geq k}\mathcal{V}^{j}_{n}.
\end{align*}
We claim that $\mathcal{V}^{\infty}_{n}$ is dense in $\operatorname{Emb}^{\infty}(S^1,\mathbb{C})$: Pick an arbitrary $\gamma_0\in\operatorname{Emb}^{\infty}(S^1,\mathbb{C})$, and an open neighbourhood $\mathcal{U}_0$ of it in $\operatorname{Emb}^{\infty}(S^1,\mathbb{C})$. Let $j\in\mathbb{N}_{\geq k}$ and $\delta>0$ such that 
\[
\left\{ \gamma\in\operatorname{Emb}^{\infty}(S^1,\mathbb{C}): d_{C^j}(\gamma,\gamma_0)<\delta \right\}\subseteq\mathcal{U}_0
\]
where $d_{C^j}$ is the metric inducing $C^j$-topology on $\operatorname{Emb}^{j}(S^1,\mathbb{C})$. Due to the density of $\mathcal{V}^{j}_{n}$ in $\operatorname{Emb}^{j}(S^1,\mathbb{C})$,
\[
\exists\gamma_1\in \mathcal{V}^{j}_{n} : d_{C^j}(\gamma_1,\gamma_0)<\frac{\delta}{2}.
\]
We know that $\operatorname{Emb}^{\infty}(S^1,\mathbb{C})$ is dense in $\operatorname{Emb}^{j}(S^1,\mathbb{C})$ with respect to $C^j$-topology, and $\mathcal{V}^{j}_{n}$ is open $\operatorname{Emb}^{j}(S^1,\mathbb{C})$, then 
\[
\exists \gamma_2\in \operatorname{Emb}^{\infty}(S^1,\mathbb{C})\cap \mathcal{V}^{j}_{n}: d_{C^j}(\gamma_2,\gamma_1)<\frac{\delta}{2}.
\]
By the triangle inequality of the metric $d_{C^j}$, we have
\[
\gamma_2\in \{ \gamma\in \operatorname{Emb}^{\infty}(S^1,\mathbb{C}): d_{C^j}(\gamma,\gamma_0)<\delta \}\cap \mathcal{V}^{\infty}_{n}\subseteq\mathcal{U}_0\cap \mathcal{V}^{\infty}_{n}.
\]
This proves the claim that $\mathcal{V}^{\infty}_{n}$ is dense in $\operatorname{Emb}^{\infty}(S^1,\mathbb{C})$. It is also open in $\operatorname{Emb}^{\infty}(S^1,\mathbb{C})$ since the $C^\infty$-topology is finer than the $C^k$-topology. Moreover, 
\begin{align*}
\bigcap_{n\in\mathbb{N}}\mathcal{V}^{\infty}_{n}&=\left( \bigcap_{n\in\mathbb{N}}\mathcal{V}^{k}_{n} \right) \cap \operatorname{Emb}^{\infty}(S^1,\mathbb{C})\\
&=  \operatorname{Emb}^{k}_{\text{reg}}(S^1,\mathbb{C}) \cap \left( \bigcap_{j\geq k} \operatorname{Emb}^{j}(S^1,\mathbb{C}) \right)\\
&= \bigcap_{j\geq k} \operatorname{Emb}^{j}_{\text{reg}}(S^1,\mathbb{C}) \\
&=\operatorname{Emb}^{\infty}_{\text{reg}}(S^1,\mathbb{C}).
\end{align*}
All together implies that $\operatorname{Emb}^{\infty}_{\text{reg}}(S^1,\mathbb{C})$ is residual in $\operatorname{Emb}^{\infty}(S^1,\mathbb{C})$ with respect to the $C^\infty$-topology.  As $\operatorname{Emb}^{\infty}_{\text{reg}}(S^1,\mathbb{C})$ is precisely the set of smooth Jordan curves with desired property in the statement, this completes the proof.
\end{proof}

\begin{rem}\label{rem_trans_pt}
One should note that this statement does not imply the existence of an intersection point away from $\gamma\times\{0\}$. As the ``empty intersection away from $\gamma\times\{0\}$'' also stays in the context of ``transversal intersection away from $\gamma\times\{0\}$''. In other words, the curves resulting in an empty intersection away from $\gamma\times\{0\}$ also belong to the residual set in the statement. Moreover, as the transversal intersection points away from $\gamma\times\{0\}$ are isolated hence giving rise to isolated distinct inscribed rectangles, the converse is not true in general: an isolated inscribed rectangle might come from an isolated tangential intersection point away from $\gamma\times\{0\}$, which is degenerate, and thus the curve in such scenario is not generic.
\end{rem}

\begin{exmp}\label{exa_ellipses}
One can easily tell that all ellipses except for the standard circle lie in such residual set in the statement of Proposition \ref{prop_gen_imm}. And one could produce more such generic curves from ellipses via small perturbation.
\end{exmp}

\section{On the existence of rectangular pegs}\label{Section 4}

We prove in this section Theorem \ref{existence_result} using the Lagrangian Floer homology introduced in Appendix \ref{App_Floer_homology}. To see how the Lagrangian tori fit into the framework, we first split the diffeomorphism \eqref{l_map} as the composition of the rescaling map by $\frac{1}{\sqrt{2}}$ and the symplectomorphism
\[
F:(z_1,z_2)\mapsto\left(\frac{z_1+z_2}{\sqrt{2}},\frac{z_1-z_2}{\sqrt{2}}\right)
\]
on $(\mathbb{C}^2,\omega_{\text{std}})$, i.e. $F^*\omega_{\text{std}}=\omega_{\text{std}}$. Namely, 
\[
l=\frac{1}{\sqrt{2}}\circ F=F\circ\frac{1}{\sqrt{2}}. 
\]
Therefore, $L=F(\frac{1}{\sqrt{2}}(\gamma\times\gamma))$ and $L_\phi=R_\phi F(\frac{1}{\sqrt{2}}(\gamma\times\gamma))$. We observe that $L$ and $L_\phi$ are monotone Lagrangian tori in $(\mathbb{C}^2,\omega_\text{std})$ with minimal Maslov number $N_L=N_{L_\phi}=2$. The description of monotone Lagrangian submanifolds is referred to \ref{App_monotone_Lag_Floer}.  This can be seen as follows. We assume without loss of generality that $\gamma$ bounds a domain $\overline{\gamma}$ with the same area as the unit disk $D^2(1)$, i.e. $\operatorname{Area}(\overline{\gamma})=\operatorname{Area}(D^2(1))=\pi$. The torus $\gamma\times\gamma$ is a monotone Lagrangian torus in $(\mathbb{C}^2,\omega_\text{std})$ with minimal Malsov index $2$ and monotonicity constant $\tau=\frac{\pi}{2}$. The verification of \eqref{monotonicity} follows analogously from the computation in \cite[\S 16.3, Example 16.3.6 (1)]{Oh2015} for the Clifford torus 
\[
S^1(1)\times S^1(1)=\{ (z_1,z_2)\in\mathbb{C}^2:|z_1|=|z_2|=1 \},
\] 
since $\overline{\gamma}$ is biholomorphic to $D^2(1)$ by the Riemann mapping theorem and the computation is topological. In particular, these properties are preserved by symplectomorphisms. According to \cite{Karlsson2013}, there is an area preserving isotopy (i.e. symplectic isotopy) on $\mathbb{C}$ that takes $\gamma$ to $D^2(1)$, further due to $H_c^1(\mathbb{R}^2;\mathbb{R})\cong\{0\}$ and by \cite{Polterovich_Rosen2014}, this can be turned into a Hamiltonian isotopy. Such a result is also shown in \cite[Lemma 2.3]{Lekili-Maydanskiy2014}. Therefore, $\gamma\times\gamma$ is Hamiltonian diffeomorphic to the Clifford torus $S^1(1)\times S^1(1)$, and $\frac{1}{\sqrt{2}}(\gamma\times\gamma)$ can be regarded as being Hamiltonian diffeomorphic to the Clifford torus $S^1(\frac{1}{\sqrt{2}})\times S^1(\frac{1}{\sqrt{2}})$ due to \cite[Excercise 12.3.6]{McDuff-Salamon2017}. Furthermore, notice that $F$ is actually given by a unitary matrix 
\[
\begin{pmatrix}
	\frac{1}{\sqrt{2}} & \frac{1}{\sqrt{2}} \\
	\frac{1}{\sqrt{2}} & -\frac{1}{\sqrt{2}}
\end{pmatrix}
\]
acting on $\mathbb{C}^2$, which can be written as $e^{iA}$ for some $2\times 2$ Hermitian matrix 
\[
A=
\begin{pmatrix}
	\frac{(-5+2\sqrt{2})\pi}{2(-2+\sqrt{2})} & \frac{\pi}{2\sqrt{2}} \\
	\frac{\pi}{2\sqrt{2}}  &   \frac{(5+2\sqrt{2})\pi}{2(2+\sqrt{2})}	
\end{pmatrix}.
\]
Thus, $F=e^{iA}$ can be regarded as the time-1 map of the Hamiltonian flow generated by the autonomous Hamiltonian function
\begin{equation}
	H_F=\frac{1}{2}\operatorname{Re}\langle\mathbf{z},A\mathbf{z}\rangle
\end{equation}
where $\mathbf{z}=(z_1,z_2)\in\mathbb{C}^2$ and $\langle\cdot,\cdot\rangle$ denotes the Hermitian inner product on $\mathbb{C}^2$. We regard $L_\phi:=R_\phi(L)$ as the image of  $L$ under a Hamiltonian diffeomorphism generated by the autonomous Hamiltonian function
\begin{equation}\label{rotation_Hamiltonian}
	H_\phi=\frac{1}{4}\phi|z_2|^2
\end{equation}
where $\phi\in(0,\frac{\pi}{2})$, whose Hamiltonian vector field is $\phi\partial_{\theta_2}$ in terms of the polar coordinate $z_2=r_2 e^{i\theta_2}$ generating the time-1 map of the Hamiltonian flow $\varphi_{H_\phi}=R_\phi$ as the rotation around $\mathbb{C}\times\{0\}$ by angle $\phi$. In conclusion, the two Lagrangian tori $L$ and $L_\phi$ are actually Hamiltonian isotopic to the Clifford torus $S^1(\frac{1}{\sqrt{2}})\times S^1(\frac{1}{\sqrt{2}})$, and their monotonicity constants are $\tau=\frac{\pi}{4}$. The fact that $L$ and $L_\phi$ have minimal Maslov number $N_L=N_{L_\phi}=2$ can also be deduced from the results \cite{Viterbo1990,Polterovich1991} that an embedded Lagrangian torus in $\mathbb{C}^2$ has minimal Maslov number 2, which can be regarded as a special case of Audin's conjecture, see e.g. \cite{Damian2012}.

We next determine a suitable absolute grading for the transveral and clean intersection points in $L\cap L_\phi$ via the Maslov--Morse index introduced in Appendix \ref{App_monotone_Lag_Floer} and Appendix \ref{App_MBFloer}. Due to \eqref{index_differ}, the grading differs by multiple of integers up to capping. Let $\ell_0\in\mathcal{P}(L,L_\phi)$ is a based path, and $\Gamma_{\ell_{0}}$ is defined by \eqref{Deck_grp}. Since $\mathcal{P}(L,L_\phi)$ is clearly path connected, we can choose without loss of generality that $\ell_0\subset\mathbb{C}\times\{0\}$, i.e. it lies in $\operatorname{Fix}(\varrho)$.

\begin{lem}\label{index_mod_2}
Let $x$ for the moment be any transverse or clean intersection point in $L\cap L_\phi$ equipped with two cappings $\overline{x}$ and $\widetilde{x}$, then we have
\[
\mu([x,\overline{x}])\equiv\mu([x,\widetilde{x}])\mod 2.
\]
\end{lem}

\begin{proof}
By \eqref{index_differ}, we know that
\[
\mu([x,\overline{x}])-\mu([x,\widetilde{x}])=I_{\mu}([-\widetilde{x}\# \overline{x}]),
\]
where $[-\widetilde{x}\# \overline{x}]\in \Gamma_{\ell_0}$. We only need to show that $I_{\mu}([-\widetilde{x}\# \overline{x}])$ takes value in $2\mathbb{Z}$. Let $u:S^1\times[0,1]\to\mathbb{C}^2$ be a representative map of the class $[-\widetilde{x}\# \overline{x}]\in\pi_2(\mathbb{C}^2,L\cup L_\phi)$ satisfying $u(\cdot,0)\in L$, $u(\cdot,1)\in L_\phi$, $u(0,\cdot)=u(1,\cdot)=x$, $u(\frac{1}{2},\cdot)=\ell_0$. And let $v:D^2\to \mathbb{C}^2$ be a disk map with $v|_{\partial D^2}=u|_{S^1\times\{1\}}$. Then $[v]\in \pi_2(\mathbb{C}^2,L_{\phi})$ and $[u\# v]\in\pi_2(\mathbb{C}^2,L)$. By the gluing property of the Maslov index, 
\[
\mu([u])+\mu([v])=\mu([u\# v]),
\]
and due to the minimal Maslov number $N_L=N_{L_\phi}=2$, we have
\[
\mu([u])=I_\mu([-\widetilde{x}\# \overline{x}])\in 2\mathbb{Z}.
\]
This proves the lemma.
\end{proof}
\noindent By \eqref{MB_grading} and \eqref{MB_index}, we fix a common capping for any point in the clean intersection $\gamma\times\{0\}$, suppose $f:\gamma\times\{0\}\to\mathbb{R}$ is a perfect Morse function, and let $x\in\operatorname{Crit}(f)$ be of Morse index $1$ while $y\in\operatorname{Crit}(f)$ has Morse index $0$, then due to Lemma \ref{index_mod_2} the absolute gradings for $x$ and $y$ can be given by
\begin{equation}\label{grading_x&y}
	|x|\equiv 1 \mod 2  ,\quad |y|\equiv 0 \mod 2.
\end{equation}
We can assign similarly the absolute grading $|p|:=\mu(p,\overline{p})$ for a transverse intersection point $p\in L\cap L_\phi\setminus(\gamma\times\{0\})$ with respect to a fixed cap $\overline{p}$ such that 
\begin{equation}\label{grading_p&q}
	|p|\equiv 0 \;\text{or} \; 1 \mod 2.
\end{equation} 
By Proposition \ref{prop_aut_of_point}, if such $p$ exists, there would be another transverse intersection point $\varrho(p)\in L\cap L_\phi\setminus(\gamma\times\{0\})$. It remains to determine the grading $\varrho(p)$.

\begin{lem}\label{index_invol_invar}
For any transverse intersection point $p\in L\cap L_\phi\setminus(\gamma\times\{0\})$, we have
\[
|\varrho(p)|:=\mu(\varrho(p),\overline{\varrho(p)})=\mu(p,\overline{p})=:|p|
\]	
in which we put $\overline{\varrho(p)}:=\varrho(\overline{p})$.
\end{lem}
\begin{proof}
This can be shown as follows. By \eqref{rel_abs_index}, we have
\[
\mu(p,\varrho(p);B)=\mu([\varrho(p),\varrho(\overline{p})])-\mu([p,\overline{p}])
\]
where $\varrho(\overline{p})=\overline{p}\# B$ with $B=[-\overline{p}\#\varrho(\overline{p})]\in\pi_2(p,\varrho(p))$. Similarly,  
\begin{align*}
	\mu(\varrho(p),\varrho(\varrho(p));\varrho_*(B))&=\mu(\varrho(p),p;-B)\\
	&=\mu([p,\overline{p}])-\mu([\varrho(p),\varrho(\overline{p})])
\end{align*}
where $\varrho_*(B)=[-\varrho(\overline{p})\#\overline{p}]=-B\in\pi_2(\varrho(p),p)$. Since $\varrho$ is a symplectomorphism and $L, L_\phi$ are both $\varrho$-invariant, by the invariance of the Malsov index of Lagrangian pair under symplectic action, c.f. \cite[Appendix C.1 (1)]{Oh2015}, we have
\[
\mu(\varrho(p),\varrho(\varrho(p));\varrho_*(B))=\mu(p,\varrho(p);B).
\]
This verifies the claim.
\end{proof}

The complex space $(\mathbb{C}^2,\omega_{\text{std}}=d\lambda_{\text{std}})$, where 
\[
\lambda_{\text{std}}=-\frac{i}{4}\sum_{i=1,2}(z_i dz_i+ \bar{z}_idz_i- z_id\bar{z}_i - \bar{z}_i d\bar{z}_i),
\] 
is a Liouville manifold, see Example \ref{Exa_Liou_mfd} (1). By Appendix \ref{App_Liouv_MfD}, we can choose a generic admissible $d\lambda_{\text{std}}$-compatible almost complex structure $J\in\mathcal{J}(\mathbb{C}^2,d\lambda_{\text{std}})$ so that the $\mathbb{F}_2$-coefficient self Lagrangian Floer homology $HF_*(L,L)$  is well-defined as a $\mathbb{Z}/2\mathbb{Z}$-graded module, whereas with respect to the Novikov ring coefficient \eqref{Nov_ring}, $HF_*(L,L;\Lambda)$ is a bounded $\mathbb{Z}$-graded $\Lambda$-module concentrated in degree $0$ and $1$. It is very easy to observe that the Lagrangian torus $L:= l(\gamma\times\gamma)$ in $\mathbb{C}^2$ can be evidently displaced by translation, e.g. translating in $z_1$-direction in $\mathbb{C}\times\mathbb{C}$:
\begin{equation}\label{translation}
	(z_1,z_2)=\mathbf{z} \mapsto \mathbf{z}+\mathbf{v}=(z_1+v_1,z_2)
\end{equation}
which is the time-1 map of Hamiltonian flow for the corresponding translation Hamiltonian function 
\begin{equation}\label{translation_Hamiltonian}
	H_{\mathbf{v}}=\frac{1}{2}\operatorname{Re}\langle -i\mathbf{v},\mathbf{z}\rangle=-\frac{1}{2}\operatorname{Re}(iv_1\bar{z}_1)=\frac{1}{2}(b_1x_1-a_1y_1)
\end{equation}
in $(\mathbb{C}^2,\omega_\text{std})$, where  $\mathbf{v}=(v_1,0)\in\mathbb{C}\times\{0\}$ and $v_1=a_1+ib_1$ with $a_1,b_1\in\mathbb{R}$. The image of $L$ under translation along $\mathbf{v}$ is denoted by $L_{\mathbf{v}}$. When $\|\mathbf{v}\|=|v_1|$ is sufficiently large, then $L\cap L_{\mathbf{v}}=\varnothing$. Therefore, the Lagrangian Floer homology vanishes:
\begin{equation}\label{HF_vanish}
HF_*(L,L_{\mathbf{v}})\cong \{0\}
\end{equation}
When $\gamma$ lies in the generic class in Proposition \ref{prop_gen_imm}, the monotone Lagrangian tori $L$ and $L_\phi$ intersect cleanly. Choose a perfect Morse function  $f$ on $\gamma\times\{0\}$ whose critical points $x,y\in\operatorname{Cirt}(f)$, these lift to $\mathbb{Z}$ copies of Morse--Bott circles in $\widetilde{\mathcal{{P}}}(L,L_\phi)$ with Morse functions $\{f^\nu\}_{\nu\in\mathbb{Z}}$ replicating $f$, and the critical points are graded by \eqref{grading_x&y} depending on the caps which also determine the labeling $\nu\in\mathbb{Z}$. Then Morse--Bott Floer homology $HF_*(L,L_\phi,\{f^\nu\})$ respectively $HF_*(L,L_\phi,f;\Lambda)$ is well-defined as a  $\mathbb{Z}/2\mathbb{Z}$-graded $\mathbb{F}_2$-module respectively a $\mathbb{Z}$-graded $\Lambda$-module, c.f. Appendix \ref{App_MBFloer}. 

However, our primary goal here is to prove \eqref{a_priori_exis}, namely the existence of inscribed rectangles with fixed aspect angle $\phi$ in any smooth Jordan curve $\gamma$, we cannot assume $\gamma$ a priori lies in the generic class such that $L\pitchfork L_\phi$ away from $\gamma\times\{0\}$, and therefore $L$ and $L_\phi$ not necessary intersect cleanly in general. Our strategy is to give an indirect proof: assuming  \eqref{a_priori_exis} was false, then $L\cap L_\phi\setminus(\gamma\times\{0\})=\varnothing$, i.e. $L$ and $L_\phi$ would intersect only cleanly in $\gamma\times\{0\}$, and in that case $HF_*(L,L_\phi,\{f^\nu\})$ is well-defined as claimed before; under such a assumption, we could compute and expect a non-vanishing result $HF_* (L,L_\phi,\{f^\nu\},J_t)\ncong \{0\}$ in the Morse--Bott setting with some specific choice of data, e.g. the almost complex structures $\{J_t\}$, leading one to anticipate a contradiction to \eqref{HF_vanish} due to the invariance of Floer homology. Such a contradiction would play a similar role as the one in the proof of \cite[Proposition 1.5]{Schmaeschke2016}.  

\begin{figure}
	\[
	\begin{tikzpicture}[x=0.75pt,y=0.75pt,yscale=-1,xscale=1]
		
		\draw  [color={rgb, 255:red, 0; green, 120; blue, 255 }  ,draw opacity=1 ] (330.85,130) .. controls (341.9,129.95) and (350.92,145.59) .. (351,164.92) .. controls (351.08,184.25) and (342.19,199.95) .. (331.15,200) .. controls (320.1,200.05) and (311.08,184.41) .. (311,165.08) .. controls (310.92,145.75) and (319.81,130.05) .. (330.85,130) -- cycle ;
		\draw  [fill={rgb, 255:red, 0; green, 0; blue, 0 }  ,fill opacity=1 ] (328.5,130.5) .. controls (328.5,129.12) and (329.62,128) .. (331,128) .. controls (332.38,128) and (333.5,129.12) .. (333.5,130.5) .. controls (333.5,131.88) and (332.38,133) .. (331,133) .. controls (329.62,133) and (328.5,131.88) .. (328.5,130.5) -- cycle ;
		\draw  [fill={rgb, 255:red, 0; green, 0; blue, 0 }  ,fill opacity=1 ] (328.5,199.5) .. controls (328.5,198.12) and (329.62,197) .. (331,197) .. controls (332.38,197) and (333.5,198.12) .. (333.5,199.5) .. controls (333.5,200.88) and (332.38,202) .. (331,202) .. controls (329.62,202) and (328.5,200.88) .. (328.5,199.5) -- cycle ;
		\draw  [color={rgb, 255:red, 0; green, 119; blue, 255 }  ,draw opacity=1 ] (330.85,234) .. controls (341.9,233.95) and (350.92,249.59) .. (351,268.92) .. controls (351.08,288.25) and (342.19,303.95) .. (331.15,304) .. controls (320.1,304.05) and (311.08,288.41) .. (311,269.08) .. controls (310.92,249.75) and (319.81,234.05) .. (330.85,234) -- cycle ;
		\draw  [fill={rgb, 255:red, 0; green, 0; blue, 0 }  ,fill opacity=1 ] (328.5,234.5) .. controls (328.5,233.12) and (329.62,232) .. (331,232) .. controls (332.38,232) and (333.5,233.12) .. (333.5,234.5) .. controls (333.5,235.88) and (332.38,237) .. (331,237) .. controls (329.62,237) and (328.5,235.88) .. (328.5,234.5) -- cycle ;
		\draw  [fill={rgb, 255:red, 0; green, 0; blue, 0 }  ,fill opacity=1 ] (328.5,303.5) .. controls (328.5,302.12) and (329.62,301) .. (331,301) .. controls (332.38,301) and (333.5,302.12) .. (333.5,303.5) .. controls (333.5,304.88) and (332.38,306) .. (331,306) .. controls (329.62,306) and (328.5,304.88) .. (328.5,303.5) -- cycle ;
		\draw  [color={rgb, 255:red, 0; green, 119; blue, 255 }  ,draw opacity=1 ] (330.85,25) .. controls (341.9,24.95) and (350.92,40.59) .. (351,59.92) .. controls (351.08,79.25) and (342.19,94.95) .. (331.15,95) .. controls (320.1,95.05) and (311.08,79.41) .. (311,60.08) .. controls (310.92,40.75) and (319.81,25.05) .. (330.85,25) -- cycle ;
		\draw  [fill={rgb, 255:red, 0; green, 0; blue, 0 }  ,fill opacity=1 ] (328.5,25.5) .. controls (328.5,24.12) and (329.62,23) .. (331,23) .. controls (332.38,23) and (333.5,24.12) .. (333.5,25.5) .. controls (333.5,26.88) and (332.38,28) .. (331,28) .. controls (329.62,28) and (328.5,26.88) .. (328.5,25.5) -- cycle ;
		\draw  [fill={rgb, 255:red, 0; green, 0; blue, 0 }  ,fill opacity=1 ] (328.5,94.5) .. controls (328.5,93.12) and (329.62,92) .. (331,92) .. controls (332.38,92) and (333.5,93.12) .. (333.5,94.5) .. controls (333.5,95.88) and (332.38,97) .. (331,97) .. controls (329.62,97) and (328.5,95.88) .. (328.5,94.5) -- cycle ;
		\draw [color={rgb, 255:red, 0; green, 116; blue, 255 }  ,draw opacity=1 ]   (311,161.08) -- (311,170.08) ;
		\draw [shift={(311,172.08)}, rotate = 270] [color={rgb, 255:red, 0; green, 116; blue, 255 }  ,draw opacity=1 ][line width=0.75]    (10.93,-3.29) .. controls (6.95,-1.4) and (3.31,-0.3) .. (0,0) .. controls (3.31,0.3) and (6.95,1.4) .. (10.93,3.29)   ;
		\draw [color={rgb, 255:red, 0; green, 118; blue, 255 }  ,draw opacity=1 ]   (351,161.08) -- (351,170.08) ;
		\draw [shift={(351,172.08)}, rotate = 270] [color={rgb, 255:red, 0; green, 118; blue, 255 }  ,draw opacity=1 ][line width=0.75]    (10.93,-3.29) .. controls (6.95,-1.4) and (3.31,-0.3) .. (0,0) .. controls (3.31,0.3) and (6.95,1.4) .. (10.93,3.29)   ;
		\draw [color={rgb, 255:red, 0; green, 114; blue, 255 }  ,draw opacity=1 ]   (311,55.08) -- (311,64.08) ;
		\draw [shift={(311,66.08)}, rotate = 270] [color={rgb, 255:red, 0; green, 114; blue, 255 }  ,draw opacity=1 ][line width=0.75]    (10.93,-3.29) .. controls (6.95,-1.4) and (3.31,-0.3) .. (0,0) .. controls (3.31,0.3) and (6.95,1.4) .. (10.93,3.29)   ;
		\draw [color={rgb, 255:red, 0; green, 116; blue, 255 }  ,draw opacity=1 ]   (351,55.08) -- (351,64.08) ;
		\draw [shift={(351,66.08)}, rotate = 270] [color={rgb, 255:red, 0; green, 116; blue, 255 }  ,draw opacity=1 ][line width=0.75]    (10.93,-3.29) .. controls (6.95,-1.4) and (3.31,-0.3) .. (0,0) .. controls (3.31,0.3) and (6.95,1.4) .. (10.93,3.29)   ;
		\draw [color={rgb, 255:red, 0; green, 116; blue, 255 }  ,draw opacity=1 ]   (311,264.08) -- (311,273.08) ;
		\draw [shift={(311,275.08)}, rotate = 270] [color={rgb, 255:red, 0; green, 116; blue, 255 }  ,draw opacity=1 ][line width=0.75]    (10.93,-3.29) .. controls (6.95,-1.4) and (3.31,-0.3) .. (0,0) .. controls (3.31,0.3) and (6.95,1.4) .. (10.93,3.29)   ;
		\draw [color={rgb, 255:red, 0; green, 115; blue, 255 }  ,draw opacity=1 ]   (351,264.08) -- (351,273.08) ;
		\draw [shift={(351,275.08)}, rotate = 270] [color={rgb, 255:red, 0; green, 115; blue, 255 }  ,draw opacity=1 ][line width=0.75]    (10.93,-3.29) .. controls (6.95,-1.4) and (3.31,-0.3) .. (0,0) .. controls (3.31,0.3) and (6.95,1.4) .. (10.93,3.29)   ;
		\draw [color={rgb, 255:red, 253; green, 0; blue, 30 }  ,draw opacity=1 ]   (331,94.5) .. controls (324,98) and (323,124) .. (331,130.5) ;
		\draw [shift={(325.81,118.53)}, rotate = 265.8] [color={rgb, 255:red, 253; green, 0; blue, 30 }  ,draw opacity=1 ][line width=0.75]    (10.93,-3.29) .. controls (6.95,-1.4) and (3.31,-0.3) .. (0,0) .. controls (3.31,0.3) and (6.95,1.4) .. (10.93,3.29)   ;
		\draw [color={rgb, 255:red, 255; green, 4; blue, 4 }  ,draw opacity=1 ]   (331.15,95) .. controls (339,99) and (339,122) .. (331.15,131) ;
		\draw [shift={(336.36,118.7)}, rotate = 276.33] [color={rgb, 255:red, 255; green, 4; blue, 4 }  ,draw opacity=1 ][line width=0.75]    (10.93,-3.29) .. controls (6.95,-1.4) and (3.31,-0.3) .. (0,0) .. controls (3.31,0.3) and (6.95,1.4) .. (10.93,3.29)   ;
		\draw [color={rgb, 255:red, 253; green, 0; blue, 30 }  ,draw opacity=1 ]   (331,198.5) .. controls (324,202) and (323,228) .. (331,234.5) ;
		\draw [shift={(325.81,222.53)}, rotate = 265.8] [color={rgb, 255:red, 253; green, 0; blue, 30 }  ,draw opacity=1 ][line width=0.75]    (10.93,-3.29) .. controls (6.95,-1.4) and (3.31,-0.3) .. (0,0) .. controls (3.31,0.3) and (6.95,1.4) .. (10.93,3.29)   ;
		\draw [color={rgb, 255:red, 255; green, 4; blue, 4 }  ,draw opacity=1 ]   (331.15,199) .. controls (339,203) and (339,226) .. (331.15,235) ;
		\draw [shift={(336.36,222.7)}, rotate = 276.33] [color={rgb, 255:red, 255; green, 4; blue, 4 }  ,draw opacity=1 ][line width=0.75]    (10.93,-3.29) .. controls (6.95,-1.4) and (3.31,-0.3) .. (0,0) .. controls (3.31,0.3) and (6.95,1.4) .. (10.93,3.29)   ;
		\draw [color={rgb, 255:red, 253; green, 0; blue, 30 }  ,draw opacity=1 ]   (331,-10.5) .. controls (324,-7) and (323,19) .. (331,25.5) ;
		\draw [shift={(325.81,13.53)}, rotate = 265.8] [color={rgb, 255:red, 253; green, 0; blue, 30 }  ,draw opacity=1 ][line width=0.75]    (10.93,-3.29) .. controls (6.95,-1.4) and (3.31,-0.3) .. (0,0) .. controls (3.31,0.3) and (6.95,1.4) .. (10.93,3.29)   ;
		\draw [color={rgb, 255:red, 255; green, 4; blue, 4 }  ,draw opacity=1 ]   (331.15,-10) .. controls (339,-6) and (339,17) .. (331.15,26) ;
		\draw [shift={(336.36,13.7)}, rotate = 276.33] [color={rgb, 255:red, 255; green, 4; blue, 4 }  ,draw opacity=1 ][line width=0.75]    (10.93,-3.29) .. controls (6.95,-1.4) and (3.31,-0.3) .. (0,0) .. controls (3.31,0.3) and (6.95,1.4) .. (10.93,3.29)   ;
		\draw [color={rgb, 255:red, 253; green, 0; blue, 30 }  ,draw opacity=1 ]   (332,304.5) .. controls (325,308) and (324,334) .. (332,340.5) ;
		\draw [shift={(326.81,328.53)}, rotate = 265.8] [color={rgb, 255:red, 253; green, 0; blue, 30 }  ,draw opacity=1 ][line width=0.75]    (10.93,-3.29) .. controls (6.95,-1.4) and (3.31,-0.3) .. (0,0) .. controls (3.31,0.3) and (6.95,1.4) .. (10.93,3.29)   ;
		\draw [color={rgb, 255:red, 255; green, 4; blue, 4 }  ,draw opacity=1 ]   (332.15,305) .. controls (340,309) and (340,332) .. (332.15,341) ;
		\draw [shift={(337.36,328.7)}, rotate = 276.33] [color={rgb, 255:red, 255; green, 4; blue, 4 }  ,draw opacity=1 ][line width=0.75]    (10.93,-3.29) .. controls (6.95,-1.4) and (3.31,-0.3) .. (0,0) .. controls (3.31,0.3) and (6.95,1.4) .. (10.93,3.29)   ;

	\end{tikzpicture}
	\]
	\caption{Morse--Bott Floer complex $(CF_*(L,L_\phi,f,J_t^{\varrho}),\partial)$.}\label{Fig_6}
\end{figure}

\begin{proof}[Proof of Theorem \ref{existence_result}]
Assume $L\cap L_\phi\setminus(\gamma\times\{0\})=\varnothing$ provided arbitrary smooth Jordan curve $\gamma$ and fixed aspect angle $\phi\in(0,\frac{\pi}{2}]$. A convenient choice of the $d\lambda_{\text{std}}$-compatible almost complex structure $J$ for computation is the $\varrho$-invariant ones, i.e. 
\[
\varrho^*J=\varrho_*\circ J\circ\varrho_*=J.
\] 
 E.g. the standard almost complex structure
\[
J_0=
\begin{pmatrix}
	i & 0\\
	0 & i
\end{pmatrix}
\]
which is clearly $\varrho$-invariant and $d\lambda_{\text{std}}$-compatible. The set of such almost complex structures is denoted by $\mathcal{J}^{\mathbb{Z}/2\mathbb{Z}}(\mathbb{C}^2,d\lambda_{\text{std}})$, and it is a non-empty and contractible submanifold of the space $\mathcal{J}(\mathbb{C}^2,d\lambda_{\text{std}})$. This can be shown by following the proof of \cite[Proposition 1.1]{Welschinger2005} analogously. Transversality theorems mentioned in Appendix \ref{App_Floer_homology} for the Lagrangian (Morse--Bott) Floer homology hold similarly with generic choices of admissible $\{J_t^{\varrho}\}_{t\in[0,1]}\subset\mathcal{J}^{\mathbb{Z}/2\mathbb{Z}}(\mathbb{C}^2,d\lambda_{\text{std}})$, and hence the Morse--Bott Floer homology $HF_*(L,L_\phi,\{f^\nu\},J_t^{\varrho})$ is well-defined and satisfies the invariance property, i.e. Theorem \ref{MB_homology_thm}. We refer to \cite[\S 5c]{Khovanov-Seidel2002} for the core technical arguments that we can adapt and utilize for such equivariant transversality. Any $J_t^{\varrho}$-holomorphic strip $u$ with Lagrangian boundary condition on $L$ and $L_\phi$ which is not entirely lying in $\mathbb{C}\times\{0\}$ (i.e. $\operatorname{im}(u)\not\subset\operatorname{Fix}(\varrho)$) will give rise to another one by applying the involution $\varrho$, i.e. $u_{\varrho}:=\varrho\circ u$ solves the non-linear Cauchy-Riemann equation: 
\begin{align*}
	\bar{\partial}_{J_t^{\varrho}}u_\varrho&:=\frac{1}{2}(d(\varrho\circ u)+J_t^{\varrho}\circ d(\varrho\circ u)\circ i)\\
	&=\frac{1}{2}(\varrho_*\circ d u+J_t^{\varrho}\circ \varrho_*\circ du\circ i)\\
	&=\frac{1}{2}(\varrho_*\circ d u+\varrho_*\circ J_t^{\varrho}\circ du\circ i)\\
	&=\varrho_*\circ \bar{\partial}_{J_t^{\varrho}} u\\
	&=0
\end{align*}
Then one would get zero contribution from such pair of holomorphic strips in the boundary operator 
\[
\partial:CF_*(L,L_\phi,\{f^\nu\},J_t^{\varrho})\to CF_{*-1}(L,L_\phi,\{f^\nu\},J_t^{\varrho}).
\]
In principle, there might be a $J_t^{\varrho}$-holomorphic strip $u_0$ filling the circle $\gamma\times\{0\}$ in $\mathbb{C}\times\{0\}$ (i.e. $\operatorname{im}(u_0)\subset\operatorname{Fix}(\varrho)$ and $\partial\operatorname{im}(u_0)=\gamma\times\{0\}$) with asymptotics $x,y\in\operatorname{Crit}(f)$. If such $u_0$ is counted in the boundary operator, then it would yield non-zero contribution. However, this is impossible for two reasons. The first crucial reason is due to the equivariant transversality result in \cite[\S 5c]{Khovanov-Seidel2002} that such $J_t^{\varrho}$-holomorphic strip $u_0$ can never be regular for a generic choice of admissible $\{J_t^{\varrho}\}_{t\in[0,1]}\subset\mathcal{J}^{\mathbb{Z}/2\mathbb{Z}}(\mathbb{C}^2,d\lambda_{\text{std}})$ since $\operatorname{im}(u_0)\subset\operatorname{Fix}(\varrho)$. The second reason is that such $u_0$ has Maslov--Viterbo index 4 since it fills the circle $\gamma\times\{0\}$ in $L$ (or $L_\phi$), which is homotopic to the product of two fundamental generators in $\pi_1(L)\cong\pi_2(\mathbb{C}^2,L)$. Despite such a straightforward computation (also carried out similarly in \cite[p.934]{Greene-Lobb2023}), this can also be deduced from the fact that $u_0$ has symplectic area $\pi$ and the monotonicity constant of $L$ or $L_\phi$ is $\frac{\pi}{4}$. Any $J_t$-holomorphic strip contributing to the boundary operator has Maslov--Viterbo index $2$ by the dimension formula \eqref{MB_dimension} and the grading \eqref{MB_grading}, and it has symplectic area $\frac{\pi}{2}$. This can be seen directly from the difference of critical values of the action functional: the critical values of all points in the same Morse--Bott component are equal, the difference of the action functional between two Morse--Bott components is computed as
\begin{align*}
\mathcal{A}([x,\overline{x}])-\mathcal{A}([x,\widetilde{x}])&=-\int_{[0,1]\times[0,1]}\overline{x}{}^*\omega_{\text{std}}+\int_{[0,1]\times[0,1]}\widetilde{x}{}^*\omega_{\text{std}}\\
&=-\int_{S^1\times[0,1]}(-\widetilde{x}\#\overline{x})^*\omega_{\text{std}}\\
&=-I_{\omega_{\text{std}}}([-\widetilde{x}\#\overline{x}])\\
&=-\frac{\pi}{4} I_\mu([-\widetilde{x}\#\overline{x}])\in \frac{\pi}{2}\mathbb{Z}
\end{align*}
where $x$ is for the moment an arbitrary point in $\gamma\times\{0\}$,  $\overline{x}$ and $\widetilde{x}$ are various cappings for it. For the fourth equality, we use the monotonicity condition $I_{\omega_{\text{std}}}=\tau I_\mu$ with $\tau=\frac{\pi}{4}$. Thus, even such $J_t^{\varrho}$-holomorphic strip $u_0$ could exist, it is irrelevant for the boundary operator. Under the assumption that there is no intersection point away from $\gamma\times\{0\}$, the Morse--Bott Floer complex regarded as a $\mathbb{Z}/2\mathbb{Z}$-graded chain complex over $\mathbb{F}_2$ reads
\[
\begin{tikzcd}
	\cdots \arrow[r, "y\mapsto 2x"] & \mathbb{F}_2\langle x\rangle \arrow[r, "x\mapsto 2y"] & \mathbb{F}_2\langle y\rangle \arrow[r, "y\mapsto 2x"] & \mathbb{F}_2\langle x\rangle \arrow[r, "x\mapsto 2y"] & \mathbb{F}_2\langle y\rangle \arrow[r, "y\mapsto 2x"] & \cdots
\end{tikzcd}
\]
which is schematically depicted in Figure \ref{Fig_6}. Thus, the Morse--Bott Floer homology would turn out to be non-vanishing:
\[ 
HF_* (L,L_\phi,\{f^\nu\},J_t^{\varrho})\cong H_*(S^1;\mathbb{F}_2)\ncong\{0\}
\] 
in which the isomorphism is regarded as between the $\mathbb{Z}/2\mathbb{Z}$-graded $\mathbb{F}_2$-modules. Due to the invariance, i.e. Theorem \ref{MB_homology_thm}, and \eqref{HF_vanish}, we get a contradiction:
\[
\{0\}\ncong HF_* (L,L_\phi,\{f^\nu\},J_t^{\varrho})\cong HF_*(L,L_{\mathbf{v}},J_t)\cong\{0\}
\]
where $\{J_t\}_{t\in[0,1]}\subset\mathcal{J}(\mathbb{C}^2,d\lambda_{\text{std}})$ is arbitrary generic choice of admissible ones (not necessarily $\varrho$-invariant). Hence, the assumption is false, and \eqref{a_priori_exis} is true.
\end{proof}
\begin{rem}
Moreover, we anticipate that a $\mathbb{Z}/2\mathbb{Z}$-equivariant Lagrangian Floer homology for these $\mathbb{Z}/2\mathbb{Z}$-invariant Lagrangian tori, analogous to the ones developed in \cite{Seidel-Smith2010,Kim-Kim-Kwon2022}, could be established, with the $\mathbb{Z}/2\mathbb{Z}$-symmetry arising from the action of $\varrho$; and it should be vanishing, i.e.
\[
HF_*^{\mathbb{Z}/2\mathbb{Z}}(L,L_\phi,\{f^\nu\})\cong HF_*^{\mathbb{Z}/2\mathbb{Z}}(L,L_{\mathbf{v}})\cong\{0\},
\]
due to the aforementioned $\mathbb{Z}/2\mathbb{Z}$-invariant translational displacement \eqref{translation} generated by the $\mathbb{Z}/2\mathbb{Z}$-invariant Hamiltonian \eqref{translation_Hamiltonian}.
\end{rem}

\begin{figure}
	\[
	\begin{tikzpicture}[x=0.75pt,y=0.75pt,yscale=-1,xscale=1]
		
		\draw [color={rgb, 255:red, 255; green, 159; blue, 0 }  ,draw opacity=1 ]   (387.85,109) -- (318.69,169.96) ;
		\draw [shift={(317.19,171.28)}, rotate = 318.61] [color={rgb, 255:red, 255; green, 159; blue, 0 }  ,draw opacity=1 ][line width=0.75]    (10.93,-3.29) .. controls (6.95,-1.4) and (3.31,-0.3) .. (0,0) .. controls (3.31,0.3) and (6.95,1.4) .. (10.93,3.29)   ;
		\draw [color={rgb, 255:red, 255; green, 0; blue, 29 }  ,draw opacity=1 ]   (287.19,57.91) .. controls (267.69,81.4) and (275.58,120.09) .. (284.33,142.4) ;
		\draw [shift={(285,144.08)}, rotate = 247.75] [color={rgb, 255:red, 255; green, 0; blue, 29 }  ,draw opacity=1 ][line width=0.75]    (10.93,-3.29) .. controls (6.95,-1.4) and (3.31,-0.3) .. (0,0) .. controls (3.31,0.3) and (6.95,1.4) .. (10.93,3.29)   ;
		\draw   (304.85,213) .. controls (315.9,212.95) and (324.92,228.59) .. (325,247.92) .. controls (325.08,267.25) and (316.19,282.95) .. (305.15,283) .. controls (294.1,283.05) and (285.08,267.41) .. (285,248.08) .. controls (284.92,228.75) and (293.81,213.05) .. (304.85,213) -- cycle ;
		\draw   (304.85,4) .. controls (315.9,3.95) and (324.92,19.59) .. (325,38.92) .. controls (325.08,58.25) and (316.19,73.95) .. (305.15,74) .. controls (294.1,74.05) and (285.08,58.41) .. (285,39.08) .. controls (284.92,19.75) and (293.81,4.05) .. (304.85,4) -- cycle ;
		\draw [color={rgb, 255:red, 253; green, 157; blue, 0 }  ,draw opacity=1 ]   (304.85,109) -- (385.85,109) ;
		\draw [shift={(387.85,109)}, rotate = 180] [color={rgb, 255:red, 253; green, 157; blue, 0 }  ,draw opacity=1 ][line width=0.75]    (10.93,-3.29) .. controls (6.95,-1.4) and (3.31,-0.3) .. (0,0) .. controls (3.31,0.3) and (6.95,1.4) .. (10.93,3.29)   ;
		\draw [color={rgb, 255:red, 130; green, 250; blue, 1 }  ,draw opacity=1 ][line width=1.5]    (331.43,109.29) -- (362.19,109.28) ;
		\draw [color={rgb, 255:red, 130; green, 250; blue, 1 }  ,draw opacity=1 ][line width=1.5]    (341,150.33) -- (370.19,124.28) ;
		\draw [color={rgb, 255:red, 253; green, 0; blue, 31 }  ,draw opacity=1 ]   (304.85,109) .. controls (289.03,96.48) and (287.33,72.56) .. (287.2,59.83) ;
		\draw [shift={(287.19,57.91)}, rotate = 90] [color={rgb, 255:red, 253; green, 0; blue, 31 }  ,draw opacity=1 ][line width=0.75]    (10.93,-3.29) .. controls (6.95,-1.4) and (3.31,-0.3) .. (0,0) .. controls (3.31,0.3) and (6.95,1.4) .. (10.93,3.29)   ;
		\draw [color={rgb, 255:red, 253; green, 0; blue, 31 }  ,draw opacity=1 ]   (291.19,219.82) .. controls (280.69,202.72) and (281.13,178.24) .. (286.41,164.75) ;
		\draw [shift={(287.19,162.91)}, rotate = 114.78] [color={rgb, 255:red, 253; green, 0; blue, 31 }  ,draw opacity=1 ][line width=0.75]    (10.93,-3.29) .. controls (6.95,-1.4) and (3.31,-0.3) .. (0,0) .. controls (3.31,0.3) and (6.95,1.4) .. (10.93,3.29)   ;
		\draw [color={rgb, 255:red, 253; green, 158; blue, 0 }  ,draw opacity=1 ]   (304.85,109) .. controls (362.61,121.7) and (376.91,114.96) .. (325.76,153.64) ;
		\draw [shift={(324.19,154.82)}, rotate = 322.96] [color={rgb, 255:red, 253; green, 158; blue, 0 }  ,draw opacity=1 ][line width=0.75]    (10.93,-3.29) .. controls (6.95,-1.4) and (3.31,-0.3) .. (0,0) .. controls (3.31,0.3) and (6.95,1.4) .. (10.93,3.29)   ;
		\draw [color={rgb, 255:red, 132; green, 255; blue, 0 }  ,draw opacity=1 ][line width=1.5]    (334.19,114.82) .. controls (364.19,119.82) and (362.19,123.82) .. (342.19,140.82) ;
		\draw [color={rgb, 255:red, 255; green, 159; blue, 0 }  ,draw opacity=1 ]   (304.85,109) .. controls (343.21,117.6) and (344.17,128.28) .. (326.4,142.79) ;
		\draw [shift={(325,143.92)}, rotate = 321.82] [color={rgb, 255:red, 255; green, 159; blue, 0 }  ,draw opacity=1 ][line width=0.75]    (10.93,-3.29) .. controls (6.95,-1.4) and (3.31,-0.3) .. (0,0) .. controls (3.31,0.3) and (6.95,1.4) .. (10.93,3.29)   ;
		\draw [color={rgb, 255:red, 132; green, 255; blue, 0 }  ,draw opacity=1 ][line width=1.5]    (331.29,118.57) .. controls (338.38,123.32) and (338.29,126.57) .. (335.29,133.57) ;
		\draw  [color={rgb, 255:red, 155; green, 155; blue, 155 }  ,draw opacity=1 ][fill={rgb, 255:red, 155; green, 155; blue, 155 }  ,fill opacity=1 ] (386.85,109) .. controls (386.85,107.62) and (387.97,106.5) .. (389.35,106.5) .. controls (390.73,106.5) and (391.85,107.62) .. (391.85,109) .. controls (391.85,110.38) and (390.73,111.5) .. (389.35,111.5) .. controls (387.97,111.5) and (386.85,110.38) .. (386.85,109) -- cycle ;
		\draw   (304.85,109) .. controls (315.9,108.95) and (324.92,124.59) .. (325,143.92) .. controls (325.08,163.25) and (316.19,178.95) .. (305.15,179) .. controls (294.1,179.05) and (285.08,163.41) .. (285,144.08) .. controls (284.92,124.75) and (293.81,109.05) .. (304.85,109) -- cycle ;
		\draw [color={rgb, 255:red, 255; green, 2; blue, 2 }  ,draw opacity=1 ]   (302.85,109) .. controls (255.67,111.89) and (241.48,190.32) .. (289.71,218.97) ;
		\draw [shift={(291.19,219.82)}, rotate = 209.17] [color={rgb, 255:red, 255; green, 2; blue, 2 }  ,draw opacity=1 ][line width=0.75]    (10.93,-3.29) .. controls (6.95,-1.4) and (3.31,-0.3) .. (0,0) .. controls (3.31,0.3) and (6.95,1.4) .. (10.93,3.29)   ;
		\draw  [fill={rgb, 255:red, 0; green, 0; blue, 0 }  ,fill opacity=1 ] (302.85,109) .. controls (302.85,107.62) and (303.97,106.5) .. (305.35,106.5) .. controls (306.73,106.5) and (307.85,107.62) .. (307.85,109) .. controls (307.85,110.38) and (306.73,111.5) .. (305.35,111.5) .. controls (303.97,111.5) and (302.85,110.38) .. (302.85,109) -- cycle ;
		\draw [color={rgb, 255:red, 132; green, 255; blue, 0 }  ,draw opacity=1 ][line width=1.5]    (265.19,137.91) .. controls (258.33,154) and (256.19,174.91) .. (267.19,195.91) ;
		\draw [color={rgb, 255:red, 132; green, 255; blue, 0 }  ,draw opacity=1 ][line width=1.5]    (283.19,179.91) .. controls (283.19,195.91) and (282.19,187.91) .. (284.19,202.91) ;
		\draw [color={rgb, 255:red, 132; green, 255; blue, 0 }  ,draw opacity=1 ][line width=1.5]    (288.33,77) .. controls (292.47,92.09) and (291.33,89) .. (297.33,102) ;
		\draw [color={rgb, 255:red, 132; green, 255; blue, 0 }  ,draw opacity=1 ][line width=1.5]    (277,79.33) .. controls (275,89.33) and (275.33,101) .. (276.33,114) ;
		
		\draw (393,103.4) node [anchor=north west][inner sep=0.75pt]  [font=\scriptsize,color={rgb, 255:red, 155; green, 155; blue, 155 }  ,opacity=1 ]  {$p$};
		\draw (355.19,142.68) node [anchor=north west][inner sep=0.75pt]  [font=\tiny,color={rgb, 255:red, 133; green, 255; blue, 0 }  ,opacity=1 ]  {$H$};

	\end{tikzpicture}
	\]
	\caption{Breaking of the interpolated Floer trajectories in the ``stretching homotopy'' procedure. The green region in an interpolated Floer trajectory resembles the finite time window in $\mathbb{R}$ where the Hamiltonian term is switched on. When $\|H\|<\frac{\pi}{2}$, the interpolated Floer trajectories depicted by the orange arrows have to break at some $p\in L\cap L_\phi\setminus(\gamma\times\{0\})$.  When $\|H\|>\frac{\pi}{2}$, the interpolated Floer trajectories can possess enough energy to break at the adjacent Morse--Bott components, these are depicted by the red arrows.}\label{Fig_7}
\end{figure}

An alternative potential approach for showing $\eqref{a_priori_exis}$ involving Floer trajectories without requiring the well-establishment of Floer homology would be following the ``stretching homotopy'' technique in \cite{Albers-Frauenfelder2010}. We can pick a copy of $\gamma\times\{0\}$ in the Novikov covering $\widetilde{\mathcal{P}}(L,L_\phi)$ of the path space as the Morse--Bott critical submanifold of the symplectic action functional $\mathcal{A}$ on which the Floer trajectories start and end, playing the role of the contact hypersurface in \cite{Albers-Frauenfelder2010}. The Hamiltonian function $H$ to be interpolated in the Floer trajectories in the procedure could be chosen as such that $\varphi_H$ pushes $L_\phi$ away from $L$, i.e. $L\cap\varphi_H(L_\phi)$. The minimal action gap among the critical values is $\frac{\pi}{2}$ from previous computation. In order to employ the ``stretching homotopy'' technique to detect an intersection point $p\in L\cap L_\phi\setminus(\gamma\times\{0\})$, it is crucial that the Hofer norm (definition refereed originally to \cite{Hofer1990}) of the push-off Hamiltonian $H$ should be strictly less than the minimal action gap, i.e. $\|H\|<\frac{\pi}{2}$. This is similar to the assumption made for the perturbation Hamiltonian in \cite{Albers-Frauenfelder2010}. In that case, due to the local Gromov--Floer convergence in the finite action window $[a-\|H\|,a+\|H\|]$ where $a$ is the critical value of the chosen Morse--Bott circle, analogous to \cite[Theorem 2.9]{Albers-Frauenfelder2010}, the interpolated Floer trajectories would have to break at an additional intersection point $p$ away from $\gamma\times\{0\}$, owing to the fact that no breaking could happen when time parameter $s\to\infty$ since $L\cap \varphi_H(L_\phi)=\varnothing$. If the Hofer norm of the push-off Hamiltonian is larger than the minimal action gap, i.e. $\|H\|>\frac{\pi}{2}$, such an argument would fail because the interpolated Floer trajectories could gain enough energy and break at the adjacent Morse--Bott circles, then such a procedure cannot detect the existence of $p\in L\cap L_\phi\setminus(\gamma\times\{0\})$. See Figure \ref{Fig_7} for the illustration. This emphasizes the importance of the assumption $\|H\|<\frac{\pi}{2}$. Nevertheless, it is well-known that the displacement energy of the Clifford torus $S^1(r)\times S^1(r)\subset\mathbb{C}^2$ is $\pi r^2$, i.e. $e(S^1(r)\times S^1(r))=\pi r^2$. C.f. \cite[p.62]{Sikorav1990}, \cite[p.359, Remark]{Polterovich1993}, and \cite[Example 2.7, Remark 3.3]{Brendel2023} for arbitrary dimensional result. By the symplectic invariance of the Hofer norm, we have
\[
e(L_\phi)=e(L)=e\left(S^1\left(\frac{1}{\sqrt{2}}\right)\times S^1\left(\frac{1}{\sqrt{2}}\right)\right)=\frac{\pi}{2}.
\]
It seems very hopeless to find a Hamiltonian $H$ pushing $L_\phi$ away from $L$ with $\|H\|<\frac{\pi}{2}$ even though assuming $L\cap L_\phi\setminus(\gamma\times\{0\})=\varnothing$ a priori. At least, this cannot be achieved by applying Hamiltonian diffeomorphism which is compactly supported closely around $\gamma\times\{0\}$ for the following reasons. By \cite[Lemma 3.1]{Cieliebak-Ekholm-Latschev2010}, a neighborhood of clean intersection $\gamma\times\{0\}\subset L\cap L_\phi\subset(\mathbb{C}^2,\omega_{\text{std}})$ is symplectomorphic to a neighborhood of clean intersection $\gamma\times\{0\}\subset 0_{L}\cap \nu^*_{\gamma\times\{0\}}\subset (T^*L,-d\lambda_{\text{can}})$, and there is no (compactly supported) Hamiltonian push-off of a conormal bundle away from the zero section in cotangent bundle. In other words, the Hofer norm $\|H\|$ cannot be too small.

Very recently, we promptly noticed new preprints \cite{Greene-Lobb2024a,Greene-Lobb2024b} by Greene and Lobb, in which they have ventured in a similar direction involving Floer homology albeit in a different geometric setup. The Lagrangian tori considered therein are
$\gamma\times\gamma$ and $\Phi_{h_t}^{1}(R_{\theta}(\gamma\times\gamma))$ where $R_{\theta}$ is the time-$\theta$ Hamiltonian flow of the Hamiltonian function 
\[
H=\frac{1}{4}|z_1-z_2|^2
\] 
which is the rotation around the diagonal $\Delta_{\mathbb{C}}\subset\mathbb{C}^2$ by angle $\theta$, and $\Phi_{h_t}^{1}$ is the time-1 flow of a generic non-degenerate time-dependent Hamiltonian function $h_t$ perturbing the Lagrangian torus $R_\theta(\gamma\times\gamma)$ so that $\gamma\times\gamma$ and $\Phi_{h_t}^1(R_{\theta}(\gamma\times\gamma))$ intersect transversely away from $\Delta_{\mathbb{C}}$, on which these two tori intersect cleanly and the clean intersection is $\Delta_{\gamma}:=\{(p,p)\in\gamma\times\gamma \}\cong\gamma$. We make a remark here that $\gamma\times\gamma$ and $R_{\theta}(\gamma\times\gamma)$ are invariant under the symplectic involution
\[
(z_1,z_2)\mapsto(z_2,z_1)
\]
whose fixed point set is $\Delta_{\mathbb{C}}$, and $(\gamma\times\gamma)\cap R_{\theta}(\gamma\times\gamma)\cap\Delta_{\mathbb{C}}=\Delta_{\gamma}$. Inspired from the construction of Heegaard Floer homology, Greene and Lobb manage to establish a version of Lagrangian Floer homology for these two Lagrangian tori as a $\mathbb{Z}$-graded $\mathbb{F}_2$-module whose boundary operator only counts the holomorphic strips avoiding the diagonal $\Delta_{\mathbb{C}}$. We expect that the Floer homology in \cite{Greene-Lobb2024a,Greene-Lobb2024b} can be established analogously in the geometric setup we have considered in this paper for the Lagrangian tori $L$ and $L_\phi$, i.e. the ones in \cite{Greene-Lobb2021}, where $\mathbb{C}\times\{0\}$ plays the role of $\Delta_{\mathbb{C}}$, and yields the same non-vanishing result as \cite[Theorem 2.23]{Greene-Lobb2024a}. A similar indirect proof can be made by using their version of Floer homology: for a general $\gamma$, assume that there was a priori no intersection point away from $\gamma\times\{0\}$, i.e. assuming \eqref{a_priori_exis} was false, then the transversality for $L\cap L_\phi$ away $\gamma\times\{0\}$ holds trivially, and the Floer homology in \cite{Greene-Lobb2024a,Greene-Lobb2024b} of $L$ and $L_\phi$ without Hamiltonian perturbation would be well-defined and vanishing since there is no generators in the Floer complex. On the other hand, however, the Floer homology of $L$ and $L_\phi$ with small Hamiltonian perturbation would be isomorphic to the relative singular homology of the pair $(L,\gamma\times\{0\})$, see the proof of \cite[Theorem 2.23]{Greene-Lobb2024a}. By the Hamiltonian invariance of the Floer homology, we would get a contradiction. This would prove \eqref{a_priori_exis} as well.

\section{The second rectangular peg}\label{Section 5}

A priori, the intersection point $p\in L\cap L_\phi\setminus(\gamma\times\{0\})$ corresponding to the first existing rectangular peg in Theorem \ref{existence_result} is not necessarily transversal for general $\gamma$, hence $\{p,\varrho(p)\}$ need not to be a pair of generators in any Floer complex. To test whether adding only one pair of generators $\{p,\varrho(p)\}$ other than $x^\nu,y^\nu\in\operatorname{Crit}(f^\nu)$ in the previous Morse--Bott Floer complex $CF_*(L,L_\phi,\{f^\nu\})$ in the well-defined case is enough to kill the Morse--Bott cirlce $\gamma\times\{0\}$ so that the Morse--Bott Floer homology $HF_*(L,L_\phi,\{f^\nu\})$ becomes vanishing as it oughts to be, we assume that $\gamma$ lies in the generic class in Proposition \ref{prop_gen_imm}, such that $L$ and $L_\phi$ intersects cleanly. It turns out that this is still not enough, and another occurring contradiction yields another pair of generators $\{q,\varrho(q)\}$, hence the existence of second rectangular peg.

\begin{proof}[Proof of Theorem \ref{main}]
Suppose $\gamma$ now lies in the generic class in Proposition \ref{prop_gen_imm} with any fixed $\phi\in(0,\frac{\pi}{2})$. Assume that there is only one $\varrho$-invariant pair of intersection points $\{p,\varrho(p)\}\subset L\cap L_\phi\setminus(\gamma\times\{0\})$ by Theorem \ref{existence_result}.  Due to Proposition \ref{prop_aut_of_point}, we aim to find another $\varrho$-invariant pair of intersection points in $L\cap L_\phi\setminus(\gamma\times\{0\})$ different from the $\{p,\varrho(p)\}$. Since $L$ is displaceable, and the Lagrangian Floer homology of $L$ vanishes by \eqref{HF_vanish}. Then after thoroughly examining all potential cascade trajectories scenarios (counting modulo 2) in the complex $CF_*(L,L_\phi,\{f^\nu\})$, we can conclude that 
\[
HF_*(L,L_\phi,\{f^\nu\})\ncong\{0\}.
\] 
By the invariance of Floer homology, i.e. Theorem \ref{MB_homology_thm}, we get a contradiction:
\[
\{0\}\ncong HF_*(L,L_\phi,\{f^\nu\})\cong HF_*(L,L_{\mathbf{v}})\cong\{0\}
\]
The tedious analysis of the Morse--Bott Floer complex can actually be bypassed by simply observing the Euler characteristic. By the gradings \eqref{grading_x&y} and \eqref{grading_p&q} and Lemma \ref{index_invol_invar}, the Euler characteristic of $HF_*(L,L_\phi,f;\Lambda)$ with Novikov ring coefficient \eqref{Nov_ring} is computed as
\begin{align*}
\chi(HF_*(L,L_\phi,f;\Lambda))&:=\sum_{i=0,1}(-1)^i\operatorname{rank}_\Lambda HF_i(L,L_\phi,f;\Lambda)\\
&=\sum_{i=0,1}(-1)^i\operatorname{rank}_\Lambda CF_i(L,L_\phi,f;\Lambda)\\
&=
\begin{cases}
	  2    &   \text{if}\; |p|=|\varrho(p)|=0   \\
	-2    &   \text{if}\;|p|=|\varrho(p)|=1
\end{cases}
\end{align*}
and meanwhile
\begin{align*}
	\chi(HF_*(L,L_{\mathbf{v}};\Lambda))&:=\sum_{i=0,1}(-1)^i\operatorname{rank}_\Lambda HF_i(L,L_{\mathbf{v}};\Lambda)\\
	&=\sum_{i=0,1}(-1)^i\operatorname{rank}_\Lambda CF_i(L,L_{\mathbf{v}};\Lambda)\\
	&=0
\end{align*}
This leads to a contradiction as well. Thus, there must be another distinct pair $\{q,\varrho(q)\}$ of intersection points away from $\gamma\times\{0\}$ in $L\cap L_\phi$. See Figure \ref{Fig_8}.
\end{proof}

\begin{figure}
	\[
	\begin{tikzpicture}[x=0.75pt,y=0.75pt,yscale=-1,xscale=1]
		
		\draw   (232.85,25) .. controls (243.9,24.95) and (252.92,40.59) .. (253,59.92) .. controls (253.08,79.25) and (244.19,94.95) .. (233.15,95) .. controls (222.1,95.05) and (213.08,79.41) .. (213,60.08) .. controls (212.92,40.75) and (221.81,25.05) .. (232.85,25) -- cycle ;
		\draw  [fill={rgb, 255:red, 0; green, 0; blue, 0 }  ,fill opacity=1 ] (230.5,25.5) .. controls (230.5,24.12) and (231.62,23) .. (233,23) .. controls (234.38,23) and (235.5,24.12) .. (235.5,25.5) .. controls (235.5,26.88) and (234.38,28) .. (233,28) .. controls (231.62,28) and (230.5,26.88) .. (230.5,25.5) -- cycle ;
		\draw  [fill={rgb, 255:red, 0; green, 0; blue, 0 }  ,fill opacity=1 ] (230.5,94.5) .. controls (230.5,93.12) and (231.62,92) .. (233,92) .. controls (234.38,92) and (235.5,93.12) .. (235.5,94.5) .. controls (235.5,95.88) and (234.38,97) .. (233,97) .. controls (231.62,97) and (230.5,95.88) .. (230.5,94.5) -- cycle ;
		\draw  [color={rgb, 255:red, 0; green, 0; blue, 0 }  ,draw opacity=1 ][fill={rgb, 255:red, 0; green, 0; blue, 0 }  ,fill opacity=1 ] (169.5,24.5) .. controls (169.5,23.12) and (170.62,22) .. (172,22) .. controls (173.38,22) and (174.5,23.12) .. (174.5,24.5) .. controls (174.5,25.88) and (173.38,27) .. (172,27) .. controls (170.62,27) and (169.5,25.88) .. (169.5,24.5) -- cycle ;
		\draw  [color={rgb, 255:red, 0; green, 0; blue, 0 }  ,draw opacity=1 ][fill={rgb, 255:red, 0; green, 0; blue, 0 }  ,fill opacity=1 ] (292.5,24.5) .. controls (292.5,23.12) and (293.62,22) .. (295,22) .. controls (296.38,22) and (297.5,23.12) .. (297.5,24.5) .. controls (297.5,25.88) and (296.38,27) .. (295,27) .. controls (293.62,27) and (292.5,25.88) .. (292.5,24.5) -- cycle ;
		\draw  [color={rgb, 255:red, 155; green, 155; blue, 155 }  ,draw opacity=1 ][fill={rgb, 255:red, 155; green, 155; blue, 155 }  ,fill opacity=1 ] (168.5,92.5) .. controls (168.5,91.12) and (169.62,90) .. (171,90) .. controls (172.38,90) and (173.5,91.12) .. (173.5,92.5) .. controls (173.5,93.88) and (172.38,95) .. (171,95) .. controls (169.62,95) and (168.5,93.88) .. (168.5,92.5) -- cycle ;
		\draw  [color={rgb, 255:red, 155; green, 155; blue, 155 }  ,draw opacity=1 ][fill={rgb, 255:red, 155; green, 155; blue, 155 }  ,fill opacity=1 ] (291.5,92.5) .. controls (291.5,91.12) and (292.62,90) .. (294,90) .. controls (295.38,90) and (296.5,91.12) .. (296.5,92.5) .. controls (296.5,93.88) and (295.38,95) .. (294,95) .. controls (292.62,95) and (291.5,93.88) .. (291.5,92.5) -- cycle ;
		\draw   (433.85,25) .. controls (444.9,24.95) and (453.92,40.59) .. (454,59.92) .. controls (454.08,79.25) and (445.19,94.95) .. (434.15,95) .. controls (423.1,95.05) and (414.08,79.41) .. (414,60.08) .. controls (413.92,40.75) and (422.81,25.05) .. (433.85,25) -- cycle ;
		\draw  [fill={rgb, 255:red, 0; green, 0; blue, 0 }  ,fill opacity=1 ] (431.5,25.5) .. controls (431.5,24.12) and (432.62,23) .. (434,23) .. controls (435.38,23) and (436.5,24.12) .. (436.5,25.5) .. controls (436.5,26.88) and (435.38,28) .. (434,28) .. controls (432.62,28) and (431.5,26.88) .. (431.5,25.5) -- cycle ;
		\draw  [fill={rgb, 255:red, 0; green, 0; blue, 0 }  ,fill opacity=1 ] (431.5,94.5) .. controls (431.5,93.12) and (432.62,92) .. (434,92) .. controls (435.38,92) and (436.5,93.12) .. (436.5,94.5) .. controls (436.5,95.88) and (435.38,97) .. (434,97) .. controls (432.62,97) and (431.5,95.88) .. (431.5,94.5) -- cycle ;
		\draw  [color={rgb, 255:red, 155; green, 155; blue, 155 }  ,draw opacity=1 ][fill={rgb, 255:red, 155; green, 155; blue, 155 }  ,fill opacity=1 ] (370.5,24.5) .. controls (370.5,23.12) and (371.62,22) .. (373,22) .. controls (374.38,22) and (375.5,23.12) .. (375.5,24.5) .. controls (375.5,25.88) and (374.38,27) .. (373,27) .. controls (371.62,27) and (370.5,25.88) .. (370.5,24.5) -- cycle ;
		\draw  [color={rgb, 255:red, 155; green, 155; blue, 155 }  ,draw opacity=1 ][fill={rgb, 255:red, 155; green, 155; blue, 155 }  ,fill opacity=1 ] (493.5,24.5) .. controls (493.5,23.12) and (494.62,22) .. (496,22) .. controls (497.38,22) and (498.5,23.12) .. (498.5,24.5) .. controls (498.5,25.88) and (497.38,27) .. (496,27) .. controls (494.62,27) and (493.5,25.88) .. (493.5,24.5) -- cycle ;
		\draw  [color={rgb, 255:red, 0; green, 0; blue, 0 }  ,draw opacity=1 ][fill={rgb, 255:red, 0; green, 0; blue, 0 }  ,fill opacity=1 ] (369.5,92.5) .. controls (369.5,91.12) and (370.62,90) .. (372,90) .. controls (373.38,90) and (374.5,91.12) .. (374.5,92.5) .. controls (374.5,93.88) and (373.38,95) .. (372,95) .. controls (370.62,95) and (369.5,93.88) .. (369.5,92.5) -- cycle ;
		\draw  [color={rgb, 255:red, 0; green, 0; blue, 0 }  ,draw opacity=1 ][fill={rgb, 255:red, 0; green, 0; blue, 0 }  ,fill opacity=1 ] (492.5,92.5) .. controls (492.5,91.12) and (493.62,90) .. (495,90) .. controls (496.38,90) and (497.5,91.12) .. (497.5,92.5) .. controls (497.5,93.88) and (496.38,95) .. (495,95) .. controls (493.62,95) and (492.5,93.88) .. (492.5,92.5) -- cycle ;
		
		\draw (56,88) node [anchor=north west][inner sep=0.75pt]   [align=left] {{\scriptsize degree 0}};
		\draw (57,19) node [anchor=north west][inner sep=0.75pt]   [align=left] {{\scriptsize degree 1}};
		\draw (167,9.4) node [anchor=north west][inner sep=0.75pt]  [font=\scriptsize,color={rgb, 255:red, 0; green, 0; blue, 0 }  ,opacity=1 ]  {$p$};
		\draw (286,9.4) node [anchor=north west][inner sep=0.75pt]  [font=\scriptsize,color={rgb, 255:red, 0; green, 0; blue, 0 }  ,opacity=1 ]  {$\varrho ( p)$};
		\draw (230,12.4) node [anchor=north west][inner sep=0.75pt]  [font=\scriptsize,color={rgb, 255:red, 0; green, 0; blue, 0 }  ,opacity=1 ]  {$x$};
		\draw (228,97.4) node [anchor=north west][inner sep=0.75pt]  [font=\scriptsize,color={rgb, 255:red, 0; green, 0; blue, 0 }  ,opacity=1 ]  {$y$};
		\draw (165.5,94.9) node [anchor=north west][inner sep=0.75pt]  [font=\scriptsize,color={rgb, 255:red, 155; green, 155; blue, 155 }  ,opacity=1 ]  {$q$};
		\draw (286,97.4) node [anchor=north west][inner sep=0.75pt]  [font=\scriptsize,color={rgb, 255:red, 155; green, 155; blue, 155 }  ,opacity=1 ]  {$\varrho ( q)$};
		\draw (368.5,95.9) node [anchor=north west][inner sep=0.75pt]  [font=\scriptsize,color={rgb, 255:red, 0; green, 0; blue, 0 }  ,opacity=1 ]  {$p$};
		\draw (486,97.4) node [anchor=north west][inner sep=0.75pt]  [font=\scriptsize,color={rgb, 255:red, 0; green, 0; blue, 0 }  ,opacity=1 ]  {$\varrho ( p)$};
		\draw (431,12.4) node [anchor=north west][inner sep=0.75pt]  [font=\scriptsize,color={rgb, 255:red, 0; green, 0; blue, 0 }  ,opacity=1 ]  {$x$};
		\draw (429,97.4) node [anchor=north west][inner sep=0.75pt]  [font=\scriptsize,color={rgb, 255:red, 0; green, 0; blue, 0 }  ,opacity=1 ]  {$y$};
		\draw (369.5,8.9) node [anchor=north west][inner sep=0.75pt]  [font=\scriptsize,color={rgb, 255:red, 155; green, 155; blue, 155 }  ,opacity=1 ]  {$q$};
		\draw (490,9.4) node [anchor=north west][inner sep=0.75pt]  [font=\scriptsize,color={rgb, 255:red, 155; green, 155; blue, 155 }  ,opacity=1 ]  {$\varrho ( q)$};
		\draw (216,53.4) node [anchor=north west][inner sep=0.75pt]  [font=\tiny]  {$\gamma \times \{0\}$};
		\draw (418,53.4) node [anchor=north west][inner sep=0.75pt]  [font=\tiny]  {$\gamma \times \{0\}$};

	\end{tikzpicture}
	\]
	\caption{Geometric generators in the Morse--Bott Floer complex $CF_*(L,L_\phi,f;\Lambda)$ with gradings \eqref{grading_x&y} and \eqref{grading_p&q}.}\label{Fig_8}
\end{figure}

\begin{rem}
Theorem \ref{main} can also be proved in a purely differential topological way via the algebraic intersection numbers.\footnote{This was suggested by Kai Cieliebak to the author after his talk in the \textit{Oberseminar Differentialgeometrie}, Wintersemester 2023-2024, at the Universität Augsburg.} We present the proof as follows. Recall that the algebraic intersection number $X\bullet Y$ of two closed oriented embedded submanifolds $X,Y$ in an oriented smooth manifold $M$ such that $X\pitchfork Y$ is defined to be the signed count of the intersection points, see Appendix \ref{App_int_numb_intro}. It is well-known that the sign of each intersection point does not depend on the choice of the positive bases, and the number $X\bullet Y$ is well-defined as a homotopy invariant. Now take $X=L$ and $Y=L_\phi=R_\phi(L)$ in $M=\mathbb{C}^2$. For a curve $\gamma$ in the generic class, $L$ intersects with $L_\phi$ cleanly along $\gamma\times\{0\}$ and transversely away from that part. To put it in the entirely transverse scenario, we perturb the autonomous Hamiltonian function \eqref{rotation_Hamiltonian} by a small Morse function:
\[
H_\phi^\varepsilon=H_\phi+\varepsilon\beta f
\]
where $f$ is a perfect Morse function on $\gamma\times\{0\}$, and $\beta:\mathbb{C}^2\to\mathbb{R}$ is a bump function supported in a small tubular neighborhood $\mathcal{U}$ of $\gamma\times\{0\}$ containing no intersection points in $L\cap L_\phi\setminus(\gamma\times\{0\})$, and $\varepsilon>0$ is sufficiently small. The symplectic gradient vector field of $H_\phi^{\varepsilon}$ is given by
\[
-J_0\operatorname{grad}_{g_0}H_{\phi}- J_0\operatorname{grad}_{g_0}\varepsilon\beta f
\]
from which we can tell that the zeros of the vector field in $\mathcal{U}$ is exactly the two critical points of $f$ on $\gamma\times\{0\}$. The time-1-flow yields a diffeomorphism $\varphi_{H_\phi^\varepsilon}$ on $L$ whose image $L_\phi^\varepsilon:=\varphi_{H_\phi^\varepsilon}(L)$ coincides with $L_\phi$ outside of a small tubular neighborhood $\mathcal{U}$ of $\gamma\times\{0\}$, and $L\cap L_\phi^\varepsilon$ has now only two transversal intersection points $x_0,x_1$ in $\mathcal{U}$ and $L\cap L_\phi^\varepsilon\setminus\mathcal{U}=L\cap L_\phi\setminus(\gamma\times\{0\})$, namely $L\pitchfork L_\phi^\varepsilon$. There is an obvious homotopy between $L_\phi$ and $L_\phi^\varepsilon$, and due to the homotopy invariance of algebraic intersection number, we put
\[
L\bullet L_\phi=L\bullet L_\phi^\varepsilon
\]
which is clearly independent from the choices of $f$, $\varepsilon$, $\beta$, etc. Assume there is only one pair of intersection points $\{p,\varrho(p)\}$ in $L\cap L_\phi\setminus(\gamma\times\{0\})$, then according to the perturbation above, up to choice of orientation, we have
\begin{align*}
L\bullet L_\phi=L\bullet L_\phi^\varepsilon&=(\operatorname{sign}(x_0)+\operatorname{sign}(x_1))+(\operatorname{sign}(p)+\operatorname{sign}(\varrho(p)))\\
&=\pm\chi(\gamma\times\{0\})\pm(1+1)\\
&=\pm\chi(S^1)\pm 2\\
&=0\pm 2 \\
&=\pm 2.
\end{align*}
The first bracket sum is due to the following. First, we notice that the orientation signs of $x_0,x_1\in L\cap L_\phi^{\varepsilon}\subset\mathbb{C}^2$ coincide with the orientation signs of them considered as intersection points of sections $(\gamma\times\{0\})\cap J_0\operatorname{grad}_{g_0}\varepsilon \beta f$ in the vector bundle $J_0T(\gamma\times\{0\})$. This is, up to a uniform sign, equal to the orientation signs of $x_0$ and $x_1$ as intersection points of sections $(\gamma\times\{0\})\cap \operatorname{grad}_{g_0} f$ in $T(\gamma\times\{0\})$, which in fact are further equal to the indices of $x_0$ and $x_1$ as zeros of the  gradient vector field $\operatorname{grad}_{g_0}  f$ in $T(\gamma\times\{0\})$. Hence the sum of the orientation signs of $x_0,x_1\in L\cap L_\phi^{\varepsilon}$ coincides up to a uniform sign with the sum of indices of zeros $x_0$ and $x_1$ of $\operatorname{grad}_{g_0}f$ in $T(\gamma\times\{0\})$, which is the Euler characteristic by the Poincar\'e--Hopf theorem. More detail and general treatment of computing algebraic intersection numbers provided appearance of clean intersection is referred to Appendix \ref{Appendix_B}. The second bracket sum is due to the fact that $\varrho$ is an orientation preserving diffeomorphism, therefore $p$ and $\varrho(p)$ must have the same orientation signs. On the other hand, we know that $L$ can be displaced away from itself by some translational smooth isotopy, hence it is apparent that
\[
L\bullet L_\phi= \operatorname{sign}(\varnothing) =0.
\]
This leads to a contradiction! Hence, the assumption must be wrong, and there must be another pair of intersection point $\{q,\varrho(q)\}$ in $L\cap L_\phi\setminus(\gamma\times\{0\})$. In fact, the Lagrangian Floer homology can be thought of as a ``categorification'' of algebriac intersection number in the sense that 
\[
\chi(HF_*(L,L_\phi,f;\Lambda))=L\bullet L_\phi
\]
See Appendix \ref{App_Cat}. In this spirit, these two proofs are well consistent.
\end{rem}

\section{Generic doubling of cyclic quadrilateral pegs}\label{Section 6}

In \cite{Greene-Lobb2023}, Greene and Lobb  consider a different geometric setup to obtain two embedded Lagrangian tori in $\mathbb{C}^2$ with respect to $\omega_\text{std}$ from a smooth Jordan curve $\gamma\subset\mathbb{C}$ as following:
\begin{align*}
	T_1:=R_\phi(F_s(\gamma\times\gamma)), \qquad T_2:=F_t(\gamma\times\gamma),
\end{align*}
where $s,t\in(0,\frac{1}{2}]$. The map
\begin{align*}
	F_r: \mathbb{C}^2 &\to \mathbb{C}^2 \\
	(z_1,z_2) &\mapsto \left( (1-r) z_1 + r z_2  , \sqrt{r(1-r)}(z_1- z_2)  \right)
\end{align*}
with $r\in(0,\frac{1}{2}]$ is a $\mathbb{C}$-linear automorphism of $\mathbb{C}^2$, and  
\begin{align*}
	R_\phi:\mathbb{C}^2&\to\mathbb{C}^2 \\
	(z_1,z_2)&\mapsto(z_1,z_2 e^{i\phi})
\end{align*}
with $\phi\in(0,\pi)$. Note that the maps $F_s$ and $F_t$ do not preserve $\omega_{\text{std}}$, the two embedded Lagrangian tori $T_1$ and $T_2$ are not related by Hamiltonian diffeomorphism in general. In particular, $T_1$ and $T_2$ are in general no longer invariant under symplectic involution $\varrho$ unless $s=t=\frac{1}{2}$. Moreover, there is no apparent symmetry for $T_1$ and $T_2$ when $s\neq \frac{1}{2}$ or $t\neq\frac{1}{2}$. Nevertheless, they have minimal Maslov number $2$ by \cite{Viterbo1990,Polterovich1991}, and intersect again cleanly along $\gamma\times\{0\}$, c.f. \cite[Lemma 2]{Greene-Lobb2023}. 

\subsubsection*{Intersection points and pegs}
By \cite[Lemma 1]{Greene-Lobb2023}, an intersection point in $T_1\cap T_2\setminus(\gamma\times\{0\})$, which is a solution of the equation
\[
R_\phi(F_s(A,C))=F_t(B,D),\qquad A\neq C \; (\Leftrightarrow B\neq D),
\]
for $A,B,C,D\in\gamma\subset\mathbb{C}$, gives rise to an inscribed cyclic quadrilateral on the curve $\gamma$. See Figure \ref{Fig_9}.
\begin{figure}
\[
\begin{tikzpicture}[x=0.75pt,y=0.75pt,yscale=-1,xscale=1]
	
	\draw   (268,66) .. controls (268,31.76) and (295.76,4) .. (330,4) .. controls (364.24,4) and (392,31.76) .. (392,66) .. controls (392,100.24) and (364.24,128) .. (330,128) .. controls (295.76,128) and (268,100.24) .. (268,66) -- cycle ;
	\draw    (287,22) -- (370,19) ;
	\draw    (287,22) -- (278,100) ;
	\draw    (278,100.5) -- (377,106.5) ;
	\draw    (377,106.5) -- (370,19) ;
	\draw  [dash pattern={on 0.84pt off 2.51pt}]  (286.5,22) -- (377,106.5) ;
	\draw  [dash pattern={on 0.84pt off 2.51pt}]  (370,19) -- (278,100) ;
	
	\draw (374,102.4) node [anchor=north west][inner sep=0.75pt]    {$A$};
	\draw (264,97.4) node [anchor=north west][inner sep=0.75pt]    {$D$};
	\draw (317,50.4) node [anchor=north west][inner sep=0.75pt]    {$X$};
	\draw (372,6.4) node [anchor=north west][inner sep=0.75pt]    {$B$};
	\draw (272,9.4) node [anchor=north west][inner sep=0.75pt]    {$C$};
	\draw (335,48.4) node [anchor=north west][inner sep=0.75pt]    {$\phi $};

\end{tikzpicture}
\]
\caption{Inscribed cyclic quadrilateral with data $s=\frac{|AX|}{|AC|}$, $t=\frac{|BX|}{|BD|}$, $\phi=\angle AXB$, and vertices $(A,B,C,D)$ ordered counterclockwise.}\label{Fig_9}
\end{figure}
The triple $(s,t,\phi)$ determines the counterclockwise oriented cyclic quadrilateral uniquely up to similarity, unless $s=\frac{1}{2}$ or $t=\frac{1}{2}$ in which case $(s,t,\phi)$ and $(t,s,\pi-\phi)$ determine the same oriented cyclic quadrilateral up to similarity. When $s=t=\frac{1}{2}$, this recovers the geometric setup in \cite{Greene-Lobb2021} as described in \S\ref{Section 1}. 

The correspondence between intersection points in  $T_1\cap T_2\setminus(\gamma\times\{0\})$ and  cyclic quadrilaterals with data $(s,t,\phi)$ inscribed on $\gamma$ is described as follows.
\begin{itemize}
\item[(1)] If $s\neq\frac{1}{2}$ or $t\neq\frac{1}{2}$, the intersection points are bijective to the inscribed cyclic quadrilaterals. 

\item[(2)] If $s=t=\frac{1}{2}$, this is the case as discussed in \S\ref{Section 2}, Proposition \ref{prop_aut_of_point}.
\end{itemize}

\subsubsection*{Existence of pegs}
After performing a Lagrangian surgery \cite[p.933-934]{Greene-Lobb2023}, the resulting object $T_1\#_{\gamma\times\{0\}}T_2$ is an immersed Lagrangian torus with minimal Maslov number $4$, c.f. \cite[p.434]{Greene-Lobb2023}. The self-intersection points in $T_1\#_{\gamma\times\{0\}}T_2$ are one-to-one corresponding to intersection points in $T_1\cap T_2\setminus(\gamma\times\{0\})$. Thus, if there were no self-intersection points in $T_1\#_{\gamma\times\{0\}}T_2$, it would be an embedded Lagrangian torus in $\mathbb{C}^2$, whose minimal Maslov number must be 2 by \cite{Viterbo1990, Polterovich1991}. This contradiction shows the existence of inscribed cyclic quadrilateral with any $(s,t,\phi)$ inscribed on any smooth Jordan curve $\gamma$, c.f. \cite[p.931]{Greene-Lobb2023}.

\subsubsection*{Transversality} 
Mimicking the discussion in \S\ref{Section 3}, a similar geometric transversality result holds as follows.

\begin{prop}\label{prop_geo_trans_cyc_quad}
For any fixed data $s,t\in(0,\frac{1}{2}]$ and $\phi\in(0,\pi)$, there is a residual set of smooth Jordan curves such that for every $\gamma$ in the set, the corresponding smooth Lagrangian tori $T_1$ and $T_2$ constructed in above way have only finitely many transversal intersection points away from $\gamma\times\{0\}$. 
\end{prop}

\noindent We provide a sketch of proof which is completely similar to that of Proposition \ref{prop_gen_imm}. The only degenerate inscribed cyclic quadrilaterals, which needs to be circumvented, is $A=B=C=D$. These amount to the clean intersection points $\gamma\times\{0\}$.  A map $\mathcal{F}$ in analogy to that in the proof of Proposition \ref{prop_gen_imm} can be constructed as:
\begin{align*}
\mathcal{F}: \operatorname{Emb}^k(S^1,\mathbb{C})\times ( (S^1)^4\setminus \Delta_{1234}^{(S^1)^4} ) &\to \mathbb{C}^2\times\mathbb{C}^2\\
(\gamma,t_1,t_2,t_3,t_4) & \mapsto (R_\phi(F_s(\gamma(t_1)\times\gamma(t_2)))  , F_t(\gamma(t_3)\times\gamma(t_4)) )
\end{align*} 
The set $\mathcal{F}^{-1}(\Delta_{\mathbb{C}^2})$ can be shown to be a Banach submanifold similarly by implciit function theorem, and the projection 
\begin{align*}
	\pi|_{\mathcal{F}^{-1}(\Delta_{\mathbb{C}^2})}:\mathcal{F}^{-1}(\Delta_{\mathbb{C}^2})&\to\operatorname{Emb}^k(S^1,\mathbb{C})\\
	(\gamma,t_1,t_2,t_3,t_4)&\mapsto\gamma
\end{align*}
is a Fredholm map with Fredholm index $0$, whose regular values amount to the geometric transversal condition $T_1\times  T_2\pitchfork\Delta_{\mathbb{C}^2}$ in $\mathbb{C}^2\times\mathbb{C}^2$ away from $(\gamma\times\{0\})\times(\gamma\times\{0\})$, in other words $T_1\pitchfork T_2$ in $\mathbb{C}^2$ away from $\gamma\times\{0\}$. By Sard-Samle theorem, the $C^k$-genericity statement follows. Further via the same argument at the end of the proof of Proposition \ref{prop_gen_imm}, the $C^k$-genericity can be promoted to $C^\infty$-genericity.

\subsubsection*{Doubling}
In the current situation, the two embedded smooth Lagrangian tori $T_1$ and $T_2$ are not Hamiltonian isotopic to each other and even not monotone. Neither the Floer homology of monotone Lagrangian submanifolds introduced in Appendix \ref{App_Floer_homology} nor the one in \cite{Greene-Lobb2024a,Greene-Lobb2024b} is applicable. Nevertheless, the computational formula for algebraic intersection number of cleanly intersecting Lagrangian submanifolds, i.e. \eqref{lagcleanint}, is still valid. Suppose now $\gamma$ lies in the generic class from Proposition \ref{prop_geo_trans_cyc_quad}. Assume that there is only one intersection point in $T_1\cap T_2\setminus(\gamma\times\{0\})$ due to the existence result in \cite{Greene-Lobb2023}, then by \eqref{lagcleanint}, we have
\begin{align*}
T_1\bullet T_2&=\pm\chi(\gamma\times\{0\})\pm 1=0\pm 1 =\pm 1.
\end{align*}
However, $T_1$ and $T_2$ can be displaced from each other by translation in $\mathbb{C}^2$, thus we also have
\[
T_1\bullet T_2=\operatorname{sign}(\varnothing)=0.
\]
This contradiction leads to the existence of a second intersection point distinct from the first one in $T_1\cap T_2\setminus(\gamma\times\{0\})$. This shows that for a generic class of smooth Jordan curve there are two distinct cyclic quadrilateral pegs with $(s,t,\phi)$ such that $s\neq\frac{1}{2}$ or $t\neq\frac{1}{2}$ by Proposition \ref{prop_geo_trans_cyc_quad} together with the previous discussion on intersection points and inscribed cyclic quadrilaterals. Combining with Theorem \ref{main}, we have proved the following extended statement.

\begin{thm}\label{thm_doubling_cyc_quad}
For any fixed data $s,t\in(0,\frac{1}{2}]$ and $\phi\in(0,\pi)$ except $s=t=\frac{1}{2}$ and $\phi=\frac{\pi}{2}$, there is a generic class of smooth Jordan curves such that any cyclic quadrilateral in $\mathbb{R}^2$ with the fixed data $(s,t,\phi)$ admits at least two similar ones whose vertices land distinctly on any curve in this generic class. 
\end{thm}

\appendix
\renewcommand{\thesection}{A} 

\section{Lagrangian Floer homology}\label{App_Floer_homology}

The Floer homology of Lagrangian submanifolds was first introduced by Andreas Floer, c.f. \cite{Floer1988-I,Floer1988-II,Floer1988-III} in aspherical setting, and further extended to more general setting, e.g. \cite{Oh1993,Oh1995,Oh1996}, see also \cite{Fukaya-Oh-Ohta-Ono2009,Oh2015}. In this appendix, we give an overview of the Lagrangian Floer homology within the scope of this paper.

\subsection{Floer homology of monotone Lagrangian submanifolds}\label{App_monotone_Lag_Floer}

We consider a \textit{monotone symplectic manifold} $(M,\omega)$, i.e. 
\[
I_{c_1}=\kappa I_\omega 
\] 
for some constant $\kappa\geq 0$, where $I_\omega:\pi_2(M)\to\mathbb{R}$ is the evaluation map of symplectic area, and $I_{c_1}:\pi_2(M)\to\mathbb{R}$ is the evaluation map of the first Chern class $c_1(TM)\in H^2(M;\mathbb{Z})$. Suppose there are two compact Lagrangian submanifolds $L_0$ and $L_1$ which are \textit{monotone}, namely
\begin{equation}\label{monotonicity}
I_\omega=\tau I_\mu 
\end{equation}
for some constant $\tau\geq 0$, where $I_\omega:\pi_2(M,L_i)\to\mathbb{R}$ is the evaluation map of the symplectic area, and $I_\mu:\pi_2(M,L_i)\to\mathbb{Z}$ is the evaluation map of the Maslov class $\mu\in H^1(\Lambda(\mathbb{C}^n);\mathbb{Z})$, for $i=0,1$ respectively. The \textit{minimal Maslov number} $N_{L_i}$ of $L_i$ for $i=0,1$ respectively is the positive generator of the image of $I_\mu$.

The symplectic action functional 
\begin{equation}\label{action_funct}
\mathcal{A}([x,\overline{x}])=-\int_{[0,1]\times [0,1]}\overline{x}{}^*\omega
\end{equation}
is well-defined on the Novikov covering space 
\[
\widetilde{\mathcal{P}}(L_0,L_1)
\] 
of the path space with Largangian boundary condition: 
\[
\mathcal{P}(L_0,L_1):=\{ x:[0,1]\to M: x(0)\in L_0, x(1)\in L_1    \}.
\] 
C.f. \cite[\S 2.2]{Fukaya-Oh-Ohta-Ono2009}, \cite[\S 13.4]{Oh2015}. The Novikov covering space $\widetilde{\mathcal{P}}(L_0,L_1)$ consists of equivalence classes of pairs $[x,\overline{x}]$ where $\overline{x}:[0,1]\times[0,1]\to M$ called the \textit{cap of $x$} satisfies $\overline{x}(1,\cdot)=x$, $\overline{x}(0,\cdot)=\ell_0\in\mathcal{P}(L_0,L_1)$ (i.e. based point in path space), $\overline{x}(\cdot,0)\in L_0$, $\overline{x}(\cdot,1)\in L_1$, and the equivalence relation is given by
\[
(x,\overline{x})\sim(x,\widetilde{x}):\Leftrightarrow I_\omega(-\overline{x}\#\widetilde{x})=0=I_\mu(-\overline{x}\#\widetilde{x}).
\]
The Deck transformation group is 
\begin{equation}\label{Deck_grp}
\Gamma_{\ell_0}:=\Gamma(L_0,L_1;\ell_0):=\frac{\pi_1(\mathcal{P}(L_0,L_1) , \ell_0 )  }{\ker (I_\omega)\cap \ker(I_\mu)},
\end{equation}
where $I_\omega:\pi_1(\mathcal{P}(L_0,L_1),\ell_0)\to\mathbb{R}$ and $I_\mu:\pi_1(\mathcal{P}(L_0,L_1),\ell_0)\to\mathbb{Z}$, c.f. \cite[\S 2.2.1]{Fukaya-Oh-Ohta-Ono2009}, \cite[\S 13.4]{Oh2015}. The group $\Gamma_{\ell_0}$ acts on $\widetilde{\mathcal{P}}(L_0,L_1)$ transitively and effectively. The critical points set of the action functional $\operatorname{Crit}(\mathcal{A})$ is the set of intersection points between these two Lagrangian submanifolds together with capping, i.e. 
\begin{equation}\label{crit_set}
\operatorname{Crit}(\mathcal{A})\cong (L_0\cap L_1) \times \pi_1(\mathcal{P}(L_0,L_1),\ell_0).
\end{equation}

The Floer complex is defined as follows. The Floer chain module is a free $\mathbb{F}_2$-module generated by equivalence classes
\[
[x,\overline{x}]\in\frac{\operatorname{Crit}(\mathcal{A})}{\ker(I_\omega)\cap\ker(I_\mu)}\cong (L_0\cap L_1)\times\Gamma_{\ell_0},
\] 
i.e. 
\[
CF_*(L_0,L_1,\ell_0):=\bigoplus_{\substack{ [x,\overline{x}]\in (L_0\cap L_1)\times\Gamma_{\ell_{0}} }} \mathbb{F}_2\langle [x,\overline{x}] \rangle
\]
and
\[
CF_*(L_0,L_1):=\bigoplus_{\ell_{0}}CF_*(L_0,L_1,\ell_{0}),
\]
in which $\bigoplus_{\ell_{0}}$ runs over the set of base points of connected components of $\mathcal{P}(L_0,L_1)$. This is an infinitely generated $\mathbb{F}_2$-vector space. There is a natural $\Gamma_{\ell_0}$-action endowing $CF_*(L_0,L_1,\ell_0)$ a $\mathbb{F}_2[\Gamma_{\ell_0}]$-module structure: For any $g\in\Gamma_{\ell_0}$,
\begin{equation}\label{Gamma_action}
g\cdot [x,\overline{x}]=[x,g\#\overline{x}]
\end{equation}
where the element $g$ can be regarded as a map $S^1\times[0,1]\to M$ satisfying $g(\cdot,0)\in L_0$, $g(\cdot,1)\in L_1$.The boundary operator is given by counting rigid holomorphic strips which is described as follows. We pick a smooth family $\{J_t\}_{t\in[0,1]}$ of $\omega$-compatible almost complex structures whose set is denoted by $\mathcal{J}(M,\omega)$. This induces a smooth family of Riemannian metrics $g_{J_t}(\cdot,\cdot)=\omega(\cdot,J_t\cdot)$, and further a $L^2$-inner product 
\[
\langle\cdot,\cdot\rangle_{J}:=\int_{0}^{1}\omega(\cdot,J_t\cdot) dt
\] 
on $\mathcal{P}(L_0,L_1)$ and hence on $\widetilde{\mathcal{P}}(L_0,L_1)$. The negative $L^2$-gradient flow equation of the symplectic action functional $\mathcal{A}$ reads: smooth map $u:\mathbb{R}\times[0,1]\to M$ solving 
\begin{equation}\label{Floer_eq}
\begin{cases}
\partial_s u+J_t\partial_t u=0,\\
u(s,0)\in L_0, u(s,1)\in L_1.
\end{cases}
\end{equation}
This is literally the non-linear Cauchy-Riemann equation 
\[
\bar{\partial}_{J_t}u=0
\]
for $J_t$-holomorphic strip with boundary condition on $L_0$ and $L_1$ in terms of the coordinate $(s,t)\in\mathbb{R}\times[0,1]$.  The energy or the $L^2$-norm of the $J_t$-holomorphic strip is equal to its symplectic area:
\begin{align*}
E(u)&:=\int_{\mathbb{R}\times[0,1]}\left\|\partial_s u\right\|_{J_t}^2  ds \wedge dt \\
&=\frac{1}{2}\int_{\mathbb{R}\times[0,1]} ( \left\|\partial_s u\right\|_{J_t}^2+\left\|\partial_t u \right\|_{J_t}^2 )ds \wedge dt \\
&=\int_{\mathbb{R}\times[0,1]} u^*\omega
\end{align*}
When the energy of $J_t$-holomorphic strip is finite, i.e. 
\[
E(u)<\infty,
\] 
then the solution has asymptotic behavior:
\[
\lim_{s\to\pm\infty}u(s,\cdot)=x_\pm\in L_0\cap L_1
\]
and 
\[
E(u)=\mathcal{A}([x_-,\overline{x}_-])-\mathcal{A}([x_+,\overline{x}_+])
\]
if $[x_+,\overline{x}_+]=[x_+,\overline{x}_-\# u]$, vice versa. C.f. \cite[\S 2.3.1]{Fukaya-Oh-Ohta-Ono2009}, \cite[\S 13.9.1]{Oh2015}. We denote by
\[
\mathcal{M}(x_-,x_+;J_t):=\widetilde{\mathcal{M}}(x_-,x_+;J_t)/\mathbb{R}
\]
the \textit{(unparameterized) moduli space of $J_t$-holomorphic strips} with $x_{\pm}\in L_0\cap L_1$ as asymptotics, where
\[
\widetilde{\mathcal{M}}(x_-,x_+;J_t):=\left\{ u\in C^\infty(\mathbb{R}\times[0,1],M) :  \text{$u$ solves \eqref{Floer_eq} with $\lim_{s\to\pm\infty}u(s,\cdot)=x_{\pm}$} \right\}
\]
is the \textit{parameterized moduli space} which carries a proper and transitive $\mathbb{R}$-action by shift of $s$-parameter. There is a decomposition 
\[
\widetilde{\mathcal{M}}(x_-,x_+;J_t)=\bigcup_{B\in\pi_2(x_-,x_+)}\widetilde{\mathcal{M}}(x_-,x_+;B,J_t)
\]
where $\pi_2(x_-,x_+)$ is the set of homotopy classes of maps $u:[0,1]\times[0,1]\to M$ relative to boundary $u(0,\cdot)=x_-$, $u(1,\cdot)=x_+$, $u(\cdot,0)\in L_0$, $u(\cdot,1)\in L_1$, it is a subset in the relative homotopy group $\pi_2(M,L_0\cup L_1)$. The space $\widetilde{\mathcal{M}}(x_-,x_+;B,J_t)$ consists of those $J$-holomorphic strips $u$ with homotopy class $[u]=B$. The virtual dimension of the moduli space is computed by Fredholm index of the Cauchy-Riemann operator $\bar{\partial}_{J_t}$ at $u$ (c.f. \cite[Theorem 15.3.2]{Oh2015}):
\begin{equation}\label{vir_dim}
\dim\widetilde{\mathcal{M}}(x_-,x_+;B,J_t)=\operatorname{ind}(u)=\mu(x_-,x_+;B)
\end{equation}
where 
\[
\operatorname{ind}(u):=\operatorname{ind}(\bar{\partial}_{J_t}|_u):=\operatorname{ind}(D_u\bar{\partial}_{J_{t}})
\]
is the Fredholm index of the Fredholm operator $\bar{\partial}_{J_t}|_u$ associated to $u$, and 
\[
\mu(x_-,x_+;B)\in\mathbb{Z}
\] 
is the Maslov--Viterbo index of $u$ introduced by \cite{Viterbo1987}, see also \cite{Floer1988-II}, \cite[Definition 2.3.9]{Fukaya-Oh-Ohta-Ono2009}, \cite[Definition 13.6.2]{Oh2015}. The Maslov--Viterbo index of $u$ can be expressed as the difference of the Maslov--Morse index (definition refers to \cite[Definition 2.2.12]{Fukaya-Oh-Ohta-Ono2009}, \cite[Definition 13.6.2]{Oh2015}) of the pairs $[x_-,\overline{x}_-]$ and $[x_+,\overline{x}_+]$:
\begin{equation}\label{rel_abs_index}
\mu(x_-,x_+;B)=\mu([x_+,\overline{x}_+])-\mu([x_-,\overline{x}_-])
\end{equation}
with $\overline{x}_+=\overline{x}_-\# B$. C.f. \cite[Proposition 2.3.9]{Fukaya-Oh-Ohta-Ono2009}, \cite[Proposition 13.6.3]{Oh2015}. This formula is independent of the choice of the cap $\overline{x}_-$. In some nice situations, one has transversality which would imply that, for a generic choice of $\{J_t\}_{t\in[0,1]}\subset\mathcal{J}(M,\omega)$, the parameterized moduli space $\widetilde{\mathcal{M}}(x_-,x_+;B,J_t)$ is a smooth manifold of expected dimension \eqref{vir_dim}, and hence the (unparameterized) moduli space $\mathcal{M}(x_-,x_+;B,J_t)$ is a smooth manifold of dimension 
\begin{equation}
\dim \mathcal{M}(x_-,x_+;B,J_t)=\mu(x_-,x_+;B)-1.
\end{equation}
E.g. refer to the statement \cite[Proposition 3.3]{Oh1993}, \cite[Theorem 2.3.11]{Fukaya-Oh-Ohta-Ono2009}, \cite[Theorem 15.1.3]{Oh2015}. In general, the moduli space $\mathcal{M}(x_-,x_+;J_t)$ can be further compactified due to the Gromov--Floer convergence, by adding broken $J_t$-holomorphic strips, $J_t$-holomorphic disks and $J_t$-holomorphic spheres (these lateral two are also known as disk and sphere bubbling), into a compact space $\overline{\mathcal{M}}(x_-,x_+;J_t)$ with respect to the Gromov--Floer topology. Provided above transversality, for a generic choice of $\{J_t\}_{t\in[0,1]}\subset\mathcal{J}(M,\omega)$, the compact space $\overline{\mathcal{M}}(x_-,x_+;J_t)$ is a stratified space, whose strata are supposed to carry manifold structures. E.g. for the purpose of defining Floer complex, we would like to have enough transversality so that the $0$-dimensional stratum $\mathcal{M}^0$ is compact $0$-dimensional manifold (i.e. finite point set), and the $1$-dimensional stratum $\overline{\mathcal{M}}{}^1$ is a compact $1$-dimensional manifold with boundary $\partial\overline{\mathcal{M}}{}^1=\mathcal{M}^0$. The details for theses topological structures of moduli spaces are rather technical, for nice expositions, we refer to \cite[Chapter 2]{Fukaya-Oh-Ohta-Ono2009}, \cite[\S 14, \S 15, \S 16]{Oh2015}. Now with all the above necessary ingredients, the boundary operator is written as
\begin{align*}
	\partial:CF_{*}(L_0,L_1)&\to CF_{*-1}(L_0,L_1)\\
	[x,\overline{x}]&\mapsto\sum_{\substack{y\in L_0\cap L_1 \\ B\in \pi_2(x,y) \\ \mu(x,y;B)=1 }}\#_2\mathcal{M}^{0}(x,y;B,J_t)  [y,\overline{y}]
\end{align*}
where $\#_2$ stands for mod 2 count of the 0-dimensional moduli space $\mathcal{M}^{0}(x,y;B,J_t)$. In general, due to the formula \eqref{rel_abs_index}, one can put a relative grading or an absolute grading modulo some integer (i.e. a periodic grading) on the Floer complex provided by the Maslov--Morse index (c.f. \cite[\S 13.9, p.65]{Oh2015}): For the pair $[x,\overline{x}]$, due to \eqref{Gamma_action}, the grading differs by
\begin{equation}\label{index_differ}
	\mu([x,g\#\overline{x}])=\mu([x,\overline{x}])-I_{\mu}(g) 
\end{equation}
for $g\in\Gamma_{\ell_0}$, note that we can view the difference of various cappings, i.e. $(-\widetilde{x})\#\overline{x}$ as such an element $g$. See also \cite{Floer1988-II}, \cite[Lemma 4.7]{Oh1993}. When the boundary operator satisfies the condition
\[
\partial\circ\partial =0,
\]
we say that the Lagrangian submanifolds are \textit{unobstructed}. In such case, the Lagrangian Floer homology is well-defined: 
\[
HF_*(L_0,L_1):=\frac{\ker(\partial:CF_{*}\to CF_{*-1})}{\operatorname{im}(\partial:CF_{*+1}\to CF_{*})}.
\]
Otherwise, we say the Lagrangian submanifolds are \textit{obstructed}, and there is \textit{anomaly}. The appearance of anomaly often concerns the disk bubbling, c.f. \cite{Oh1995}, \cite[\S 2.4]{Fukaya-Oh-Ohta-Ono2009}, \cite[16.4.2]{Oh1995}. 

In the monotone case, as we have indicated at the beginning of this section, the Lagrangian Floer homology is established by Yong-Geun Oh.

\begin{thm}[\cite{Oh1993,Oh1995},{\cite[\S 16]{Oh2015}}]\label{thm_HF_monotone}
	Let $L_0,L_1$ be compact monotone Lagrangian submanifolds in a compact monotone symplectic manifold $(M,\omega)$ with minimal Malsov number $N_{L_0},N_{L_1}$, and put $N:=\min\{N_{L_0},N_{L_1} \}$. 
	\begin{itemize}
		\item[(i)] Suppose $L_0\pitchfork L_1$ and $N\geq 3$. Then the Floer homology $HF_*(L_0,L_1)$ is well-defined as a $\mathbb{F}_2$-module with $\mathbb{Z}/N\mathbb{Z}$-grading. 
		
		\item[(ii)] Suppose $N= 2$, and $L_1=\varphi_H(L_0)$ with $\varphi_H$ a Hamiltonian diffeomorphism generated by a non-degenerate Hamiltonian, i.e. $L_0\pitchfork\varphi_H(L_0)$. Then the Floer homology $HF_*(L_0,\varphi_H(L_0))$ is well-defined as a $\mathbb{F}_2$-module with $\mathbb{Z}/N\mathbb{Z}$-grading. 
	\end{itemize}
	In particular, the Floer homology $HF_*(L_0,L_1)$ is independent of the choice of $\omega$-compatible almost complex structures and invariant under Hamiltonian diffeomorphisms.
\end{thm}

\begin{figure}
	\[
	\begin{tikzpicture}[x=0.75pt,y=0.75pt,yscale=-1,xscale=1]
		
		\draw    (183,57) .. controls (223,30) and (243,30) .. (283,57) ;
		\draw    (183,51) .. controls (223,73) and (243,74) .. (283,51) ;
		\draw    (363,57.59) .. controls (383,36) and (399,32) .. (423,57.59) ;
		\draw    (363,52) .. controls (385,75) and (398,74) .. (423,52) ;
		\draw    (423,52) .. controls (442,34) and (459,39) .. (476,57) ;
		\draw    (423,57.59) .. controls (441,74) and (456,76) .. (475,51) ;
		
		\draw (225,75.4) node [anchor=north west][inner sep=0.75pt]    {$L_{0}$};
		\draw (225,10.4) node [anchor=north west][inner sep=0.75pt]    {$L_{1}$};
		\draw (228,47.4) node [anchor=north west][inner sep=0.75pt]    {$2$};
		\draw (415.4,73.15) node [anchor=north west][inner sep=0.75pt]    {$L_{0}$};
		\draw (415.4,12.65) node [anchor=north west][inner sep=0.75pt]    {$L_{1}$};
		\draw (386.6,49.26) node [anchor=north west][inner sep=0.75pt]    {$1$};
		\draw (440.6,50.26) node [anchor=north west][inner sep=0.75pt]    {$1$};
		\draw (309,48.4) node [anchor=north west][inner sep=0.75pt]    {$\rightarrow $};

	\end{tikzpicture}
	\]
	\caption{$J_t$-holomorphic strip of Maslov index $2$ breaking into two $J$-holomorphic strips of Maslov index $1$.}\label{Fig_10}
\end{figure}

\begin{figure}
	\[
	\begin{tikzpicture}[x=0.75pt,y=0.75pt,yscale=-1,xscale=1]
		
		\draw    (153,74) .. controls (193,47) and (213,47) .. (253,74) ;
		\draw    (153,68) .. controls (193,90) and (213,91) .. (253,68) ;
		\draw    (339,74.3) .. controls (367.8,53.2) and (382.2,53.2) .. (411,74.3) ;
		\draw    (339,69.61) .. controls (367.8,86.8) and (382.2,87.59) .. (411,69.61) ;
		\draw   (361,43) .. controls (361,34.72) and (367.72,28) .. (376,28) .. controls (384.28,28) and (391,34.72) .. (391,43) .. controls (391,51.28) and (384.28,58) .. (376,58) .. controls (367.72,58) and (361,51.28) .. (361,43) -- cycle ;
		\draw    (440,75.7) .. controls (468.8,54.6) and (483.2,54.6) .. (512,75.7) ;
		\draw    (440,71.01) .. controls (468.8,88.21) and (483.2,88.99) .. (512,71.01) ;
		\draw   (461,99.4) .. controls (461,91.12) and (467.72,84.4) .. (476,84.4) .. controls (484.28,84.4) and (491,91.12) .. (491,99.4) .. controls (491,107.69) and (484.28,114.4) .. (476,114.4) .. controls (467.72,114.4) and (461,107.69) .. (461,99.4) -- cycle ;
		
		\draw (195,92.4) node [anchor=north west][inner sep=0.75pt]    {$L_{0}$};
		\draw (195,27.4) node [anchor=north west][inner sep=0.75pt]    {$L_{1}$};
		\draw (198,63.4) node [anchor=north west][inner sep=0.75pt]    {$2$};
		\draw (279,65.4) node [anchor=north west][inner sep=0.75pt]    {$\rightarrow $};
		\draw (367.58,99.8) node [anchor=north west][inner sep=0.75pt]    {$L_{0}$};
		\draw (363.58,5) node [anchor=north west][inner sep=0.75pt]    {$L_{1}$};
		\draw (370.72,37.14) node [anchor=north west][inner sep=0.75pt]    {$2$};
		\draw (370.72,64.14) node [anchor=north west][inner sep=0.75pt]    {$0$};
		\draw (471.58,120.2) node [anchor=north west][inner sep=0.75pt]    {$L_{0}$};
		\draw (467.58,35.4) node [anchor=north west][inner sep=0.75pt]    {$L_{1}$};
		\draw (470.72,95.54) node [anchor=north west][inner sep=0.75pt]    {$2$};
		\draw (471.72,65.54) node [anchor=north west][inner sep=0.75pt]    {$0$};

	\end{tikzpicture}
	\]
	\caption{$J_t$-holomorphic strip of Maslov index $2$ bubbling off $J$-holomorphic disks of Maslov index $2$.}\label{Fig_11}
\end{figure}

\begin{rem}\label{rem_bubbling}
Since the compactified moduli space $\overline{\mathcal{M}}{}^1$ of dimension 1 is involved in the definition of boundary operator $\partial$ and the boundary stratum $\mathcal{M}^0=\partial\overline{\mathcal{M}}{}^1$ is relevant for the condition $\partial^2=0$. These are all smooth manifolds in monotone case, c.f.  \cite[16.4.2]{Oh2015}. When $N\geq 3$, there are only once-broken $J_t$-holomorphic strips, whose each component is of Maslov index $1$, in the boundary stratum due to dimensional reason. See Figure \ref{Fig_10}. When $N=2$, there are additionally $J_t$-holomorphic disks of Maslov index $2$, i.e. disk bubbling, in the boundary stratum. See Figure \ref{Fig_11}. C.f. \cite[Theorem 16.4.4]{Oh2015}. In particular, there is obstruction due the appearance of disk bubbling, namely the boundary operator satisfies
\begin{equation}\label{Oh's_formula}
\partial\circ\partial = \Phi_{L_0}(\text{pt.})-\Phi_{L_1}(\text{pt.}).
\end{equation} 
where $\Phi_{L_i}(\text{pt.})$ is the one-point open Gromov--Witten invariant for Lagrangian submanifolds $L_{i}$, $i=0,1$, and ``$\text{pt}.$'' represents an arbitrary point in $L_i$ respectively, c.f. \cite[Example 15.4.3, p.163]{Oh2015}. If $L_1=\varphi_H(L_0)$ for some Hamiltonian diffeomorphism $\varphi_H$, then it is shown that the right hand side of \eqref{Oh's_formula} cancels, hence the obstruction vanishes. Moreover, this also holds when working with the (universal) Novikov coefficient, which is to be introduced succeedingly.  C.f. \cite{Oh1995}, \cite[\S 16.4.2]{Oh2015}, and \cite[\S 2.4.5]{Fukaya-Oh-Ohta-Ono2009}.
\end{rem}

It is more delightful to work with an absolutely $\mathbb{Z}$-graded bounded chain complex and homology instead of the $N$-periodic unbounded chain complex and homology. For this purpose, we introduce the universal Novikov ring (field) $\Lambda$ over $\mathbb{F}_2$ and define the Floer homology with repsec to such coefficient:
\begin{equation}\label{Nov_ring}
\Lambda:=\left.\left\{ \sum_{i=0}^{\infty}c_i T^{a_i} e^{\frac{d_i}{2}} \; \right|\; c_i\in\mathbb{F}_2,\;d_i\in\mathbb{Z},\; a_i\in\mathbb{R},\;a_i\leq a_{i+1},\;\lim_{i\to\infty}a_i=\infty  \right\} 
\end{equation}
where $T$ and $e$ are formal parameters whose gradings are put to be $|T|=0$ and $|e|=2$. A systematic discussion of Novikov rings can be found in \cite[\S 2.2.1, \S 2.4.3, \S 5.1]{Fukaya-Oh-Ohta-Ono2009}, \cite[\S 13.4, \S 13.8]{Oh2015}. More precisely, the ring $\Lambda$ here is actually the extended universal Novikov field $\Lambda_{nov}[e^{\frac{1}{2}}]$ in terms of the notion of \cite[\S 5.1.1, \S 5.1.3]{Fukaya-Oh-Ohta-Ono2009}. The Floer complex is now defined as the free $\Lambda$-module finitely generated by $L_0\cap L_1$, i.e. 
\[
CF_*(L_0,L_1;\Lambda)=\bigoplus_{\substack{ x\in L_0\cap L_1 }} \mathbb{F}_2\langle x \rangle\otimes_{\mathbb{F}_2}\Lambda.
\]
The boundary operator is defined by
\begin{align*}
	\partial:CF_{*}(L_0,L_1;\Lambda)&\to CF_{*-1}(L_0,L_1;\Lambda)\\
	x&\mapsto\sum_{\substack{  y\in L_0\cap L_1  \\ B\in\pi_2(x,y)  \\  \mu(x,y;B) =1 } }\#_2\mathcal{M}^{0}(x,y;B,J_t) y\otimes T^{\omega(B)}e^{\frac{\mu(B)}{2}}
\end{align*}
via the mod 2 count of $0$-dimensional moduli space $\mathcal{M}^0(x,y;B,J_t)$. In particular, $e^{\frac{\mu(B)}{2}}=e^{\frac{1}{2}}$ has degree $1$ provided $\mu(x,y;B)=1$. Fixing a cap for each intersection point $x\in L_0\cap L_1$, we can equip an absolute $\mathbb{Z}$-grading for the Floer complex, c.f. \cite[\S 5.1.3]{Fukaya-Oh-Ohta-Ono2009}. This should be distinguished from the canonical absolute $\mathbb{Z}$-grading introduced in \cite{Seidel2000} when $2c_1(TM)=0$ and the Maslov classes of the oriented Lagrangian submanifolds $L_0$ and $L_1$ vanish, see \cite[Remark 5.1.18]{Fukaya-Oh-Ohta-Ono2009}. In this sense, we obtain a $\mathbb{Z}$-graded Floer homology 
\[
HF_*(L_0,L_1;\Lambda):=H_*(CF_{*}(L_0,L_1;\Lambda),\partial)
\] 
as $\Lambda$-module which is also invariant under Hamiltonian diffeomorphisms. In particular, when the set of intersection points $L_0\cap L_1$ is finite, the Floer complex is a bounded $\mathbb{Z}$-graded chain complex over $\Lambda$, and so is the Floer homology a bounded $\mathbb{Z}$-graded $\Lambda$-module.

\subsection{Morse--Bott Floer homology}\label{App_MBFloer}
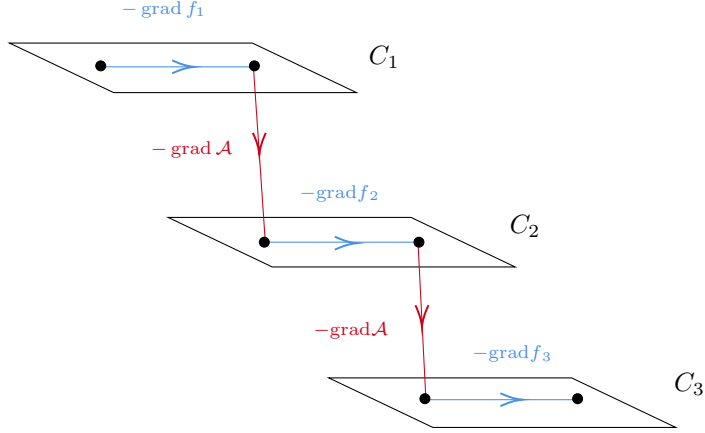
\begin{figure}
\[
\begin{tikzpicture}[x=0.75pt,y=0.75pt,yscale=-1,xscale=1]
	
	\draw   (266.5,30) -- (144,30) -- (196.5,55) -- (319,55) -- cycle ;
	\draw   (346.5,118) -- (224,118) -- (276.5,143) -- (399,143) -- cycle ;
	\draw   (427.5,199) -- (305,199) -- (357.5,224) -- (480,224) -- cycle ;
	\draw [color={rgb, 255:red, 74; green, 144; blue, 226 }  ,draw opacity=1 ]   (190,42) -- (267,42) ;
	\draw [color={rgb, 255:red, 74; green, 144; blue, 226 }  ,draw opacity=1 ]   (273,130.5) -- (350,130.5) ;
	\draw [color={rgb, 255:red, 74; green, 144; blue, 226 }  ,draw opacity=1 ]   (354,210) -- (431,210) ;
	\draw [color={rgb, 255:red, 74; green, 144; blue, 226 }  ,draw opacity=1 ]   (237,42) ;
	\draw [shift={(237,42)}, rotate = 180] [color={rgb, 255:red, 74; green, 144; blue, 226 }  ,draw opacity=1 ][line width=0.75]    (10.93,-3.29) .. controls (6.95,-1.4) and (3.31,-0.3) .. (0,0) .. controls (3.31,0.3) and (6.95,1.4) .. (10.93,3.29)   ;
	\draw [color={rgb, 255:red, 74; green, 144; blue, 226 }  ,draw opacity=1 ]   (319,131) ;
	\draw [shift={(319,131)}, rotate = 180] [color={rgb, 255:red, 74; green, 144; blue, 226 }  ,draw opacity=1 ][line width=0.75]    (10.93,-3.29) .. controls (6.95,-1.4) and (3.31,-0.3) .. (0,0) .. controls (3.31,0.3) and (6.95,1.4) .. (10.93,3.29)   ;
	\draw [color={rgb, 255:red, 74; green, 144; blue, 226 }  ,draw opacity=1 ]   (402.5,210) ;
	\draw [shift={(402.5,210)}, rotate = 180] [color={rgb, 255:red, 74; green, 144; blue, 226 }  ,draw opacity=1 ][line width=0.75]    (10.93,-3.29) .. controls (6.95,-1.4) and (3.31,-0.3) .. (0,0) .. controls (3.31,0.3) and (6.95,1.4) .. (10.93,3.29)   ;
	\draw [color={rgb, 255:red, 208; green, 2; blue, 27 }  ,draw opacity=1 ]   (267,42) -- (273,130.5) ;
	\draw [color={rgb, 255:red, 208; green, 2; blue, 27 }  ,draw opacity=1 ]   (270,82) -- (270,84.25) ;
	\draw [shift={(270,86.25)}, rotate = 270] [color={rgb, 255:red, 208; green, 2; blue, 27 }  ,draw opacity=1 ][line width=0.75]    (10.93,-3.29) .. controls (6.95,-1.4) and (3.31,-0.3) .. (0,0) .. controls (3.31,0.3) and (6.95,1.4) .. (10.93,3.29)   ;
	\draw [color={rgb, 255:red, 208; green, 2; blue, 27 }  ,draw opacity=1 ]   (350,130.5) -- (354,210) ;
	\draw [color={rgb, 255:red, 208; green, 2; blue, 27 }  ,draw opacity=1 ]   (352,169) -- (352,171.25) ;
	\draw [shift={(352,173.25)}, rotate = 270] [color={rgb, 255:red, 208; green, 2; blue, 27 }  ,draw opacity=1 ][line width=0.75]    (10.93,-3.29) .. controls (6.95,-1.4) and (3.31,-0.3) .. (0,0) .. controls (3.31,0.3) and (6.95,1.4) .. (10.93,3.29)   ;
	\draw  [fill={rgb, 255:red, 0; green, 0; blue, 0 }  ,fill opacity=1 ] (187.5,41.5) .. controls (187.5,40.12) and (188.62,39) .. (190,39) .. controls (191.38,39) and (192.5,40.12) .. (192.5,41.5) .. controls (192.5,42.88) and (191.38,44) .. (190,44) .. controls (188.62,44) and (187.5,42.88) .. (187.5,41.5) -- cycle ;
	\draw  [fill={rgb, 255:red, 0; green, 0; blue, 0 }  ,fill opacity=1 ] (265,41.5) .. controls (265,40.12) and (266.12,39) .. (267.5,39) .. controls (268.88,39) and (270,40.12) .. (270,41.5) .. controls (270,42.88) and (268.88,44) .. (267.5,44) .. controls (266.12,44) and (265,42.88) .. (265,41.5) -- cycle ;
	\draw  [fill={rgb, 255:red, 0; green, 0; blue, 0 }  ,fill opacity=1 ] (270,130.5) .. controls (270,129.12) and (271.12,128) .. (272.5,128) .. controls (273.88,128) and (275,129.12) .. (275,130.5) .. controls (275,131.88) and (273.88,133) .. (272.5,133) .. controls (271.12,133) and (270,131.88) .. (270,130.5) -- cycle ;
	\draw  [fill={rgb, 255:red, 0; green, 0; blue, 0 }  ,fill opacity=1 ] (348,130.5) .. controls (348,129.12) and (349.12,128) .. (350.5,128) .. controls (351.88,128) and (353,129.12) .. (353,130.5) .. controls (353,131.88) and (351.88,133) .. (350.5,133) .. controls (349.12,133) and (348,131.88) .. (348,130.5) -- cycle ;
	\draw  [fill={rgb, 255:red, 0; green, 0; blue, 0 }  ,fill opacity=1 ] (351,209.5) .. controls (351,208.12) and (352.12,207) .. (353.5,207) .. controls (354.88,207) and (356,208.12) .. (356,209.5) .. controls (356,210.88) and (354.88,212) .. (353.5,212) .. controls (352.12,212) and (351,210.88) .. (351,209.5) -- cycle ;
	\draw  [fill={rgb, 255:red, 0; green, 0; blue, 0 }  ,fill opacity=1 ] (428,209.5) .. controls (428,208.12) and (429.12,207) .. (430.5,207) .. controls (431.88,207) and (433,208.12) .. (433,209.5) .. controls (433,210.88) and (431.88,212) .. (430.5,212) .. controls (429.12,212) and (428,210.88) .. (428,209.5) -- cycle ;
	
	\draw (324,30.4) node [anchor=north west][inner sep=0.75pt]    {$C_{1}$};
	\draw (395,115.4) node [anchor=north west][inner sep=0.75pt]    {$C_{2}$};
	\draw (478,195.4) node [anchor=north west][inner sep=0.75pt]    {$C_{3}$};
	\draw (199,7.4) node [anchor=north west][inner sep=0.75pt]  [font=\scriptsize]  {$\textcolor[rgb]{0.29,0.56,0.89}{-\operatorname{grad} f_{1}}$};
	\draw (214,78.4) node [anchor=north west][inner sep=0.75pt]  [font=\scriptsize]  {$\textcolor[rgb]{0.82,0.01,0.11}{-\operatorname{grad}\mathcal{A}}$};
	\draw (296,169.4) node [anchor=north west][inner sep=0.75pt]  [font=\scriptsize]  {$\textcolor[rgb]{0.82,0.01,0.11}{-}\textcolor[rgb]{0.82,0.01,0.11}{\operatorname{grad}}\textcolor[rgb]{0.82,0.01,0.11}{\mathcal{A}}$};
	\draw (289,100.4) node [anchor=north west][inner sep=0.75pt]  [font=\scriptsize]  {$\textcolor[rgb]{0.29,0.56,0.89}{-}\textcolor[rgb]{0.29,0.56,0.89}{\operatorname{grad}}\textcolor[rgb]{0.29,0.56,0.89}{f}\textcolor[rgb]{0.29,0.56,0.89}{_{2}}$};
	\draw (376,181.4) node [anchor=north west][inner sep=0.75pt]  [font=\scriptsize]  {$\textcolor[rgb]{0.29,0.56,0.89}{-}\textcolor[rgb]{0.29,0.56,0.89}{\operatorname{grad}}\textcolor[rgb]{0.29,0.56,0.89}{f}\textcolor[rgb]{0.29,0.56,0.89}{_{3}}$};

\end{tikzpicture}
\]
\caption{A cascade trajectory.}\label{Fig_12}
\end{figure}

Lagrangian Floer theory in the Morse--Bott setting was introduced and developed in \cite{Pozniak1994,Frauenfelder2003,Cornea-Lalonde2006,Biran-Cornea2009a,Biran-Cornea2009b,Fukaya-Oh-Ohta-Ono2009,Johns2010,Sheridan2011,Seidel2011,Seidel2014,Schmaeschke2016}. We follow the cascade approach in \cite{Frauenfelder2003,Schmaeschke2016}. Assume $L_0$ and $L_1$ intersect cleanly, namely $L_0\cap L_1$ consists of a disjoint union of connected submanifolds. Suppose $C\subseteq L_0\cap L_1$, then $T_pL_0\cap T_pL_1=T_p C$ $\forall p\in C$. The action functional \eqref{action_funct} is a Morse--Bott function on the path space in the sense that
\[
\ker\left(\operatorname{Hess}\mathcal{A}|_{T_p\mathcal{P}(L_0,L_1)}\right)=T_p C, \qquad \forall p\in C.
\]
Suppose 
\[
L_0\cap L_1=\bigsqcup_{i\in I}C_i
\]
where each $C_i$ is a connected closed submanifold. We choose auxiliary Morse functions $f_i$ on $C_i$. The Floer chain $\Lambda$-module is now freely generated by critical points of Morse function $x\in\operatorname{Crit}(f_i)$, $i\in I$:
\[
CF_*(L_0,L_1,\{f_i\};\Lambda):=\bigoplus_{i\in I}\bigoplus_{  x\in\operatorname{Crit}(f_i) }\Lambda\langle x \rangle.
\]
The grading of $x\in\operatorname{Crit}(f_i)$ is now given by
\begin{equation}\label{MB_grading}
	|x|=\operatorname{ind}_{f_i}(x)+\mu([x,\overline{x}])-\frac{1}{2}\dim C_i, \qquad i\in I,
\end{equation}
where $\operatorname{ind}_{f_i}(x)$ is the Morse index of $x$, and $\overline{x}$ is a fixed capping of $x$ in a canonical way (c.f.\cite[\S3.7 (3.7.49)]{Fukaya-Oh-Ohta-Ono2009}) such that $\mu([x,\overline{x}])$ stays the same for all points $x\in C_i$. Therefore, one can assign an integral grading $\mu(C_i)$ for each Morse--Bott component $C_i$ in $L_0\cap L_1$ by \cite[Appendix C, p.177]{Frauenfelder2003}, \cite[\S 9.1.2]{Schmaeschke2016}:
\begin{equation}\label{MB_index}
\mu(C_i):= \mu([x,\overline{x}])-\frac{1}{2}\dim C_i, \qquad i\in I.
\end{equation}
To define the boundary operator, a generalization of Floer trajectories called \textit{cascades} is counted. The formal definition of trajectories with cascades can be found in \cite[Definition C.22]{Frauenfelder2003}, \cite[\S 9.1.1]{Schmaeschke2016}. These are tuples of $J$-holomorphic strips (or Floer trajectories) and Morse trajectories connecting the critical points in sequel. See Figure \ref{Fig_12}.
For generic choice of $\{J_t\}_{t\in[0,1]}\subset\mathcal{J}(M,\omega)$ and Riemannian metric $g$ on $L_0\cap L_1$, the (unparameterized) moduli space $\mathcal{M}(x,y;J_t,g)$ of cascades connecting $x\in \operatorname{Crit}(f_i)$ and $y\in\operatorname{Crit}(f_j)$ is a smooth manifold of dimension
\begin{equation}\label{MB_dimension}
\dim\mathcal{M}(x,y;J_t,g) = |x|-|y|-1.
\end{equation}
C.f. \cite[Theorem C.23]{Frauenfelder2003}, \cite[\S 9.2.1]{Schmaeschke2016}. There is also a Gromov--Floer compactification for the moduli space of cascades, c.f. \cite[\S 9.2.2]{Schmaeschke2016}. And in the monotone setting, despite the broken trajectories, disk bubbling can also appear, refer to Remark \ref{rem_bubbling}. The boundary operator is defined by
\begin{align*}
\partial:CF_*(L_0,L_1,\{f_i\};\Lambda)&\to CF_{*-1}(L_0,L_1,\{f_i\};\Lambda) \\
x&\mapsto \sum_{y} \#_2\mathcal{M}^0(x,y;J_t,g) y \otimes T^{\mathcal{A}([y,\overline{y}])-\mathcal{A}([x,\overline{x}])} e^{\frac{\mu([x,\overline{x}])-\mu([y,\overline{y}])}{2}}
\end{align*}
When working with ordinary ring coefficient, e.g. $\mathbb{F}_2$, there are more generators in the Morse--Bott Floer complex since $\mathcal{A}$ is actually defined on the Novikov covering space $\widetilde{\mathcal{{P}}}(L_0,L_1)$ and it yields $\Gamma_{\ell_0}$ copies (up to equivalence) of those Morse--Bott components in $L_0\cap L_1$ as critical submanifolds due to \eqref{crit_set}. The Morse--Bott Floer complex with $\mathbb{F}_2$ coefficient is defined analogously to Appendix \ref{App_monotone_Lag_Floer} with auxiliary Morse functions $\{f_i^{\nu}\}^{ \nu\in\Gamma_{\ell_0}}_{i\in I }$ and metric $g$ on all Morse--Bott critical submanifolds $\{C_i^{\nu}\}^{ \nu\in\Gamma_{\ell_0}}_{i\in I }$ of $\mathcal{A}$ in $\widetilde{\mathcal{{P}}}(L_0,L_1)$ (replicating the data for all copies labeled by $\nu\in\Gamma_{\ell_0}$ with respect to a fixed $i\in I$), i.e. 
\[
CF_*(L_0,L_1,\{f_i^\nu\}):=\bigoplus_{\ell_0}\bigoplus_{\nu,i}\bigoplus_{  x^\nu\in\operatorname{Crit}(f_i^{\nu}) }\mathbb{F}_2\langle x^\nu \rangle
\]
where the grading is given by \eqref{MB_grading}. The caps determine the labeling $\nu\in\Gamma_{\ell_0}$ for the copies in the Novikov covering space, namely we have identification $x^\nu\equiv[x,\overline{x}]$. The boundary operator is defined similarly by mod 2 counts of 0-dimensional moduli spaces. The Morse--Bott Floer homology $HF_*(L_0,L_1,\{f_i^\nu\})$ is then a $\mathbb{Z}/N\mathbb{Z}$-graded $\mathbb{F}_2$-module, while $HF_*(L_0,L_1,\{f_i\};\Lambda)$ is a bounded $\mathbb{Z}$-graded $\Lambda$-module.

In the monotone case as in Theorem \ref{thm_HF_monotone}, we have the invariance theorem:
\begin{thm}[{\cite[Appendix C.4]{Frauenfelder2003},\cite[\S 9.3]{Schmaeschke2016}}]\label{MB_homology_thm}
The Floer homology 
\[
HF_*(L_0,L_1,\{f_i\},J_t,g;\Lambda):=H_*(CF_*(L_0,L_1,\{f_i\},J_t,g;\Lambda),\partial)
\] is well-defined, independent of the choice of $(\{f_i\},J_t,g)$. Moreover it is invariant under Hamiltonian diffeomorphism. Namely, there is a canonical isomorphism
\[
HF_*(L_0,L_1,\{f_i\},J_t,g;\Lambda)\cong HF_*(L_0,\varphi_H(L_1),\{\tilde{f}_j\},\tilde{J}_t,\tilde{g};\Lambda)
\]
for any $\varphi_H\in\operatorname{Ham}(M,\omega)$ such that $L$ intersects with $\varphi_H(L)$ cleanly,  $(\{f_i\},J_t,g)$ and $(\{\tilde{f}_j\},\tilde{J}_t,\tilde{g})$ are two arbitrary generic choices of collection of admissible data. These hold also for ordinary ring coeffcients, e.g. $\mathbb{F}_2$.
\end{thm}

\subsection{Liouville manifolds}\label{App_Liouv_MfD} The previous discussions for Floer homology are all settled in compact (closed) symplectic manifolds. There is suitable generalization to non-compact symplectic manifolds, e.g. exact symplectic manifolds. In particular, we focus on the most studied case called the Liouville manifold. A \textit{Liouville manifold} is an exact symplectic manifold $(M,d\lambda)$ with a primitive $\lambda$ called Liouville 1-form, and a complete Liouville vector field $Z$, i.e. $i_Zd\lambda=\lambda$, and an exhaustion by Liouville domains. More precisely, a \textit{Liouville domain} is a compact exact symplectic manifold $(M,d\lambda)$ with contact type boundary $\partial M$ (i.e. $\alpha:=\lambda|_{\partial M}$ is a contact 1-form) such that the Liouville vector field $Z$, which is uniquely determined by the Liouville 1-form, is pointing outwards along $\partial M$. The \textit{completion of a Liouville domain} $(M,d\lambda)$ denoted by $(\widehat{M},d\widehat{\lambda})$ is obtained by gluing the \textit{symplectization} of $\partial M$ to the Liouville domain $M$ in the following way. The negative part of symplectization of $\partial M$ is embedded into a tubular neighborhood of $\partial M$ in $M$ via the flow of the Liouville vector field
\begin{align*}
	( (0,1]\times \partial M , d(r \alpha) ) &\hookrightarrow (M,d\lambda)\\
	(r,x) & \mapsto \varphi^Z_{\ln r}(x) 	
\end{align*} 
while the positive part $([1,\infty) \times \partial M , d(r\alpha))$ is attached along the boundary $\partial M$ by identifying $\partial M$ with $\{1\}\times \partial M$, i.e. 
\[
\widehat{M}:=M\cup_{\partial M}([1,\infty)\times\partial M).
\]
The 1-form $\widehat{\lambda}$ is globally constructed in $\widehat{M}$ such that 
\[
\begin{cases}
\widehat{\lambda}|_{M}=\lambda\\
\widehat{\lambda}|_{[1,\infty) \times \partial M}=d(r\alpha)
\end{cases}
\]
and the extended Liouville vector field $\widehat{Z}$ on $[1,\infty)\times\partial M$ is given by $r\partial_{r}$. Then a Liouville manifold is nothing else than the completion of a Liouville domain. 

\begin{exmp}\label{Exa_Liou_mfd} The most common examples are the following.
\begin{itemize}
\item[(1)] $(\mathbb{R}^{2n},d\lambda_{\text{std}})$ has Liouville domains as $2n$-disk $(D^{2n},d\lambda_{\text{std}})$ with Liouville 1-form
\[
\lambda_{\text{std}}=\frac{1}{2}\sum_{i=1}^{n}(x_idy_i-y_idx_i),
\] 
and Liouville vector field
\[
Z=\frac{1}{2}\sum_{i=1}^{n}(x_i\partial_{x_i}+y_i\partial_{y_i}),
\]
where $(x_1,y_1,\dots,x_n,y_n)\in\mathbb{R}^{2n}$.
	
\item[(2)] $(T^*Q,d\lambda_{\text{can}})$ the cotangent bundle of a connected closed Riemannian $n$-manifold $Q$ has Liouville domain as codisk bundle $(D^*Q,d\lambda_{\text{can}})$ with Liouville 1-form
\[
\lambda_{\text{can}}=\sum_{i=1}^{n}p_i dq_i,
\]
and Liouville vector field
\[
Z=\sum_{i=1}^{n}p_i\partial_{p_i},
\]
where $(q_1,\dots,q_n,p_1,\dots,p_n)\in T^*Q$.
\end{itemize}
\end{exmp}

Lagrangian Floer theory for Liouville manifolds are introduced and established in \cite{Abouzaid-Seidel2010} with certain exact Lagrangian submanifolds as objects. One of the most important feature is that they provide a natural choice of $d\widehat{\lambda}$-compatible almost complex structure $J$, which is said to be of \textit{contact type} on $[r_0,\infty)\times\partial M$ for some $r_0\geq  1$: 
\[
\widehat{\lambda}\circ J=rd r
\]
where $r\in[1,\infty)$. This is equivalent to $J(r\partial_r)=R_\alpha$ where $R_\alpha$ is the Reeb vector field of the contact form $\alpha$ on $\partial M$, i.e. $i_{R_\alpha}\alpha=1$ and $i_{R_\alpha} d\alpha=0$. Then an $d\widehat{\lambda}$-compatible almost complex structure $J$ is \textit{admissible} if it is of contact type for $r_0\gg 1$. By maximum principle and similar type of argument, c.f. \cite[\S 7]{Abouzaid-Seidel2010}, there is $C^0$-bound of the Floer trajectories or $J$-holomorphic strips, which is required to obtain the Gromov--Floer compactness for the moduli spaces. See also \cite[\S 12.3]{Oh2015} for the treatment in cotangent bundles specially. Another important feature in Liouville manifold is that there is no non-constant $J$-holomorphic sphere due to exactness of the symplectic manifold, hence sphere bubbling is forbidden by geometric reason. 

The Floer theory for compact monotone Lagrangian submanifolds in a Liouville manifold instead of exact Lagrangian submanifolds can be established similarly in the geometric setting with convex type $d\widehat{\lambda}$-compatible almost complex structure. This is studied e.g. in $T^*S^n$ \cite{Abouzaid-Diogo2023}, in which Morse--Bott technique is also involved. For our purpose, we do not need a full Lagrangian Floer theory with higher algebraic structures, instead we only need to established the Lagrangian Floer homology for closed monotone Lagrangian submanifolds in Liouville manifold. Due to above discussions, there is a generic choice of admissible compatible almost complex structure such that all the statements for moduli spaces in Appendix \ref{App_monotone_Lag_Floer} and Appendix \ref{App_MBFloer} still hold in Liouville manifold, and the Lagrangian Floer homology is well-defined with the desired properties, e.g. Hamiltonian invariance.

\renewcommand{\thesection}{B} 
\section{Algebraic intersection numbers for clean intersections}\label{Appendix_B}
\begin{center}
\author{Urs Frauenfelder, Zhen Gao}
\end{center}

In this appendix, we introduce the notion of clean normal bundle over a clean intersection. The
algebraic intersection number can be computed directly in terms of the Euler number of the clean normal bundle. In the special case of Lagrangian intersections the clean normal bundle coincides with the tangent bundle of the intersection. 

\subsection{Introduction}\label{App_int_numb_intro}
Suppose throughout the appendix that $M$ is an oriented manifold and $L_1$ and $L_2$ are two closed connected oriented submanifolds of $M$ of complementary dimension, i.e.
$$\dim M=\dim L_1+\dim L_2.$$
If $L_1$ and $L_2$ intersect transversally in the sense that for each intersection point
$p \in L_1 \cap L_2$ it holds that
$$T_pM= T_p L_1 + T_p L_2$$
then one defines the algebraic intersection number 
$$L_1 \bullet L_2 \in \mathbb{Z}$$
as the signed count of intersection points of $L_1$ and $L_2$. The sign is determined as follows.
If $p \in L_1 \cap L_2$ choose positively ordered basis $[v_1,\ldots, v_{l_1}]$ of $T_p L_1$ and
$[w_1,\ldots, w_{l_2}]$ of $T_p L_2$, where we abbreviate by $l_1$ the dimension of $L_1$ and
by $l_2$ the dimension of $L_2$. If $[v_1,\ldots, v_{l_1}, w_1,\ldots, w_{l_2}]$ is a positively  oriented basis of $T_p M$, we put $\operatorname{sign}(p)=1$ and otherwise in the case where it is a negatively oriented basis we put $\operatorname{sign}(p)=-1$. Then we define
\begin{equation}\label{int}
	L_1 \bullet L_2=\sum_{p \in L_1 \cap L_2} \operatorname{sign}(p).
\end{equation}
The algebraic intersection number is a homotopy invariant. C.f. \cite[Chapter 5, Theorem 2.1]{Hirsch2012}. Therefore one can define it as well for non-transverse intersections. Namely, after small perturbation one can achieve that the two submanifolds intersect transversally so that one can define it for a transverse perturbation which by homotopy invariance does not depend on the choice of the perturbation. 

Although generically two submanifolds intersect transversally, in practise one quite often faces the situation where this is not the case. The reason for this is that we often consider the situtation where a Lie group $G$ acts on the manifold $M$. If the two submanifolds are invariant under the action of the Lie group, their intersection $L_1 \cap L_2$ is invariant as well and therefore intersection points are usually not isolated but come as whole 
$G$-orbits. For that reason one frequently has to deal with the case of \emph{clean intersections}, in the sense that the intersection $L_1\cap L_2$ is itself a disjoint union of connected submanifolds (possibly of various dimensions) in $M$ with the property that for each
$p \in L_1\cap L_2$ one has
\[
T_p (L_1\cap L_2)=T_p L_1 \cap T_p L_2.
\]
In this appendix, we want to discuss a formula which allows one to compute the algebraic intersection number directly in the case of a clean intersection without perturbing. For that purpose we introduce the \emph{clean normal bundle} (see Section \ref{App_clean_normal_bundle})
\[
\nu^\text{cl}_C \to C
\]
for each connected component $C$ in $L_1\cap L_2$. The rank of the vector bundle $\nu^{\text{cl}}_C$ satisfies
\begin{equation}\label{cleanrk}
\operatorname{rank}(\nu^\text{cl}_C)=\dim C.
\end{equation}
Moreover, total space of the clean normal bundle has a canonical orientation induced from the orientations of $L_1$, $L_2$, and $M$. The the total space of the clean normal bundle is canonically oriented as manifold, however, the clean normal bundle does not need to be orientable as a vector bundle, since we do not assume that the intersection $C$ is orientable.\footnote{A (real) vector bundle over a manifold $E\xrightarrow{\pi} B$ is orientable if and only if $w_1(E)\in H^1(B;\mathbb{Z}/2\mathbb{Z})$ vanishes. The total space $E$ is oriented as a manifold means that its tangent bundle $TE$ is an oriented vector bundle over the space $E$. The tangent bundle of $E$ has a non-canonical splitting $$TE\cong \pi^* TB\oplus\pi^*E.$$ By the Whitney product formula: 
\begin{align*}
w_1(TE)&=w_0(\pi^* TB)\smile w_1(\pi^*E)+w_1(\pi^* TB)\smile w_0(\pi^*E)\\
&= \pi^* w_1(E)+\pi^*w_1(TB).
\end{align*} 
If $B$ is orientable, i.e. $w_1(TB)=0$, then $w_1(TE)=\pi^*w_1(E)$, as $TE$ being an oriented vector bundle, we have $w_1(TE)=0=w_1(E)$, hence $E$ is an oriented vector bundle over $B$. If $B$ is non-orientable, then we have $0=w_1(E)+w_1(TB)$, this implies that $w_1(E)\equiv w_1(TB) \mod 2$, thus $E$ is a non-orientable vector bundle over $B$.} We can in view of (\ref{cleanrk}) define
the Euler characteristic 
\[
\chi(\nu^\text{cl}_C) \in \mathbb{Z}
\]
for every connected component $C\subset L_1\cap L_2$. For that, we do not need an orientation of $C$. See \S\ref{App_Euler} for more detail. We derive a formula for the algebraic intersection number in terms of this Euler characteristic
\begin{equation}\label{cleanint}
	L_1\bullet L_2=\sum_{C\in\mathcal{C}(L_1\cap L_2)}\chi(\nu^\text{cl}_C).
\end{equation}
In particular, when $L_1$ and $L_2$ are two Lagrangian submanifolds intersecting cleanly in a symplectic manifold $(M,\omega)$, we obtain:
\begin{equation}\label{lagcleanint}
	L_1\bullet L_2=\sum_{C\in\mathcal{C}(L_1\cap L_2)}\operatorname{sign}(C)\chi(C).
\end{equation}

\subsection{Clean normal bundles}\label{App_clean_normal_bundle}

We assume that $L_1$ and $L_2$ intersect cleanly, i.e. $L_1 \cap L_2$ is a disjoint union of connected submanifolds. We denote by $T_{C}L_1$ and $T_{C}L_2$ the restriction vector bundle of $TL_1$ and $TL_2$ to $C$, i.e. 
\[
T_{C}L_1:=TL_1|_{C}=\bigsqcup_{p\in C}T_pL_1,
\] 
and 
\[
T_{C}L_2:=TL_2|_{C}=\bigsqcup_{p\in C}T_pL_2;
\] 
one can also view them as pullback bundles of $TL_1$ and $TL_2$ along the inclusions $C\hookrightarrow L_1$ and $C\hookrightarrow L_2$. Over $C\in\mathcal{C}(L_1\cap L_2)$ we have the following three vector bundles
\[
\begin{tikzcd}
	T_{C}L_1 \arrow[rd] & TC \arrow[l, hook'] \arrow[r, hook] \arrow[d] & T_{C}L_2 \arrow[ld] \\
	& C                                             &                     
\end{tikzcd}
\]
We perform the fiber sum (distinguishing from Whitney sum) of the vector bundles $T_{C}L_1$ and $T_{C}L_2$ over $C$ to obtain a new vector bundle
\[
\begin{tikzcd}
	T_{C}L_1+T_{C}L_2 \arrow[d] \\
	C                                     
\end{tikzcd}
\]
whose fiber over a point $p \in C$ is given by the sum of vector spaces
\[
T_pL_1+T_pL_2 \subset T_p M.
\]
In particular, the vector bundle $	T_{C}L_1+T_{C}L_2$ is a subbundle of the vector bundle
$T_{C}M:=TM|_{C}$. We define the \emph{clean normal bundle} as the quotient bundle:
\begin{equation}
\nu^\text{cl}_{C}:=\frac{T_CM}{T_{C}L_1+T_{C}L_2 }.
\end{equation}
It is clear that $\nu^\text{cl}_C\subseteq\nu_C$ where $\nu_C=T_CM/TC$ is the ordinary normal bundle of $C$ in $M$. Abbreviating the dimensions as
\[
n=\dim C, \quad m=\dim M, \quad l_1=\dim L_1, \quad l_2=\dim L_2
\]
we have
\[
\operatorname{rank}(	T_{C}L_1+T_{C}L_2)=l_1+l_2-n
\]
and therefore using that $L_1$ and $L_2$ are of complementary dimension in $M$, so that
$l_1+l_2=m$, we obtain for the rank of the clean normal bundle over $C$:
\begin{align*}
	\operatorname{rank}(\nu^\text{cl}_{C})&=\operatorname{rank}(T_CM)-\operatorname{rank}(	T_{C}L_1+T_{C}L_2)\\
	&=m-(l_1+l_2-n)\\
	&=n.
\end{align*}
In order to avoid the necessity to discuss trivial cases, we have assumed that
\[
l_1>0, \quad l_2>0.
\]
\subsubsection*{Orientation}\label{ori_clean_norm_bun} We next explain how the orientations of $L_1$, $L_2$ and $M$ induce an orientation on the total space of $\nu^\text{cl}_{C}$. For that purpose we assume that $p \in C$. We consider several cases:
\begin{description}
	\item[Case\,1:] \emph{Assume that $0<n<\min\{l_1,l_2\}$.}
	
	\noindent In other words, $n$ is positive and smaller than $l_1$ as well as $l_2$.
	Choose an ordered basis of $T_p C$:
	\[
	\mathfrak{B}(T_pC)=[v_1,\ldots,v_{n}].
	\]
	Extend this basis to a positively oriented basis of $T_p L_1$:
	\[
	\mathfrak{B}(T_pL_1)=[v_1,\ldots, v_{n},w_1,\ldots w_{l_1-n}],
	\]
	and to a positively oriented basis of $T_p L_2$:
	\[
	\mathfrak{B}(T_pL_2)=[u_1,\ldots u_{l_2-n},v_1,\ldots, v_{n}].
	\]
	In contrast to $T_pL_1$, where the additional basis vectors are added after the basis vectors of $T_pC$, the basis vectors for $T_pL_2$ are added before the ones of $T_p C$.
	Now extend the union of these two bases to a positively oriented basis of $T_p M$:
	\[
	\mathfrak{B}(T_pM)=[ u_1,\ldots, u_{l_2-n},v_1,\ldots, v_{n},w_1,\ldots w_{l_1-n},\eta_1,\ldots ,\eta_{n}].
	\]
	Abbreviate by
	\[
	\pi : T_p M \to \nu^\text{cl}_{C}|_{p}
	\]
	the projection map of $T_p M$ onto the space $\nu^\text{cl}_{C}|_{p}$. Then 
	\[
	\pi:(u_1,\dots, u_{l_2-n},v_1,\dots, v_{n},w_1,\dots w_{l_1-n},\eta_1,\dots ,\eta_{n})\mapsto(\eta_1,\dots ,\eta_{n})
	\]
	and 
	\[
	\mathfrak{B}(T_{(p,0)}\nu^\text{cl}_{C})=[v_1,\ldots, v_{n}, \eta_1,\ldots, \eta_{n}]
	\]
	is an ordered basis of $T_{(p,0)} \nu^\text{cl}_{C}\cong T_pC\oplus \nu^\text{cl}_{C}|_{p}$. We define an orientation on the total space of $\nu^\text{cl}_{C}$ by declaring this basis to be positively oriented. This does not depend on the choices involved. Indeed, let us change the orientation of $T_p C$ by replacing $v_i$ by $-v_i$ for any $i\in\{1,\dots,n\}$. In order to meet the requirement that $\mathfrak{B}(T_pL_1)$ and $\mathfrak{B}(T_pL_2)$ are ordered bases with respect to the given orientations on $T_p L_1$ respectively $T_pL_2$ we replace $w_i$ by $-w_i$ respectively $u_i$ by
	$-u_i$, so that two of the basis vectors change sign. In order that $\mathfrak{B}(T_pM)$ is still an ordered basis of $T_p M$ we replace in addition $\eta_i$ by $-\eta_i$ so that four of the basis vectors change sign. Now in the basis $\mathfrak{B}(T_{(p,0)}\nu^\text{cl}_{C})$, the vectors $v_i$ and $\eta_i$ change sign, so that
	the orientation of $T_{(p,0)}\nu^\text{cl}_{C}$ does not change.
	
	\item[Case\,2:] \emph{Assume that $n=l_1<l_2$.} 
	
	\noindent In this case the orientation of $L_1$ induces an orientation of $C$, so that we start with a positively oriented basis of $T_p L_1=T_p C$. We extend this basis to a positively oriented basis of $T_p L_2$ as in Case\,1 by adding the additional basis vector before the ones of
	$T_pL_1=T_p C$. We extend this basis further to a positively oriented basis of $T_p M$. As in Case\,1 this induces an orientation on $T_{(p,0)} \nu^\text{cl}_{C}$.
	
	\item[Case\,3:] \emph{Assume that $n=l_2<l_1$.}
	
	\noindent This case is analogous to Case\,2 just with the roles of $L_1$ and $L_2$ interchanged. 
	
	\item[Case\,4:] \emph{Assume that $n=l_1=l_2$.} 
	
	\noindent In this case $\nu^\text{cl}_{C}$ coincides with the usual normal bundle $\nu_C$ of $C$ in $M$. An orientation of $C$ together with the orientation of $M$ induces an orientation on the normal bundle which together with the orientation on $C$ gives an orientation on the total space of the normal bundle. This orientation of the total space is independent of the choice of the orientation of $C$ since flipping the orientation of $C$ also flips the orientation of the normal bundle so that the orientation on the total space does not change. If the orientations of $L_1$ and $L_2$ along $C$ agree, endow the total space of the normal bundle with its canonical orientation. Otherwise in the case where the orientations of $L_1$ and $L_2$ along $C$ are opposite, endow the total space with the orientation which is opposite to the canonical one. 
	
	\item[Case\,5:] \emph{Assume that $n=0$.}
	
	\noindent In this case $C=\{p\}$ is a single point at which the manifolds $L_0$ and $L_1$ intersect transversely. The clean normal bundle $\nu^\text{cl}_{C}$ as well just consists of the point $\{(p,0)\}$. An orientation of a point is just the assignment of a sign to the point. We endow it with the sign
	$\operatorname{sign}(p)$ as discussed before (\ref{int}). This finishes the discussion  of the orientation of the total space of $\nu^\text{cl}_{\{p\}}$.
\end{description}

Suppose now that $(M,\omega)$ is a symplectic manifold and $L_1$ and $L_2$ are oriented Lagrangian submanifolds of $M$ intersecting cleanly. Because the symplectic form is non-degenerate, the normal bundle of any submanifold can be identified with its symplectic orthogonal complement bundle:
\begin{equation}\label{normal_symp_orth}
\nu L_1\cong TL_1^{\perp\omega}, \quad \nu L_2\cong TL_2^{\perp\omega}.
\end{equation}
Moreover, since $L_1$ and $L_2$ are Lagrangian, we have in particular
\begin{equation}\label{symp_orth_Lag}
TL_1=TL_1^{\perp\omega}, \quad TL_2=TL_2^{\perp\omega}.
\end{equation}
\begin{lem}\label{normal}
Let $L_1$ and $L_2$ be Largangian submanifolds intersecting cleanly. There is a vector bundle isomorphism:
\[
\nu^{\normalfont\text{cl}}_{C}\cong TC
\]
for each $C\in\mathcal{C}(L_1\cap L_2)$.
\end{lem}

\begin{proof}
We only need to verify everything point-wise at each $p\in C$. By the identifications \eqref{normal_symp_orth}, we have
\[
\nu_pL_1\cong T_pL_1^{\perp\omega},\qquad \nu_pL_2\cong T_pL_2^{\perp\omega}.
\]
By the definition $\nu_pL_i:=T_pM/T_pL_i$, $i=1,2$, there is an identification
\begin{align*}
\nu_pL_1\cap \nu_pL_2&=\left(\frac{T_pM}{T_pL_1}\right)\cap\left(\frac{T_pM}{T_pL_2}\right)\\
&=\frac{T_pM}{T_pL_1+T_pL_2}\\
&=:\nu_C^\text{cl}|_p.
\end{align*}
Thus, we obtain the desired isomorphism in the statement
\begin{align*}
\nu_C^\text{cl}|_p&=\nu_pL_1\cap \nu_pL_2\\
&\cong T_pL_1^{\perp\omega} \cap T_pL_2^{\perp\omega}\\
&\cong T_pL_1 \cap T_pL_2\\
&=T_pC
\end{align*}
where we have applied \eqref{normal_symp_orth} and \eqref{symp_orth_Lag} to get the second and the third isomorphism respectively. 

In fact, we can also show that there is a vector bundle isomorphism 
\begin{equation}\label{nu---TL_1+TL_2_perp_omega}
\nu_C^{\text{cl}}\cong (T_CL_1+T_CL_2)^{\perp\omega}.
\end{equation} 
Recall that for a symplectic vector space $(V,\omega)$ and subspace $W\subset V$, we have 
\[
\dim V=\dim W+\dim W^{\perp\omega}.
\]
In particular, 
\[
V=W\oplus W^{\perp\omega},
\]
i.e. 
\[
W\cap W^{\perp\omega}=\{0\}
\]
holds if and only if $W$ is a symplectic subspace. Note that this is not a necessary condition for above dimension formula. However, $T_pL_1+T_pL_2\subset T_pM$ is clearly not a symplectic subspace in general, namely 
\[
(T_pL_1+T_pL_2)\cap(T_pL_1+T_pL_2)^{\perp\omega}\neq\{0\}.
\] 
The way we show \eqref{nu---TL_1+TL_2_perp_omega} is by making use of the canonical tangent-cotangent isormophisms via $\omega$: 
\[
T_pM\cong T^*_pM
\]
such that 
\[
T_p^*M=\{ \omega(X,\cdot) : X\in T_pM   \};
\]
and via the (non-natural) isomorphisms of vector spaces and their algebraic duals
\[
T_pL_1\cong T_p^*L_1, \qquad T_pL_2\cong T_p^*L_2;
\]  
we thus have an isormophism 
\begin{align*}
\nu_C^\text{cl}|_p:=&\frac{T_pM}{T_pL_1+T_pL_2}\\
\cong&\frac{T_p^*M}{T_p^*L_1+T_p^*L_2}.
\end{align*}
At each $p\in C$, 
\[
T_p^*L_1+T^*_pL_2=\{ \alpha+\beta : \alpha\in T_p^*L_1, \;\beta\in T_p^*L_2  \}  ,
\]
where the covectors can be expressed as the linear cobinations
\begin{align*}
\alpha=\sum_{i=1}^{n}a_iv_i^*+\sum_{i=1}^{l_1-n}\tilde{a}_iw_i^* ,\\
\beta=\sum_{i=1}^{l_2-n}\tilde{b}_iu_i^*+\sum_{i=1}^{n}b_iv_i^* ,\\
\end{align*}
in terms of the dual basis 
\[	
\mathfrak{B}(T_p^*L_1)=[v_1^*,\ldots, v_{n}^*,w_1^*,\ldots w_{l_1-n}^*],
\]
and
\[
\mathfrak{B}(T_p^*L_2)=[u_1^*,\ldots u_{l_2-n}^*,v_1^*,\ldots, v_{n}^*],
\]
with $a_i,\tilde{a}_i,b_i,\tilde{b}_i\in\mathbb{R}$ for all $i$. Therefore, we can see explicitly that
\[
\frac{T^*_pM}{T^*_pL_1+T^*_pL_2}=\{ \omega(X,\cdot): X\in T_pM,\; \omega(X,Y)=0\; \forall Y\in T_pL_1+T_pL_2 \}.
\]
Applying the cotangent-tangent isomorphism via $\omega^{-1}$, we have
\begin{align*}
\frac{T^*_pM}{T^*_pL_1+T^*_pL_2}&\cong\{ X\in T_pM: \omega(X,Y)=0 \;\forall Y\in T_pL_1+T_pL_2 \}\\
&=(T_pL_1+T_pL_2)^{\perp\omega}.
\end{align*}
Hence, we have verified \eqref{nu---TL_1+TL_2_perp_omega}, and
\begin{align*}
\nu_C^{\text{cl}}|_p&\cong(T_pL_1+T_pL_2)^{\perp\omega}\\
&=(T_pL_1)^{\perp\omega}\cap(T_pL_2)^{\perp\omega}\\
&=T_pL_1 \cap T_pL_2\\
&=T_pC.
\end{align*}
This completes the proof.
\end{proof}

\begin{rem}
Alternatively, Lemma \ref{normal} can be also proved by using a compatible almost complex structure. By definition, at each $p\in C$, then we have the splitting:
\[
T_pM=(T_pL_1+T_pL_2)\oplus\nu^\text{cl}_C|_{p}.
\]
On the other hand, let $J$ be a $\omega$-compatible almost complex structure, we have
\[
T_pL_1\cap JT_pL_1=\{ 0\}, \qquad T_pL_2\cap JT_pL_2=\{0\}.
\]
This yields that 
\[
(T_pL_1+T_pL_2)\cap JT_pC=\{0\}.
\] 
Thus, by dimension reasons, we have
\[
T_pM=(T_pL_1+T_pL_2)\oplus JT_pC.
\]
Therefore, we have the identification of vector bundles
\[
\nu^\text{cl}_{C}=JTC.
\]
We conclude that there is an isomorphism of vector bundles
\[
\nu^\text{cl}_C\cong TC
\] 
given by the almost complex structure $J$.
\end{rem}

\subsection{Euler characteristics and algebraic intersection numbers}\label{App_Euler}

\subsubsection*{Euler characteristics of vector bundles} Let $E \to B$ be a real vector bundle (not necessarily orientable) over a closed connected manifold $B$ satisfying
\[
\operatorname{rank}(E)=\dim B=:n.
\]
Assume in addition that the total space $E$ is oriented as manifold. We then define the Euler characteristic $\chi(E) \in \mathbb{Z}$ as follows: If $n=0$, the total space $E$ is just a single point and an orientation is the assignment of a sign to the point. In this case the Euler characteristic is just this sign. Suppose now that $n>0$. Choose a smooth section 
\[
\sigma : B \to E
\]
transverse to the zero section, i.e. $\sigma\pitchfork B$. This is in fact a generic property for $\sigma\in\Gamma(E)$ and the zero locus $\sigma^{-1}(0)$ is a discrete finite set given $B$ compact. If $p \in B$ is a zero of $\sigma$, define the sign of $p$, $\operatorname{sign}(p)$, as follows. Choose an ordered basis of $T_p B$
\[
\mathfrak{B}(T_pB)=[v_1,\ldots, v_n].
\]
Since $s$ is transverse to the zero section, the ordered set
\[
\mathfrak{B}(T_pE)=[v_1,\ldots,v_n, D\sigma_p(v_1),\ldots,D\sigma_p(v_n)]
\]
is a basis of $T_p E$. If this ordered basis is positively oriented with respect to the given orientation on $E$ we define $\operatorname{sign}(p)=+1$ and otherwise $\operatorname{sign}(p)=-1$. This definition does not depend on the choice of the basis $\mathfrak{B}(T_pB)$. Indeed, pick any $i\in\{1,\dots,n\}$, if $v_i$ changes sign, then as well
$D\sigma_p(v_i)$ changes sign, so that we have two sign changes in $\mathfrak{B}(E_p)$. In particular, we need only the orientation of the total space $E$ to define its Euler characteristic but no orientation of the base $B$. The Euler characteristic of $E$ is given by the formula
\begin{equation}
\chi(E):=\sum_{p\in \sigma^{-1}(0)}\operatorname{sign}(p)=B\bullet S
\end{equation}
where $S:=\operatorname{im}(\sigma)$. Due to homotopy invariance of the algebraic intersection number, this quantity is a well-defined homotopy invariant and independent of the choice of transversal section. This is a mild generalization of \cite[Chapter 5, \S 2]{Hirsch2012}. The rank $n$ vector bundle over a smooth $n$-manifold is only required to have oriented total space instead of being an oriented vector bundle. In particular, when $E\to B$ is an oriented (real) vector bundle, we have the following expression in terms of characteristic class 
\[ 
\chi(E)=\langle e(E),[B]\rangle 
\] 
where $e(E)\in H^n(B;\mathbb{Z})$ is the Euler class of $E$ and $[B]\in H_n(B;\mathbb{Z})$ is the fundamental class of $B$ given it is compact and oriented. C.f. \cite[Chapter VI, \S 12]{Bredon1993}.

\subsubsection*{Algebraic intersection number}
We are now in position to write down a formula for the algebraic intersection number of 
a clean intersection. Abbreviate by $\mathcal{C}(L_1\cap L_2)$ the set of connected components of
$L_1\cap L_2$. Then the formula is
\begin{equation*}
	L_1 \bullet L_2=\sum_{C \in \mathcal{C}(L_1 \cap L_2)} \chi(\nu^\text{cl}_C).\tag{\ref{cleanint}}
\end{equation*}
Using a transverse section of the clean normal bundle, one can perturb $L_2$ locally around $C$ to make the intersection transverse so that the signed count of intersection points of $L_1$ with the perturbed $L_2$ corresponds to the signed count of zeros of this transverse section which defines the Euler characteristic of the clean normal bundle, referring to the discussion of orientation of the clean normal bundle determined by that of $L_1,L_2,M$ in \S\ref{ori_clean_norm_bun}. Therefore, one sees that the right hand side of the formula above indeed computes the algebraic intersection number. 

\begin{exmp}
We first confirm that this formula coincides with \eqref{int} in the special case of a transversal intersection. For a transversal intersection we have
\[
\mathcal{C}(L_1 \cap L_2)=\{\{p\}:p\in L_1 \cap L_2\}
\]
and for every $p \in L_1 \cap L_2$ it holds that
\[
\nu^\text{cl}_{\{p\}}=\{p\}\times\{0\}
\]
Thus, the clean normal bundle is just a signed point for which one gets
\[
\chi(\nu^\text{cl}_{\{p\}})=\operatorname{sign}(p).
\]
\end{exmp}

\begin{exmp}
Another well-known case of formula (\ref{cleanint}) is the case $L_1=L_2=L$ giving rise to the algebraic self-intersection number. In this case the clean normal bundle is just the usual normal bundle endowed with its canonical orientation if the orientations of $L_1$ and $L_2$ agree and otherwise endowed with the opposite orientation. In this case (\ref{cleanint}) simplifies to
\[
L_1 \bullet L_2=L\bullet L=\chi(\nu^\text{cl}_{L})=\chi(\nu_L).
\]
In particular, when $L_1=L_2=\Delta_M\subset M\times M$, we have
\[
\chi(TM)=\Delta_M \bullet \Delta_M=\chi(\nu^\text{cl}_{\Delta_M})=\chi(\nu_{\Delta_M}).
\] 
In this case, one can see that actually 
\[
\nu^\text{cl}_{\Delta_M}=\nu_{\Delta_M}=\bigsqcup_{(p,p)\in M\times M}\frac{T_pM\times T_pM}{T_{(p,p)}\Delta_M}\cong\bigsqcup_{p\in M}T_pM= TM
\]
by referring to section \ref{App_clean_normal_bundle} for the definition of clean normal bundle in this case.
\end{exmp}

We finally discuss the case of Lagrangian intersections. We assume that $(M,\omega)$ is a symplectic manifold and $L_1$ and $L_2$ are oriented Lagrangian submanifolds of $M$ intersecting cleanly. As a symplectic manifold, $M$ is endowed with a canonical orientation determined by the volume form $\frac{1}{n!}\omega^{n}$. In this case, the clean normal bundle coincides with the tangent bundle of $C$. This correspondence is not necessarily an isomorphism of oriented manifolds. Note that the tangent bundle of any manifold $C$ carries a canonical orientation as manifold. If $C$ is closed, the Euler characteristic
of the canonically oriented tangent bundle coincides with the Euler characteristic of the manifold, i.e. $\chi(TC)=\chi(C).$

Assume now all connected components of the intersection of $L_1$ and
$L_2$ are closed. For $C\in \mathcal{C}(L_1 \cap L_2)$, we associate to it a sign
$$\operatorname{sign}(C) \in \{\pm 1\}$$
as follows. If the orientations of $\nu^\text{cl}_{C}$ agrees with the canonical orientation of
$TC$ under the isomorphism in Lemma \ref{normal}, we put $\operatorname{sign}(C)=+1$, and otherwise $\operatorname{sign}(C)=-1$. Thus, we have
\begin{equation}
\chi(\nu^\text{cl}_C)=\operatorname{sign}(C)\chi(TC).
\end{equation}
With these conventions, \eqref{cleanint} simplifies in the case of a clean Lagrangian intersection to
\begin{equation*}
L_1 \bullet L_2=\sum_{C \in \mathcal{C}(L_1 \cap L_2)} \operatorname{sign}(C)\chi(C). \tag{\ref{lagcleanint}}
\end{equation*}

\subsection{Categorification}\label{App_Cat}
We give a direct proof of the formula
\begin{equation}\label{categorification}
\chi(HF_*(L_1,L_2,\{f_i\};\Lambda))=L_1\bullet L_2
\end{equation}
utilizing \eqref{lagcleanint} whenever the Lagrangian Floer homology $HF_*(L_1,L_2,\{f_i\};\Lambda)$ is well-defined for the pair of oriented Lagrangian submanifolds $L_1$ and $L_2$ intersecting cleanly whose chain complex $CF_*(L_1,L_2,\{f_i\};\Lambda)$ is $\mathbb{Z}$-graded via the Maslov index and finitely generated over $\Lambda$ in each degree. 

Suppose
\[
L_1\cap L_2=\bigsqcup_{i\in I} C_{i}
\]
where $I$ is a finite set, and each $C_i$ is connected oriented closed submanifold. The Floer complex in the Morse--Bott setting is generated by the critical points of the auxiliary Morse functions $f_i$ defined on each $C_i$:
\[
CF_*(L_1,L_2,\{f_i\};\Lambda):=\bigoplus_{k\in\mathbb{Z}}\bigoplus_{\substack{ x\in \underset{i\in I}{\bigsqcup} \operatorname{Crit}(f_i)  \\  |x|=k } }\Lambda \langle x \rangle
\]
where we recall that the grading is given by \eqref{MB_grading} and \eqref{MB_index}. The Euler characteristic of the homology, due to finiteness assumption, is computed by
\begin{align*}
\chi(HF_*(L_1,L_2,\{f_i\};\Lambda))&:= \sum_{k\in\mathbb{Z}}(-1)^k \operatorname{rank}_{\Lambda} HF_k(L_1,L_2,\{f_i\};\Lambda)\\ 
&= \sum_{k\in\mathbb{Z}}(-1)^k \operatorname{rank}_{\Lambda} CF_k(L_1,L_2,\{f_i\};\Lambda) \\
&=\sum_{k\in\mathbb{Z}}\sum_{\substack{ x\in \underset{i\in I}{\bigsqcup}\operatorname{Crit}(f_i)  \\  |x|=k } }(-1)^{|x|} \\
&=\sum_{i\in I}  (-1)^{\mu(C_i)} \sum_{x\in\operatorname{Crit}(f_i)} (-1)^{\operatorname{ind}_{f_i}(x)}
\end{align*}
It is well-known that
\begin{equation}\label{Morse_Euler}
\sum_{x\in\operatorname{Crit}(f_i)} (-1)^{\operatorname{ind}_{f_i}(x)}=\chi(C_i).
\end{equation}
The parity of the Maslov--Morse index $(-1)^{\mu(C_i)}$ is equal to the orientation sign of $C_i$ up to a global sign, i.e.
\begin{equation}\label{Maslov_sign}
(-1)^{\mu(C_i)}=\pm\operatorname{sign}(C_i).
\end{equation}
This can be seen as follows. According to \cite[\S 2.2.2]{Fukaya-Oh-Ohta-Ono2009}, when the Lagrangian submanifolds $L_1$ and $L_2$ are oriented, the definition of the Maslov--Morse index involves a choice of section of the oriented Lagrangian Grassmannian bundle fitting in the orientations of $L_1$ and $L_2$, therefore the Maslov--Morse index at a intersection point in $L_1\cap L_2$ is the orientation sign of the intersection point up to a global sign depending on the convention, this is also discussed by Paul Seidel, see \cite[\S 2.d, (vi)]{Seidel2000}. Further due to \cite[\S3.7 (3.7.49)]{Fukaya-Oh-Ohta-Ono2009}, the parity of the Maslov--Morse index for every point in Morse--Bott component $C_i$ is the same, thus we obtain the desired formula \eqref{Maslov_sign}. The formulas \eqref{Morse_Euler} and \eqref{Maslov_sign} together imply that
\[
\chi(HF_*(L_1,L_2,\{f_i\};\Lambda))=\pm\sum_{i\in I}\operatorname{sign}(C_i)\chi(C_i)
\]
which is \eqref{categorification} by \eqref{lagcleanint}.

\end{document}